\documentclass[12pt,reqno,twoside]{amsart}
\usepackage{amsmath}
\usepackage{amssymb}
\usepackage{bbm}
\usepackage{mathrsfs}
\usepackage{epsfig}
\usepackage{color}
\usepackage{pdfsync}
\usepackage{tensor}
\usepackage{mathabx}
\usepackage{tikz-cd}
\usepackage[english]{babel}
\usepackage[unicode=true,pdfusetitle,
 bookmarks=true,bookmarksnumbered=false,bookmarksopen=false,
 breaklinks=false,pdfborder={0 0 1},backref=false,colorlinks=false]
 {hyperref}

\title{Classification and statistics of cut-and-project sets}  
\author{Ren\'e R\"uhr}
\address{Weizmann Institute
{\tt rene.ruhr@weizmann.ac.il}}

\author{Yotam Smilansky}
\address{Rutgers University
{\tt yotam.smilansky@rutgers.edu}}
\author{Barak Weiss}
\address{Tel Aviv University 
{\tt barakw@tauex.tau.ac.il}}

\font\sn = cmssi8 scaled \magstep0
\font\si = cmti8 scaled \magstep0
\long\def\comrene#1{\ifdraft{{\color{blue}\si #1 }}\else\ignorespaces\fi}
\long\def\combarak#1{\ifdraft{\color{red}\sn #1 }\else\ignorespaces\fi}

\numberwithin{equation}{section}

\newcommand{\Liew}[1]{\lowercase{\mathfrak{#1}}}

\DeclareMathOperator{\Sp}{Sp}

\newcommand{\CC}{\mathbb{C}}
\newcommand{\RR}{\mathbb{R}}
\newcommand{\QQ}{\mathbb{Q}}

\newcommand{\sset}{\mathfrak{S}}

\newcommand{\Var}{\mathrm{Var}}

\newcommand{\covol}{\mathrm{covol}}

\newif\ifdraft\drafttrue
\draftfalse
\newcommand\name[1]{\label{#1}{\ifdraft{\sn [#1]}\else\ignorespaces\fi}}

\newcommand\eq[2]{{\ifdraft{\ \tt [#1]}\else\ignorespaces\fi}\begin{equation}\label{#1}{#2}\end{equation}}
\newcommand {\equ}[1]{\eqref{#1}}

\newcommand{\Q}{{\mathbb {Q}}}

\newcommand{\R}{{\mathbb{R}}}

\newcommand{\Res}{{\mathrm{Res}}}
\newcommand{\Z}{{\mathbb{Z}}}

\newcommand{\on}{\mathrm}

\newcommand{\bR}{\mathbb{R}}
\newcommand{\bQ}{\mathbb{Q}}

\newcommand{\LLN}{{\mathscr{X}}}
\newcommand{\ALN}{{\mathscr{Y}}}
\newcommand{\Cl}{{\mathscr{C}}}
\newcommand{\Prob}{{\mathrm{Prob}}}

\newcommand{\Propto}{\, \propto \,}

\newcommand{\N}{{\mathbb{N}}}

\newcommand{\Symp}{{\operatorname{Sp}}}

\newcommand{\Ad}{{\operatorname{Ad}}}
\newcommand{\ad}{{\operatorname{ad}}}

\newcommand{\GL}{\operatorname{GL}}

\newcommand{\SL}{\operatorname{SL}}

\newcommand{\SO}{\operatorname{SO}}

\newcommand{\diag}{{\rm diag}}

\newcommand{\Gal}{{\rm Gal}}
\newcommand {\ignore}[1]  {}

\newcommand{\spa}{{\rm span}}
\newcommand{\rank}{{\mathrm{rank}}}

\newcommand{\E}{{\mathfrak{E}}}
\newcommand{\ee}{{\mathfrak{e}}}

\newcommand{\ASL}{{\operatorname{ASL}}}

\newcommand{\PP}{{\mathcal P}}

\newcommand{\OO}{{\mathcal O}}

\newcommand{\LL}{{\mathcal L}}
\newcommand{\cL}{{\LL}}

\newcommand{\df}{{\, \stackrel{\mathrm{def}}{=}\, }}

\newcommand{\KK}{{\mathbb{K}}}

\newcommand{\hs}{\kern 0.8pt}

\newcommand{\til}{\widetilde}

\newcommand{\supp}{{\rm supp}}

\newcommand{\fphat}{\tensor[^p]{\widehat{f}}{}}

\newcommand{\Mat}{{\operatorname{Mat}}}

\newcommand{\sm}{\smallsetminus}

\newcommand{\vre}{\varepsilon}

\newcommand{\Vphys}{{V_{\mathrm{phys}}}}
\newcommand{\Vphyss}{{V'_{\mathrm{phys}}}}
\newcommand{\Vint}{{V_{\mathrm{int}}}}
\newcommand{\Vintt}{{V'_{\mathrm{int}}}}
\newcommand{\piphys}{{\pi_{\mathrm{phys}}}}
\newcommand{\piint}{{\pi_{\mathrm{int}}}}
\newcommand{\Leb}{{\mathrm{vol}}}
\newcommand{\vol}{{\mathrm{vol}}}
\newcommand{\OK}{{\mathcal{O}_{\KK}}}

\newtheorem{thm}{Theorem}[section]

\newtheorem{lem}[thm]{Lemma}

\newtheorem{prop}[thm]{Proposition}

\newtheorem{cor}[thm]{Corollary}

\newtheorem{remark}[thm]{Remark}

\newtheorem{example}[thm]{Example}

\newtheorem{dfn}[thm]{Definition}

\begin{document}
\date{\today}
\maketitle

\begin{abstract}
We define {\em Ratner-Marklof-Str\"ombergsson
measures} (following \cite{MS}). These are probability measures supported on
cut-and-project sets in $\R^d \ (d \geq 2)$ which are invariant and
ergodic for the action of 
the groups $\ASL_d(\R)$ or $\SL_d(\R)$. We classify the measures that
can arise in terms of algebraic groups and homogeneous dynamics. Using the
classification, we prove analogues of results of Siegel, Weil and
Rogers about a Siegel summation formula and identities and bounds
involving higher moments. We deduce results about
asymptotics, with error estimates, of point-counting and patch-counting for typical
cut-and-project sets. 
\end{abstract}

\section{Introduction}
A {\em cut-and-project set} is a discrete subset of $\R^d$
obtained by the following construction. Fix a direct sum decomposition
$\R^n = \R^d \oplus \R^{m}$, where the two summands in this
decomposition are denoted respectively $\Vphys, \Vint,$ so that
$$\R^n = \Vphys \oplus \Vint, $$
and the corresponding projections are 
$$\piphys: \R^n
\to \Vphys, \ \ \ \ \ \piint: \R^n \to \Vint. $$
Also fix a lattice $\LL \subset \R^n$ and a {\em window} $W \subset \Vint$;
then the corresponding cut-and-project set $\Lambda = \Lambda(\LL,W)$ is
given by
\eq{eq: def c and p set}{
\Lambda(\LL,W) \df \piphys\left(\LL \cap \pi^{-1}_{\mathrm{int}}(W) \right).
}
We sometimes allow $\LL$ to be a {\em grid},
i.e., the image of a lattice 
under a translation in $\R^n$, and sometimes require $\Lambda$ to be
{\em irreducible}, a notion we define in \S \ref{sec: 
  basics}.  
Cut-and-project sets are prototypical
aperiodic sets exhibiting 
long-term-order, and are sometimes referred to as {\em model sets} or {\em
  quasicrystals}.
Beginning with work of Meyer \cite{Meyer_book_french} in connection to
Pisot numbers, they have been intensively studied from various points
of view. See \cite{BaakeGrimm} and the references therein.

Given a
cut-and-project set, a natural
operation is to take the closure (with respect to a natural topology)
of its orbit under translations. This yields a dynamical system
for the translation group and has been studied by many
authors under different names. In recent years several investigators
have become interested in the orbit-closures under the
group $\SL_d(\R)$ (respectively $\ASL_d(\R)$), which is the
group of orientation- and volume-preserving linear (resp., affine) 
transformations of $\R^d$. In particular, in the important paper
\cite{MS}, 
motivated by problems in mathematical physics, Marklof and Str\"ombergsson 
introduced a class of natural probability measures on these
orbit-closures. The goal of this
paper is to classify and analyze such measures, and derive
consequences for the statistics and large 
scale geometry of cut-and-project sets. 

\subsection{Classification of Ratner-Marklof-Str\"ombergsson measures}
We say that a cut-and-project set is {\em irreducible} if it arises
from the above construction, where the data satisfies 
the assumptions {\bf (D), (I)} and {\bf (Reg)} given in \S \ref{subsec: c
  and p sets}. Informally speaking, {\bf (D)} and {\bf (I)} imply that
the set cannot be presented as a finite union of sets whose
construction involves smaller groups in the
cut-and-project construction, and {\bf (Reg)} is a regularity  
assumption on the window set $W$. We denote by $\Cl(\R^d)$ the space
of closed subsets of $\R^d$, equipped with the Chabauty-Fell topology.
This is a compact metric topology whose definition is recalled in \S
\ref{subsec: CF}, and which is also
referred to in the quasicrystals literature as the {\em
	local rubber  topology} or the  {\em natural topology}. Since
the groups $\ASL_d(\R)$  and $\SL_d(\R)$
act on $\R^d$, 
they also 
act on $\Cl(\R^d)$. We say that a Borel probability measure
$\mu$ on  $\Cl(\R^d)$ is a {\em
  Ratner-Marklof-Str\"ombergsson measure}, or {\em RMS measure} for
short, if it is invariant and ergodic under $\SL_d(\R)$ and gives full
measure to the set of irreducible cut-and-project sets.
We call it
{\em affine} if it is also invariant under $\ASL_d(\R)$, and {\em
  linear}
otherwise (i.e., 
if it is invariant under $\SL_d(\R)$ but not under $\ASL_d(\R)$). 

A construction of RMS measures was given in \cite{MS}, as
follows. Let $\ALN_n$ denote the space of grids of covolume one in
$\R^n$, equipped with the Chabauty-Fell topology, or equivalently with
the 
topology it 
inherits from its identification with the homogeneous space
$\ASL_n(\R)/\ASL_n(\Z)$.  Similarly, let $\LLN_n$ denote the space of
lattices of covolume one in $\R^n$, which is identified with the
homogeneous space 
$\SL_n(\R)/\SL_n(\Z)$. Fix the data $d,m, \Vphys \simeq \R^d, 
\Vint \simeq \R^m, \piphys, \piint$, as well as a set $W \subset \Vint$, and
choose $\LL$ randomly according to a probability measure $\bar \mu$ on 
$\ALN_n$. This data determines a cut-and-project set $\Lambda$, which
is random since $\cL$ is. The resulting probability measure $\mu$ on
cut-and-project sets can thus be written as the pushforward of $\bar
\mu$ under the map $\cL \mapsto \Lambda(W, \cL)$, and is easily seen to be
invariant and ergodic under 
$\SL_d(\R)$ or $\ASL_d(\R)$ if the same is true for $\bar \mu$. One
natural choice for $\bar \mu$ is
the so-called {\em Haar-Siegel measure,} which is the unique Borel
probability measure invariant under the 
group $\ASL_n(\R)$. Another is the Haar-Siegel measure on $\LLN_n$
(i.e., the unique $\SL_n(\R)$-invariant measure). It is also possible to consider other measures on
$\ALN_n$ which are $\ASL_d(\R)$- or $\SL_d(\R)$-invariant. As observed in
\cite{MS}, a fundamental result of Ratner
\cite{Ratner-annals} makes it possible to give a precise description
of such measures on $\ALN_n$. They correspond to certain algebraic
groups which are 
subgroups of $\ASL_n(\R)$ and contain $\ASL_d(\R)$ (or
$\SL_d(\R)$).

Our first result is a classification of such measures.
We refer to \S \ref{sec: basics} and \S  
\ref{sec: classification} for more precise statements, and for
definitions of the terminology.
\begin{thm}\name{thm: classification}
  Let $\mu$ be an RMS measure on $\Cl(\R^d)$. Then, up to rescaling, there are fixed
  $m$ and $W \subset \R^m$ such that $\mu$ is the pushforward via the
  map
  $$
\ALN_n \to \Cl(\R^d), \ \ 
\cL \mapsto \Lambda(\cL, W)$$
of a measure $\bar \mu$ on $\ALN_n$, 
where $n = d+m, \, W$ satisfies {\bf (Reg)}, the measure
$\bar  \mu$ is supported on a closed orbit $H  \cL_1 \subset
\ALN_n$ for a connected real algebraic group $H \subset \ASL_n(\R)$ and
$\LL_1 \in \ALN_n$.
There is an integer $k \geq d$,
a real number field $\KK$ and a 
$\KK$-algebraic group $\mathbf{G}$, such that the Levi subgroup of $H$
arises via restriction of 
scalars from $\mathbf{G}$ and $\KK$, and one of the following holds
for $\mathbf{G}$:
\begin{itemize}
\item $\mathbf{G} = \SL_k$ (as a $\KK$-group) and $n = k \cdot \deg (\KK/\Q).$
    \item $\mathbf{G} = \Symp_{2k}$ (as a $\KK$-group), and $d =2, \
      n=2k \cdot  \deg(\KK/\Q) $.  
  \end{itemize}
  Furthermore,
  in the linear (resp.\,affine) case $\mu$ is invariant under
 none of (resp., all of) the translations by nonzero elements of $\Vphys$.
\end{thm}

Here the group $\Symp_{2k}$ is the
group preserving the standard symplectic form in $2k$ variables; we
caution the reader that this group is sometimes denoted by $\Symp_k$ in
the literature. As we will see in Proposition \ref{prop: new report U}, any choice of
$\KK$ and $\mathbf{G}$ satisfying the 
description in Theorem \ref{thm: classification} gives rise to an
affine and a linear RMS measure. We note that the vertex sets of the famous Ammann-Beenker and Penrose tililngs, which are well-known to have representations as cut-and-project constructions, are associated with the real quadratic fields $\KK=\QQ(\sqrt{2})$ and $\KK=\QQ(\sqrt{5})$, resepctively, with $d=2$ and $\mathbf{G} = \SL_2$, see also \S \ref{sec:consequences}.

Theorem \ref{thm: classification} is actually a combination of two
separate results. The first extends work of Marklof and Str\"ombergsson
\cite{MS}. They introduced the pushforward $\bar \mu \mapsto \mu$
described above, where $\bar \mu$ is a homogeneous measure on
$\ALN_n$, and 
noted that the measures $\bar \mu$ could be classified using Ratner's
work. Our contribution in this regard (see Theorem \ref{thm: MS
  classification}) is to give a full list of the
measures $\bar \mu$ which can arise. The second result,
contained in our Theorem \ref{thm: MS surjective}, is that this construction is the only
way to obtain RMS measures according to our definition (which is given
in terms of $\Vphys$ rather than $\ALN_n$).

\subsection{Formulae of Siegel-Weil and Rogers}
In geometry of numbers, computations with the Haar-Siegel probability
measure on $\LLN_n$ are greatly simplified by the Siegel summation
formula \cite{Siegel}, according to which for $f \in C_c(\R^n)$,
\eq{eq: Siegel summation formula1}{
\int_{\LLN_n} \hat{f}(\LL) \, dm(\LL) = \int_{\R^n} f(x) \, d\Leb(x), \ \
\text{ where } \hat{f}(\LL) 
= \sum_{v \in \LL\sm \{0\}} f(v).
}
Here 
$m$ is the Haar-Siegel probability measure on $\LLN_n$, and $\Leb$ is
the Lebesgue measure on $\R^n$. The
analogous formula for RMS measures was proved in
\cite{MS}. Namely\footnote{Our notations differ  slightly from
  those of \cite{MS}, but the result as stated here
  can be easily shown to be equivalent to the one in \cite{MS}.}, 
suppose $\mu$ is an RMS measure, and for each $\Lambda \in \supp\, \mu$, and
for $f \in C_c(\R^d)$, set
\eq{eq: Siegel Veech transform}{
  \hat{f}(\Lambda) \df \left\{ \begin{matrix} 
      \displaystyle{\sum_{v \in \Lambda \sm \{0\}} f(v)} & \mu \text{ is linear } \\
      \displaystyle{\sum_{v \in \Lambda} f(v)} & \mu \text{ is affine. } 
      \end{matrix} \right.
  }
We will refer to $\hat{f}$ as
the {\em Siegel-Veech transform} of $f$. Then it is shown in
\cite{MS, MS_correction}, 
that for an explicitly computable constant $c>0$, for any 
$f \in C_c(\R^d)$ one has 
\eq{eq: Siegel summation MS}{
\int \hat{f}(\Lambda) \, d\mu(\Lambda) = c \, \int_{\R^d} f(x) \, d\vol(x).
}
A first step in the proof of \equ{eq: Siegel summation MS} is to show
that $\hat{f}$ is integrable, i.e., belongs to $L^1(\mu)$. As a
corollary of Theorem \ref{thm: classification}, and using reduction theory for lattices
in algebraic groups, we strengthen this and obtain the precise
integrability exponent of the Siegel-Veech transform, as follows:

\begin{thm}\name{thm: reduction theory}
Let $\mu$ be an RMS measure, let $\mathbf{G}$ and $\KK$ be as in Theorem
\ref{thm: classification}, let $r \df \rank_{\KK}(\mathbf{G})$ denote
the $\KK$-rank of $\mathbf{G}$, and define 
\eq{eq: p bounded}{ q_\mu \df \left\{ 
\begin{matrix} r +1 & \mu \text{ is linear} \\ r+2 &  \mu \text{ is affine.} \end{matrix} \right.
}
Then for any $f \in C_c(\R^d)$ and any $p<q_\mu$ we have $\hat{f} \in
L^p(\mu)$. Moreover, if the window $W$ contains a neighborhood of the
origin in $\Vint$, there are $f \in
  C_c(\R^d)$ for which $\hat f \notin L^{q_\mu}(\mu)$. 
\end{thm}
The proof involves integrating some characters over a Siegel set for
a homogeneous subspace of $\LLN_n$. 
The special case for  which $\KK=\QQ, \, \mathbf{G}=\SL_k$ and the measure
$\mu$ is linear was carried out in \cite[Lemma 3.10]{EMM}.
Note that
\eq{eq: K rank}{
\rank_{\KK}(\mathbf{G}) = \left\{ \begin{matrix} k-1 & \text{ if }
    \mathbf{G}=\SL_k\\ k  & \ \text{ if }  \mathbf{G}= \Symp_{2k}. \end{matrix} \right.
}
We will say that the RMS measure $\mu$ 
is {\em of higher rank} when $q_\mu \geq 3$; in light of the above this
happens unless $d=2, \, \mathbf{G} = \SL_2$, and $\mu$ is linear. It
follows immediately from 
Theorem \ref{thm: reduction 
  theory} that $\hat f \in L^1(\mu)$, and in the higher-rank case, that $\hat{f} \in L^2(\mu)$.

The proof of \equ{eq: Siegel summation MS} given in \cite{MS}  follows
a strategy of Veech \cite{Veech_siegel_measures}, and relies on a
difficult result of Shah \cite{Shah}. Following Weil \cite{Weil}, we
will reprove the result with a more elementary argument. Combined with
Theorem \ref{thm: reduction theory}, the argument gives a
strengthening of \equ{eq: Siegel summation MS}.


Given $p \in \N$, write 
$\bigoplus_1^p \R^{d} = 
\R^{dp}$, and for a compactly supported function $f$ on $\R^{dp}$, define
\eq{eq: def fphat}{
  \fphat
(\Lambda) \df \left\{ \begin{matrix} \displaystyle{\sum_{v_1, \ldots, v_p \in
      \Lambda \sm \{0\}} f(v_1,
\ldots, v_p) }& \mu \text{ is linear} \\ \displaystyle{\sum_{v_1, \ldots, v_p \in \Lambda} f(v_1,
\ldots, v_p)} & \mu \text{ is affine.} \end{matrix} \right.}

\begin{thm}\name{thm: weil}
  Let $\mu$ be an RMS measure, 
  and suppose $p <q_\mu$ where $q_\mu$ is as in \equ{eq: p
    bounded}. Then there is a countable collection $\{\tau_{\ee} : \ee
  \in \E\}$ of Borel measures on $\R^{dp}$ such that 
$\tau \df \sum \tau_\ee$ is locally finite, and 
for every $f \in L^1(\tau)$ we have
$$
\int \fphat \,
d\mu = \int_{\R^{dp}} f \, d\tau < \infty.
$$
The measures $\tau_\ee$ are  $H$-c\&p-algebraic, for the group $H$
appearing in Theorem \ref{thm: MS classification} (see Definition \ref{def:
   c and p algebraic measure}). 
  \end{thm}
This result is inspired by several results of Rogers for lattices, see
e.g. \cite[Thm. 4]{Rogers-acta}. Loosely speaking, c\&p-algebraic measures are
images of algebraically 
 defined measures on $\R^{np}$ under a natural map associated with the
 cut-and-project construction.

 Theorems \ref{thm: reduction theory} and \ref{thm: weil} will be
 deduced from their more general counterparts
 Theorems \ref{thm: reduction theory more general} and
 \ref{thm: weil more general}, which deal with the homogenous subspace $H
  \LL_1 \subset \ALN_n$ arising in Theorem \ref{thm: classification}.

 \subsection{Rogers-type bound on the second moment}
A fundamental problem in geometry of numbers is to control the higher
moments of random variables associated with the Haar-Siegel measure on
the space $\LLN_n$. In particular, regarding the second moment, the
following important estimate was proved in \cite{Rogers-acta,
  Rogers_number, Schmidt_metrical}: 
for the Haar-Siegel measure $m$ on $\LLN_n$, $n \geq 3$ there is a
constant $C>0$ such that for any function $f \in C_c(\R^n)$ taking
values in $[0,1]$ we have 
$$
\int_{\LLN_n}  \left| \hat f (x) - \int_{\LLN_n} \hat f \, dm \right|^2  \,
dm(x) \leq C \int_{\R^n} f \, d\vol ,
$$
where $\hat f$ is as in \equ{eq: Siegel summation formula1}.
We will prove an analogous result for RMS measures of higher rank.

\begin{thm}\name{thm: Rogers bound RMS}
Let $\mu$ be an RMS measure of higher rank. For $p=2$ let $\tau$ be
the measure as in Theorem \ref{thm: weil}. In the notation of Theorem
\ref{thm: classification}, assume that
\eq{eq: new condition}{
\mathbf{G} = \SL_k, \ \text{ or } \mu \text{ is affine.
  }}
Then there is 
$C>0$ such that for any Borel function $f: \R^d
\to [0,1]$ belonging to $L^1(\tau)$ we have 
\eq{eq: Rogers variance bound}{
\int_{\Cl(\R^d)}  \left| \hat f(x) - \int_{\Cl(\R^d)}
  \hat f \, d\mu \right|^2  \,
d\mu(x) \leq C \int_{\R^d} f \, d\vol.
}
\end{thm}
 The
case in which \equ{eq: new condition} fails, that is, $\mu$ is linear and
$\mathbf{G} = \Symp_{2k}$, and in which in addition $\KK = \QQ$, is treated in
\cite{Kelmer_Yu}, where a similar bound is obtained. 
The symplectic case with $\KK$ a proper field extension of $\Q$
is more involved, and we hope to investigate it further in future
work. 

There have been several recent papers proving an estimate like \equ{eq:
  Rogers variance bound}
for homogeneous measures associated with various algebraic groups. See
\cite{Kleinbock_Skenderi} and references
therein. The alert reader will
have noted that, even though the measure $\mu$ is the pushforward of a
measure supported on a homogeneous space $H\LL_1$, we prove the bound
\equ{eq: Rogers
  variance bound} for functions defined on $\Cl(\R^d)$ rather than
on $H\LL_1$. Indeed, while we expect such a stronger
result to be true, it requires a more
careful analysis than the one needed for our
application.

 \subsection{The Schmidt theorem for cut-and-project sets, and
   patch-counting} 
It is well-known that every irreducible cut-and-project set
$\Lambda$ has a  {\em
   density}
 \eq{eq: density exists}{
   D(\Lambda) \df \lim_{T \to \infty} \frac{\# \left(\Lambda
       \cap B(0,T) \right)}{\vol(B(0,T))}=\frac{\vol(W)}{\covol(\LL)},
   }
   where $\Lambda = \Lambda(\LL, W)$, $\vol(W)$ is the volume of
   $W$, and $\covol(\LL)$ is the covolume of $\LL$ (for two proofs,
   which are valid for a larger class of nice sets in place of
   $B(0,T)$, see
   \cite{Moody_uniform} and \cite[\S3]{MS}, and see references therein).  In particular, the
   limit exists and is positive. 
Following Schmidt \cite{Schmidt_metrical},
we would like to strengthen this result and allow counting in even more
general shapes, and with a bound on the rate of convergence. We say that a
collection of Borel subsets $\{\Omega_T : T \in \R_+\}$ of $\R^d$ is an {\em
  unbounded ordered family} if 
\begin{itemize}
  \item
$
0 \leq T_1 \leq T_2 \ \implies \Omega_{T_1} \subset \Omega_{T_2};
$
\item
  For all $T>0$, $\vol(\Omega_T) < \infty$;
\item
  $\vol(\Omega_T) \to_{T \to
   \infty} \infty
 $; and
 \item For all large enough $V>0$ there is $T$ such that
   $\vol(\Omega_T) = V$. 
\end{itemize}

\begin{thm}
  \name{thm: Schmidt analogue 1}
  Let $\mu$ be an RMS measure of higher rank, such that \equ{eq: new
    condition} holds.
   Then for every $\vre>0$, for every unbounded ordered family
$\{\Omega_T\}$, for $\mu$-a.e.\,cut-and-project set $\Lambda$,
\eq{eq: density with rate}{
\# \left(\Omega_T\cap \Lambda \right)
=
D(\Lambda) \cdot
\vol(\Omega_T) + O\left(\vol(\Omega_T)^{\frac12 + \vre} \right). 
}
\end{thm}

This result is a direct analogue of Schmidt's result for lattices, and
its proof follows \cite{Schmidt_metrical}. In the special case
$\Omega_T = B(0,T)$, we obtain an estimate for the rate of convergence in
\equ{eq: density exists}, valid for $\mu$-a.e.\,cut-and-project
set. For related work see
\cite{HKW}. Note that for $B(0,T)$, and for lattices, G\"otze
\cite{gotze_icm} has conjectured that an error estimate
$O\left(\vol(B(0,T))^{\frac{1}{2} - \frac{1}{2d}+\vre} \right)$
should hold.

Even for $\Omega_T = B(0,T)$, one
cannot expect \equ{eq: density with rate} to hold for {\em all}
cut-and-project sets; in fact, a Baire category argument as in
\cite[\S 9]{HKW} can be used to show that 
for any error function $E(T)$ with $E(T) 
= o(T^d)$ there are cut-and-project sets for which, along a
subsequence $T_n \to \infty$,
$$
\left| \# \left(B(0,T_n) \cap \Lambda \right) - 
D(\Lambda) \cdot
\vol(B(0,T_n)) \right| \geq E(T_n).
$$
Thus, it is an interesting open problem to obtain error estimates like
\equ{eq: density with rate} for explicit cut-and-project sets. Note
that for explicit cut-and-project sets which can also be described via
substitution tilings, such as the vertex set of a Penrose tiling,
there has been a lot of work in this direction, 
see \cite{Sol} and references therein. 

We now discuss patch counting, which is a refinement
which makes sense for cut-and-project sets
but not for lattices. For any discrete set $\Lambda \subset \R^d$, any
point $x \in
\Lambda$  and any $R>0$, we refer to the set
$$\PP_{\Lambda, R}(x) \df B(0, R) \cap (\Lambda -x)$$
as the {\em $R$-patch of $\Lambda$
  at $x$.} Two points $x_1, x_2 \in \Lambda$ are said to be {\em 
  $R$-patch equivalent} if $\PP_{\Lambda, R}(x_1) = \PP_{\Lambda,
  R}(x_2) $. It is well-known that any cut-and-project set $\Lambda$
is of {\em finite
local complexity}, which means that for any $R>0$,
$$
\# \{\PP_{\Lambda, R}(x) : x \in \Lambda\} < \infty. 
$$
Furthermore, it is known that whenever $\PP_0 = \PP_{\Lambda, R}(x_0)$ for
some $x_0 \in \Lambda$ and some $R>0$, the {\em density} or {\em
  absolute frequency} 
\eq{eq: patch density}{
D(\Lambda, \PP_0) = \lim_{T \to \infty} \frac{\# \{x \in \Lambda \cap
  B(0,T) : \PP_{\Lambda, R}(x)  = \PP_0 \}}{\vol(B(0,T))}
}
exists; in fact, the set in the numerator of \equ{eq:
  patch density} is itself a cut-and-project set, see
\cite[Cor.\,7.3]{BaakeGrimm}. 
Our analysis makes it possible to obtain an analogue of Theorem
\ref{thm: Schmidt analogue 1}
for counting patches, namely:

\begin{thm}\name{thm: Schmidt in patches}
 Let $\mu$ be an RMS measure of higher rank, for which \equ{eq: new
   condition} holds. For any $\delta>0$, set $\theta_0 \df
 \frac{\delta}{m+2\delta}$, where $m = \dim  \Vint$. 
Suppose the window $W \subset \Vint$ in the cut-and-project
construction satisfies $\dim_B(\partial W) \leq m
-\delta$, where $\dim_B$ denotes the upper box dimension (see \S
\ref{sec: counting}). Then for every
unbounded ordered  family $\{\Omega_T\}$ in $\R^d$, for 
$\mu$-a.e.\  $\Lambda$, for any patch $\PP_0 = \PP_{\Lambda, R}(x_0)$,
and any $\theta \in (0, \theta_0)$, we
have
\eq{eq: slightly weaker}{
\# \{x \in \Omega_T \cap \Lambda: \PP_{\Lambda, R}(x) = \PP_0\} =
D(\Lambda, \PP_0) \, \vol(\Omega_T) + O\left( 
\vol\left(\Omega_T \right)^{1-\theta}\right).}
  \end{thm}

  For additional results on effective error terms for patch-counting
  in cut-and-project sets, see \cite{HJKW}.
  
  \subsection{Acknowledgements} We are grateful to Mikhail
  Borovoi, Manfred Einsiedler, Dmitry Kleinbock, Henna Koivusalo, Jens Marklof, Dave
  Morris, Michel Skenderi, and Andreas Str\"ombergsson for useful
  discussions. Specifically, Morris supplied most of the arguments of
  Theorem \ref{thm: from morris}, 
  Borovoi supplied an argument used in Step 4 in the proof of Lemma
  \ref{lem:classification}, and Einsiedler supplied arguments for
  Lemmas \ref{lem:classification} and \ref{lem:assref}. 
  We
  gratefully acknowledge support of BSF grant 
  2016256, ISF grants 2919/19, 1570/17, 1149/18 and 264/22, 
  Swiss National Science
  Foundation 168823 and European Research Council 754475,  and the David and Rosa
  Orzen Endowment Fund. We thank
  the anonymous referees for a careful reading of the paper and for
  many helpful comments and suggestions. 

\section{Basics}
\name{sec: basics}
\subsection{Cut-and-project sets}\name{subsec: c and p sets}
In the literature, different authors impose slightly different
assumptions on the data in the cut-and-project construction.  For related
discussions, see \cite{BaakeGrimm, Moody_survey, MS}.  Here are the
assumptions which will be relevant in this paper:

\begin{itemize}
        \item[\textbf{(D)}] $\piint(\cL)$ is dense in $\Vint$.
        \item[\textbf{(I)}] $\piphys|_{\cL}$ is injective.
          \item[\textbf{(Reg)}] 
The window $W$ is Borel measurable, bounded, has non-empty interior, and its boundary
$\partial W$ has zero measure with 
respect to Lebesgue measure on $\Vint$. 
\end{itemize}
We will say that the construction is {\em irreducible} if {\bf (D), (I)}
and {\bf (Reg)} hold. 

In the literature, a more general cut-and-project scheme is discussed,
in which the groups $\Vphys \simeq \R^d, \Vint \simeq \R^m$ may be
replaced  with general locally
compact abelian groups. Note that if {\bf (D)} fails, we can replace
$\Vint$ with $\overline {\piint(\LL)}$, which is a proper subgroup of
$\Vint$, while if {\bf (I)} fails, we can replace $\Vint$ with $\Vint
/ (\LL \cap \ker \piphys)$. In both cases one can obtain the same set
using smaller groups. Note that when {\bf (D)} fails, the group $\overline
{\piint(\LL)}$ might be disconnected, and in that case, using {\bf (Reg)} we see that only finitely
many of its connected components will intersect $W$, and
$\Lambda(\mathcal{L},W)$ will have a description as a finite union of
cut-and-projects sets with an internal space of smaller dimension. 

Regarding the regularity assumptions on $W$, note that  if
no 
regularity assumptions are imposed, one can let $\Lambda$ be an
arbitrary subset of $\piphys(\LL)$ by letting $W$ be equal to $\piint
\left(\LL \cap \pi^{-1}_{\mathrm{phys}}(\Lambda) \right )$. Also, the assumption
that $W$ is bounded (respectively, has nonempty interior) implies that
$\Lambda$ is uniformly discrete (respectively, relatively dense).

Finally, note that it is not
$W$ that plays a role in \equ{eq: def c and p set}, but rather
$\pi^{-1}_{\mathrm{int}}(W)$. In particular, if convenient, one can replace the space
$\Vint$ with any space $V'_{\mathrm{int}}$ which is complementary to $\Vphys$,
and with the obvious notations, replace $W$ with $W' \df \pi'_{\mathrm{int}}
(\pi^{-1}_{\mathrm{int}}(W))$. Put otherwise, it would have been more natural to
think of $W$ as being a subset of the quotient space $\R^n/\Vphys$. We
refrain from doing so to avoid conflict with established
conventions. 

\subsection{Chabauty-Fell topology}\name{subsec: CF}
Let $\Cl(\R^d)$ denote the collection of all closed subsets of
$\R^d$. Equip $\Cl(\R^d)$ with the topology induced by the following
metric, which we will call the {\em Chabauty-Fell
 metric}: for $Y_0, Y_1
\in \Cl(\R^d)$, $d(Y_0, Y_1)$ is the infimum 
of all $\vre \in (0,1)$ for which, for both $ i=0,1,$
$$
Y_i \cap B \left(0, \vre^{-1} \right)
\text{ is contained in the } \vre\text{-neighborhood of } Y_{1-i},
$$
and $d(Y_0, Y_1) =1$ if there is no such $\vre$. 
It is known that
with this metric, $\Cl(\R^d)$ 
is a compact metric 
space. In this paper,
closures of collections in $\Cl(\R^d)$ and 
continuity of maps with image in $\Cl(\R^d)$ will always refer to this
topology, and all measures will be regular measures on the Borel
$\sigma$-algebra induced by this topology. We note that in the quasicrystals literature this topology is often referred to  as the {\em local rubber  topology} or the  {\em natural topology} .

We note that there are many topologies on the set of closed subsets
$\Cl(X)$ of a topological space $X$. The Chabauty-Fell metric was
introduced by Chabauty \cite{Chabauty} for $X = \R^d$ as well as for $X$ a locally
compact second countable group, and by Fell 
\cite{Fell} for general spaces $X$, particularly spaces arising in
functional analysis. See also \cite{Lenz_Stollmann}, where the
connection to the Hausdorff metric is elucidated via stereographic
projection. Many of the different topologies in 
the literature coincide on $\Cl\left(\R^d\right)$. Two notable
exceptions are the
Hausdorff topology, which is defined on the collection of {\em
  nonempty} closed subsets of $X$, and the weak-* topology of Borel
measures on $\R^d$, studied in
\cite{Veech_siegel_measures, MSnew}, satisfying a certain growth
condition and restricted to point processes. 
See \cite{Beer} for a comprehensive discussion of topologies on $\Cl(X)$.

We will need the following fact, which is well-known to experts, but
for which we could not find a reference (see \cite[\S 5.3]{MSnew} for
a related discussion):
\begin{prop}\name{prop: continuity Chabauty-Fell}
  Suppose $W$ is Borel measurable and bounded. Then the map
  \eq{eq: def Psi}{
    \Psi: \ALN_n \to \Cl(\R^d), \ \ \ \ \Psi(\LL)  \df \Lambda(\LL, W)
    }
  is a Borel map, and is continuous at any $\LL$
for which $\piint(\LL) \cap \partial W = \varnothing. $
  \end{prop}
  \begin{proof} 
We first prove the second assertion, that is, we assume that $\piint(\LL)
\cap \partial W = \varnothing $ and suppose by contradiction that
$\LL_j \to \LL$ in $\ALN_n$ but 
$\Psi(\LL_j) \not \to \Lambda \df \Psi(\LL).$ By passing to a
subsequence and using the definition of the Chabauty-Fell metric on
$\Cl(\R^d)$, we can assume that there is $\vre>0$ such that for all $j$,
one of the following holds:
\begin{itemize}\item[(a)]
There is $v \in \Lambda, \, \|v\| \leq \vre^{-1}$ such that for all $j$,
$\Psi(\LL_j)$ does not contain a point within distance $\vre$ of $v$.
\item[(b)]
  There is $v_j \in \Psi(\LL_j)$ such that $v_j \to v$, where $\|v\|
  \leq \vre^{-1}$, and $v \notin \Lambda$. 
\end{itemize}
In case (a), there is $u \in \LL$ such that $v = \piphys(u)$ and
$\piint(u) \in W$. By assumption $\piint(u)$ is in the interior of
$W$. Since $\LL_j \to \LL$ there is $u_j \in \LL_j$ such that $u_j \to
u$ and for large enough $j$, $\piint(u_j) \in W$ and hence $v_j =
\piphys(u_j) \in \Psi(\LL_j)$. Clearly $v_j \to v$ and we have a
contradiction.

In case (b), we let $u_j \in \LL_j$ such that $v_j =\piphys(u_j)
$. Then the images of $v_j$ under both projections $\piphys, \,
\piint$ are bounded sequences, and hence the sequence $(u_j)$ is also
bounded. Passing 
to a subsequence and using that $\LL_j \to \LL$ we can assume  $u_j 
\to u$ for some $u \in \LL$. Since $\piint(u_j) \in W$ for each $j$,
$\piint(u) \in \overline{W}$ and hence, by our assumption, $\piint(u)$
belongs to the interior of $W$, and in particular to $W$. This implies
that $v = \piphys(u) \in \Lambda$, a contradiction.

We now prove that $\Psi$ is a Borel measurable map. For this it is
enough to show that $\Psi^{-1}(B)$ is measurable in $\ALN_n$,
whenever $B=B(\Lambda, \vre)$ is the $\vre$-ball with respect to the
Chabauty-Fell metric centered at $\Lambda
= \Psi(\LL) \in \Cl(\R^d)$. Let 
$$F_1 \df \left\{x \in \LL : \piphys(x) \in  B\left(0, \vre^{-1}\right), \,
  \piint(x) \in W \right\}$$ and
$$
F_2 \df \Lambda \cap B\left(0, \vre^{-1}+
  \vre\right).$$
Then the definition of the Chabauty-Fell metric 
gives that $\LL'$ belongs to $\Psi^{-1}(B)$ if and only if for any $u_1
\in F_1$, there is $u'_1 \in \LL'$ with $\piint(u'_1) \in W$ and
$\|\piphys(u_1) - \piphys(u_1')\|< \vre$, and additionally, for any
$u'_1 \in \LL'$ with $\piint(u'_1) \in W$ and $\|\piphys(u'_1)\| <
\vre^{-1}$ there is $v \in F_2$ with $\|\piphys(u'_1) - v\|< \vre$. Since
lattices are countable, $F_1, F_2$ are finite, and $W \subset \Vint$
is Borel measurable, this shows that
$\Psi^{-1}(B)$ is described by countably
many measurable conditions. 
\end{proof}

We use this to obtain a useful continuity property for measures. Given
a topological space $X$, we denote by $\Prob(X)$ the space of regular
Borel probability measures. We equip $\Prob(X)$ with the weak-*
topology. Any Borel map $f : X \to Y$ induces a map $f_* : \Prob(X) \to
\Prob(Y)$ defined by $f_* \mu = \mu \circ f^{-1}$. 

\begin{cor}\name{cor: continuity} Let $\Psi$ 
  be as in \equ{eq: def Psi}.
  Then any $\bar \mu \in \Prob(\ALN_n) $ for which 
   \eq{eq: condition on bar mu}{\bar
  	\mu \left( \left\{\LL \in \ALN_n: \piint(\LL) \cap \partial W \neq
  	\varnothing \right\} \right) =0.}
is a continuity point for $\Psi_*$. In particular, this holds if $\bar \mu$ 
  is invariant under
 translations by elements of $\Vint \simeq \R^m$ and $\partial W$ has
 zero Lebesgue measure.
  \end{cor}

  \begin{proof}
    Suppose $\bar \mu_j \to \bar \mu$ in $\Prob(\ALN_n)$,
and let $\mu_j,\, \mu$ denote respectively the pushforwards $\Psi_* \bar
\mu_j, \, \Psi_* \bar \mu$. To establish continuity of $\Psi_*$ we
need to show $\mu_j \to \mu$.   Since $\bar \mu_j \to \bar \mu$, we
have     $\int g \, d\bar \mu_j \to \int g\, d\bar \mu$ for any $g \in
    C_c(\ALN_n)$. By the Portmanteau theorem this also holds for any
    $g$ which is bounded, compactly supported,  and for which the set
    of discontinuity points 
    has $\bar \mu$-measure zero. 
    Let $f$ be a continuous function on
    $\Cl(\R^d)$ and let $\bar f = f \circ \Psi$. 
Then $\bar f$ is
    continuous at $\bar\mu$-a.e.\  point, by Proposition \ref{prop:
      continuity Chabauty-Fell}.
    The Portmanteu theorem then ensures that
    $$
\int_{\Cl(\R^d)} f d\mu_j  = \int_{ \ALN_n} \bar f d\bar \mu_j \to
\int_{ \ALN_n} \bar f d\bar \mu  = \int_{\Cl(\R^d)} f d\mu.
$$
That is, $\mu_j \to \mu$, as required.

    For the last assertion, assuming that $\bar \mu$ is invariant under
translations by elements of $\Vint$, we need to show that \equ{eq:
  condition on bar mu} is satisfied. Letting $\mathbbm{1}_{\partial
  W}, \, m_{\Vint}$ denote respectively the indicator of $\partial W$
and Lebesgue measure on $\Vint$, and letting $B \subset \Vint$ be a
measurable set of finite and positive measure, we have by Fubini that 
\[
  \begin{split}
& \bar
    \mu \left( \left\{\LL \in \ALN_n: \piint(\LL) \cap \partial W \neq
        \varnothing \right\} \right)
    \\
    = &
\int \left[ \frac{1}{m_{\Vint}(B)}\int_{B} 
\mathbbm{1}_{\partial W}\circ \piint(\LL+x) \, dm_{\Vint}(x) \right]
\, d\bar\mu(\LL). 
\end{split}
\]
It therefore suffices to show that for any $\LL$, 
$$
m_{\Vint} \left(\{x \in \Vint : \piint(\LL + x)  \cap \partial W \neq
  \varnothing\}  \right) =0;  
$$
and indeed, this follows immediately from the countability of $\LL$
and the assumption that
$m_{\Vint}(\partial W)=0$. 
    \end{proof}

\subsection{Ratner's Theorems}\name{subsec: ratner}
Ratner's measure classification and orbit-closure theorems
\cite{Ratner-annals} are fundamental results in homogeneous
dynamics. We recall them here, in the special cases which will be
important for us. A Borel probability measure $\nu$ on 
$\ALN_n$ (respectively, $\LLN_n$) is called {\em homogeneous} if there is $x_0$ in
$ \ALN_n$ (respectively, $\LLN_n$) and
a closed subgroup $H$ of $ \ASL_n(\R)$ (respectively, $\SL_n(\R)$) such that the
$H$-action preserves $\nu$, the orbit $Hx_0$ is closed and equal to
$\supp \, \nu,$ and $H_{x_0} \df \{h \in H: 
h x_0 =x_0\}$ is a lattice in $H$. When we want to stress the role
of $H$ we will say that $\nu$ is {\em $H$-homogeneous}. 

Recall that $\ASL_n(\R)$ (respectively, 
$\ASL_n(\Z)$) denotes the group of affine
transformations of $\R^n$ whose derivative has determinant one
(respectively, and which map the integer lattice $\Z^n$ to
itself), and that $\ALN_n$ is identified with 
$\ASL_n(\R) / \ASL_n(\Z)$, via the map which identifies the coset
represented by the affine map $\varphi$ with the grid
$\varphi(\Z^n)$. Similarly, we have an identification of $\LLN_n$ with
$\SL_n(\R)/\SL_n(\Z)$. 
We 
view the elements of $\ASL_n(\R)$ concretely as pairs $(g,v)$, where $g \in 
\SL_n(\R)$ and $x \in \R^n$ determine the map $x \mapsto gx +v$. In what follows two
subgroups of $\ASL_n(\R)$ play an important role, namely the groups
$\SL_d(\R)$ and $\ASL_d(\R)$, which we will denote alternately by
$F$, and 
embed concretely in $\ASL_n(\R)$ in the upper left hand corner. That is, 
in the case $F = \SL_d(\R)$,  $g \in F$ is identified with 
\eq{eq: embedding linear}{
\left(\left(\begin{matrix} g & \mathbf{0}_{d,m}\\ \mathbf{0}_{m,d} &
      \mathrm{Id}_m \end{matrix} 
  \right),  \mathbf{0}_n \right) 
}
and in the case $F = \ASL_d(\R)$,  $(g, v) \in F$ is identified 
with
\eq{eq: embedding affine}{
\left(\left(\begin{matrix} g & \mathbf{0}_{d,m}  \\ \mathbf{0}_{m,d}
      & \mathrm{Id}_m \end{matrix} 
  \right), \left (
  \begin{matrix} v \\ \mathbf{0}_m \end{matrix}\right ) \right). 
}
Here $\mathrm{Id}_m , \, \mathbf{0}_{k, \ell}, \, \mathbf{0}_k$ denote
respectively an identity  
matrix of size $m \times m$, a zero matrix of size $k \times \ell$,  
and the zero vector in $\R^k$. We will refer to the embeddings of
$\SL_d(\R)$ and $\ASL_d(\R)$ in $\ASL_n(\R)$, given by \equ{eq:
  embedding linear} and \equ{eq: embedding affine}, as the {\em top-left corner
  embeddings}. 

The following
is a special case of Ratner's result.

\begin{thm}[Ratner]\name{thm: Ratner measure classification}
Let $2 \leq d \leq n$, and let $F$ be equal to either $\ASL_d(\R)$
or $\SL_d(\R)$ (with the top-left corner embedding in $\ASL_n(\R)$).
Then any $F$-invariant 
ergodic measure $\nu$ on 
$\ALN_n$ is $H$-homogeneous, where $H$ is a closed
connected subgroup 
of $\ASL_n(\R)$ containing $F$. Every orbit-closure
$\overline{Fx}$ is equal to $\supp \, \nu$ for some homogeneous
measure $\nu$. The same conclusion holds for $\LLN_n$ and $F =
\SL_d(\R)$. 
\end{thm}

The following additional results were obtained in \cite{Shah-Annalen,
  Tomanov-Japan}: 

\begin{thm}[Shah, Tomanov]\name{thm: properties of L}
Let $\nu, H$ be as in Theorem \ref{thm: Ratner measure
  classification}, and let $x_0 = g_0\Z^n $ in $\ALN_n$ or $\LLN_n$ such that
$\supp \, \nu = Hx_0$.
Let $\mathbf{H}'$ be the smallest algebraic subgroup of $\ASL_n$ which is defined
over $\Q$ and contains $g_0^{-1}Fg_0$.
The solvable radical of $\mathbf{H}'$ is equal to the unipotent
radical of $\mathbf{H}'$, and letting
$\mathbf{H} = g_0\mathbf{H}' g_0^{-1}$, 
$H$ is equal to the connected component of the identity
in $\mathbf{H}_\R$.
\end{thm}

We will need a result of Shah which relies on Ratner's work (once more
this is a special case of a more general result).

\begin{thm}[\cite{Shah} ] \name{thm: Shah}
Let $F $ be equal to either $\ASL_d(\R)$
or $\SL_d(\R)$ as above, let $\{g_t\}$ be a one-parameter diagonalizable subgroup
of $\SL_d(\R)$, and let $U = \{g \in F: \lim_{t \to \infty} g_{-t} g
g_t \to e\}$ be the corresponding expanding horospherical
subgroup. Let $\Omega \subset U$ be a relatively compact open subset
of $U$ and let $m_U$ be the restriction of Haar measure to $U$, normalized so that
$m_U(\Omega)=1$. Then for every $x_0 
\in \ALN_n$, letting $\nu$ be the homogeneous 
measure such that $\supp \, \nu =\overline{F x_0}$, we have
$$
\int_{\Omega} (g_tu)_* \delta_{ x_0} \, dm_U(u) 
\to_{t\to \infty} \nu,
$$
where $\delta_{x_0}$ is the Dirac measure at $x_0$ and the
convergence is weak-* convergence in $\Prob(\ALN_n)$. 
  \end{thm}


%

\subsection{Number fields, geometric embeddings, and restriction of
  scalars}\name{subsec: number fields}
For more details on the material in this subsection we refer the reader to \cite{Weil,
  Rapinchuk_Platonov, Morris, EW_homdyn}. 

Let $\KK$ be a number field of degree $D = \deg(\KK/\QQ)$, and let 
$\OO= \OK$ be its ring of integers. Let $\sigma_1, \ldots,
\sigma_r, \sigma_{r+1}, \overline{\sigma_{r+1}},\ldots, 
\sigma_{r+s}, \overline{\sigma_{r+s}}$ be the field embeddings of
$\KK$ in $\CC$ where $r +2s 
= D$, $\sigma_1, \ldots, \sigma_r$ are real embeddings and
$\sigma_{r+1}, \ldots, \sigma_{r+s} $ are complex (non-real)
embeddings. An {\em order} in $\KK$ is a subring of
$\OO$ which is of rank $D$ as an additive group.
The {\em geometric embedding} or {\em Minkowski embedding} of
an order $\Delta$ is the set 
$$
\left\{\left(\sigma_1(x), \ldots, \sigma_r(x), \sigma_{r+1}(x), 
\ldots, \sigma_{r+s}(x) \right): x \in \Delta\right\}.
$$
It is a lattice in $\R^D \simeq \R^r
\times \CC^s.$ Note that the geometric embedding depends on a choice
of ordering of the field embeddings, 
and on representatives of each pair of complex conjugate
embeddings. Thus, when we speak of `the' geometric embeddings we will
consider this data as fixed.

An algebraic group $\mathbf{G}$ defined over $\KK$ (or $\KK$-algebraic
group) is a variety defined over $\KK$ such
that the multiplication and inversion maps $\mathbf{G} \times
\mathbf{G} \to \mathbf{G}, \ \mathbf{G} \to \mathbf{G}$ are
$\KK$-morphisms. A $\KK$-homomorphism of algebraic groups is a group
homomorphism which is a $\KK$-morphism of algebraic
varieties. We will work
only 
with {\em linear algebraic groups} which means that they are affine
varieties, i.e., for some $N$, they are the subset of affine space 
$\mathbb{A}^N$ satisfying a system of polynomial equations in $N$
variables. We will omit 
the word `linear' in the rest of the paper. A
typical example of a 
$\KK$-algebraic group is a Zariski closed
matrix group, that is, a subgroup of the
matrix group $\SL_m(\CC)$ for some $m$ described by
polynomial equations in  the matrix entries, with coefficients
in $\KK$. If $\mathbf{G}_i$ are $\KK$-algebraic groups
realized as subgroups of $\SL_{m_i}(\CC)$ for $i=1,2$, and $\varphi:
\mathbf{G}_1 \to \mathbf{G}_2$ is a $\KK$-homomorphism, then there is a map
$\hat \varphi : \SL_{m_1}(\CC) \to  \SL_{m_2}(\CC) $ which is polynomial 
in the matrix entries, with coefficients
in $\KK$, such that $\hat \varphi|_{\mathbf{G}_1} =
\varphi$. For any field $L \subset \CC$ containing $\KK$, we will
denote by $\mathbf{G}_L$ the collection of
$L$-points of $\mathbf{G}$. It is a subgroup of $\SL_m(L)$, if
$\mathbf{G}$ is realized as subgroup of $\SL_m(\CC)$.

We will do
the same for rings $L = \Z$ or $L = \OO$. In this case the group
$\mathbf{G}_L$ depends on the concrete realization of $\mathbf{G}$ as
a matrix group but the commensurability class of $\mathbf{G}_L$ is
independent of choices (recall that two subgroups $\Gamma_1, \Gamma_2$ of
some ambient group $G$ are {\em commensurable } if $[\Gamma_i: \Gamma_1
\cap \Gamma_2]<\infty$ for $i=1,2$). 
By a {\em real algebraic group} we will mean a subgroup of finite index
in $\mathbf{G}_{\R}$ for some $\KK$-algebraic group $\mathbf{G}$,
where $\KK \subset \R$.

The {\em restriction of scalars $\Res_{\KK/\Q}$} is a functor from the
category of $\KK$-algebraic groups to $\Q$-algebraic groups. 
Given an algebraic group $\mathbf{G}$ defined over $\KK$, there is an
algebraic group $\mathbf{H} = \Res_{\KK/\Q}(\mathbf{G})$ defined over
$\Q$, such that $\mathbf{H}_\Q$ is naturally identified with
$\mathbf{G}_\KK$. For any $\KK$-homomorphism of $\KK$-algebraic
groups $\varphi: \mathbf{G}_1
\to \mathbf{G}_2$ we have a $\Q$-homomorphism $\Res_{\KK/\Q}(\varphi):
\Res_{\KK/\Q}(\mathbf{G}_1) \to  \Res_{\KK/\Q}(\mathbf{G}_2)$. Given a
matrix representation of $\mathbf{G}$ there is a corresponding matrix
representation of $\Res_{\KK/\Q}(\mathbf{G})$, defined as
follows. We can realize $\KK$
(as a ring) as a subalgebra of the $\Q$-algebra of $D \times D$
matrices with entries in $\Q$, and this leads 
to a corresponding identification of $\SL_m(\KK)$ with a subgroup of $\SL_{mD}(\QQ)$.
A different choice of basis will produce a group that differs by a
$\SL_{mD}(\QQ)$-conjugate. 
Now suppose
$\mathbf{G} \subset \SL_m(\CC)$ is the solution set 
of polynomial equations $P_1, \ldots, P_\ell$ in the matrix
entries, with coefficients in $\KK$. Let $\hat P_1, \ldots, \hat
P_\ell$ be the 
  matrix valued polynomials where each $\KK$-coefficient is replaced
  by its $\operatorname{Mat}_{D\times D}(\QQ)$ representative, and
  each variable (previously a matrix coefficient of $\SL_{m}(\CC)$) is
  an $\operatorname{Mat}_{D\times D}(\CC)$-block of $\SL_{mD}(\CC)$. 
These polynomials together with the (linear) polynomials that ensure
that each $D\times D$ block is an element of the $\Q$-algebra $\KK$, 
have 
coefficients in $\Q$, and $\Res_{\KK/\Q}(\mathbf{G})$
is their solution set. 

The $\R$-points of $\mathbf{H} =
\Res_{\KK/\Q}(\mathbf{G})$ can be represented 
concretely as
\eq{eq: shape of restriction}{
  {}^{\sigma_1}\mathbf{G}_\R \times \cdots \times {}^{\sigma_r}\mathbf{G}_\R \times
  {}^{\sigma_{r+1}}\mathbf{G}_{\CC} \times \cdots \times {}^{\sigma_{r+s}}\mathbf{G}_\CC ,}
where 
${}^{\sigma_j}\mathbf{G}$ is the algebraic group defined by applying the field embedding
$\sigma_j$ to the polynomials in the matrix entries, with coefficients
in $\KK$, defining $\mathbf{G}$. Here, for a $\CC$-algebraic group
$\mathbf{M}$, $\mathbf{M}_{\CC}$ is a shorthand
notation for the $\CC$-points of $\mathbf{M}$, thought of as an
$\R$-group 
via the isomorphism $\CC \cong \R^2$. More explicitly, a polynomial
equation involving $m^2$ complex matrix entries $z_{ij} =
a_{ij}+\mathbf{i}b_{ij}$, where $i,j \in \{1, \ldots, m\}$, is replaced
with the same polynomial in the matrix algebra of $2 \times 2$ real
matrices, with each appearance of $z_{ij}$ replaced by 
$A^{(ij)} \df \begin{pmatrix} a_{ij}&b_{ij} \\ -b_{ij}&a_{ij}\end{pmatrix} \in
\mathrm{Mat}_{2 \times 2} (\R)$, and with the
$2m^2$ additional equations $(A^{(ij)})_{12} = -(A^{(ij)})_{21} ,
(A^{(ij)})_{11} = (A^{(ij)})_{22}.$ 
Furthermore, denoting by $\bar \Q$ the
algebraic closure of $\Q$, there is a conjugation of $\SL_{mD}(\bar
\Q)$ by an element with coefficients in the Galois closure of $\KK$,
so that
$\mathbf{H}(\bar \Q)$ is embedded in $\SL_{mD}(\bar \Q)$ in block form
with $r+s$ blocks, where each block contains one of the factors in
\equ{eq: shape of restriction}.  

Similarly, for a $\KK$-morphism
$\varphi: \mathbf{G}_1 \to \mathbf{G}_2$, the restriction to the
factor ${}^{\sigma_j}\mathbf{G}_\R$ in formula \equ{eq: shape of
  restriction}, of the $\Q$-morphism 
$\Res_{\KK/\Q}(\varphi): \Res_{\KK/\Q}(\mathbf{G}_1) \to
\Res_{\KK/\Q}(\mathbf{G}_2)$, 
is the map $\varphi_j$ obtained from $\varphi$ by applying the field 
embedding $\sigma_j$ to its coefficients. 
Thus, after writing both $\Res_{\KK/\Q}(\mathbf{G}_1)$ and
$\Res_{\KK/\Q}(\mathbf{G}_2)$ in product form as in \equ{eq:
  shape of restriction}, we have
\eq{eq: shape of representation}{
\Res_{\KK/\Q}(\varphi)(g_1, \ldots, g_{r+s}) = (\varphi_1(g_1),
\ldots, \varphi_{r+s}(g_{r+s})). 
}

We now note a connection between restriction of scalars,
 geometric embeddings of lattices, and the action on
$\LLN_n$. Suppose that  $\OO =
\OO_\KK$, $\Delta$ is an order in $\OO$, and let $\LL$ be the geometric embedding of
$\Delta$ in $\R^D$. For $m \in \N$ set $n = Dm$ and let
$$\LL' = c \cdot 
\underbrace{\LL \oplus \cdots  \oplus \LL}_{m \text{ copies}
} 
, $$
where we choose
the dilation factor $c$ so that $\LL' \in \LLN_n$, and we choose the
ordering of the indices so that
\eq{eq: def of lattice from order}{
\LL' \df c \left \{(\sigma_1(x), \ldots, \sigma_{r+s}(x)) : x \in \Delta^m \right \}.
}
Now suppose 
$\mathbf{G}$ is an algebraic $\KK$-group without $\KK$-characters,  $\varphi:
\mathbf{G} \to \SL_m$ is a $\KK$-morphism, and $\mathbf{H} \df
\Res_{\KK/\Q}(\mathbf{G})$. Since $\varphi$ is a $\KK$-morphism, there
is a finite-index subgroup of $\mathbf{G}_{\OO}$ whose image under
$\varphi$ is contained in $\SL_m(\OO)$, and hence 
preserves $\OO^m$. This implies that a finite index subgroup of
$\mathbf{H}_{\Z}$ preserves $\LL'$. Since $\mathbf{H}_{\Z}$ is a
lattice in $H \df \mathbf{H}_{\R}$ (see \cite[\S 13]{Borel_arithmetiques}), 
we find that $H\LL'$ is a closed orbit in $\LLN_n$ which is the
support of an $H$-homogeneous measure. 

\section{Classification of invariant measures}
\name{sec: classification}
Recall from the introduction that an affine (respectively, linear) RMS
measure $\mu$ is a probability measure on 
$\Cl(\R^d)$ which gives full measure to the collection of all
irreducible cut-and-project sets, 
and is invariant and ergodic under $F$, where
\eq{eq: def G0}{
  F \df \left \{\begin{matrix} \SL_d(\R) & \ \ \text{ if } \mu \text{
      is linear} \\ \ASL_d(\R) & \ \ \text{ if } \mu \text{
      is affine}  \end{matrix}  \right. }
is the stabilizer group of $\mu$. In this section we will give some
more background on RMS measures, and two assertions (Theorem \ref{thm:
MS classification}
and \ref{thm: MS surjective}) which together imply Theorem \ref{thm:
  classification}. The careful reader will have noticed that we gave
here a seemingly weaker definition of an affine RMS measure compared
to the introduction, by requiring it to be ergodic under $\ASL_d(\R)$
instead of $\SL_d(\R)$. However, these two definitions are equivalent
by the Howe-Moore 
ergodicity theorem (see \cite{ew}).

\subsection{RMS measures --- background and basic
  strategy}\label{subsec: strategy}
In order to motivate the definition of an RMS measure,
we recall some crucial observations of
\cite{MS}. 
Let $F$ be as in \equ{eq: def G0}.
Let $\R^n =
\Vphys \oplus \Vint, \piphys, \piint, \LL, W$ be the data involved in
a cut-and-project construction. 

The observations of \cite{MS} consist 
of the following:

\begin{itemize}
  \item
From the fact that $\piphys$ 
intertwines the action of $F$ on $\R^n$ (via the top-left corner
embedding in 
$\ASL_n(\R)$) and on $\R^d$, for the map $\Psi$  defined in \equ{eq: def Psi},
one obtains the equivariance property 
\eq{eq: equivariance Psi}{\Psi \circ g = g \circ \Psi
} for all $g \in F$; in other words, $g \Lambda(\LL, W) =
\Lambda(g\LL, W).$
\item
In particular, if we fix the data $\R^n = \Vphys \oplus \Vint,
\, W$, then the map $\Psi_*: \Prob(\ALN_n) \to \Prob(\Cl(\R^d))$
considered in Corollary \ref{cor: continuity} maps $F$-invariant measures to 
$F$-invariant measures.
\item
  Due to Ratner's work described in \S\ref{subsec: ratner}, 
  ergodic 
  $F$-invariant measures on $\ALN_n$ can be described in detail, in
  terms of certain real algebraic subgroups of $\ASL_n(\R)$. 
\item 
Theorem \ref{thm: Shah} and other results from homogeneous dynamics
can then be harnessed as a powerful tool for deriving information about
cut-and-project sets. 
\end{itemize}

In order to analyze measures on $\ALN_n$, a basic strategy is to work
first with the simpler space $\LLN_n$. 
 Let $$M \df \ASL_n(\R), \ \  \Gamma \df \ASL_n(\Z).$$ 
Recall that $\ALN_n$ is
 identified with $M/\Gamma$ and under this identification, a closed
 orbit $H \LL$ is identified with $H g
 \Gamma = gH_1
 \Gamma$, where $g \in M$ is such that $\LL = g\Z^n$, and $H_1 = g^{-1}Hg $. 
 Also let $$\underline M \df \SL_n(\R) \ \ \text{ and } \underline \Gamma
 \df \SL_n(\Z).$$ 
 We think of $\underline M$ concretely as the
 stabilizer of the origin in the action of $M$ on $\R^n$.
 Recall also that $\LLN_n$ is identified with
 $\underline M/\underline \Gamma$. Let
 \eq{eq: def projection}{
   \pi: M\to
 \underline M, \ \ \ \ \ \ \underline \pi: \ALN_n \to \LLN_n}
 denote respectively the natural quotient map, and the induced map on
 the quotients (which is well-defined since $\pi(\Gamma) = \underline 
 \Gamma$). The map $\pi$ is a $\Q$-morphism, and the map $\underline
 \pi$ is realized concretely by mapping a 
 grid $\LL$ to the 
 underlying lattice $\LL - \LL$ obtained by translating   $\LL$ so
 that it has a point at the origin. It satisfies an equivariance
 property
 \eq{eq: equivariance property}{
   \underline \pi(g \LL) = \pi(g) \underline \pi(\LL) \ \ \ \
   \text{(where $ g \in M, \ \LL \in \ALN_n $)}.}
 Every fiber of $\underline \pi$ is a torus and
 thus $\underline \pi$ is a proper map.

 We summarize the spaces and maps we use in the
 following diagram.

 \[
  \begin{tikzcd}[column sep = large]
   \ALN_n =M/\Gamma \arrow{dr}{\Psi} \arrow[swap]{d}{\underline \pi} 
  &  \\  \LLN_n = \underline M/\underline \Gamma &  \Cl(\R^d) 
  \end{tikzcd}
\]
Extending the terminology in the
introduction, a homogeneous measure $\bar \mu$ on $\ALN_n$ will be
called {\em affine} if it 
is $\ASL_d(\R)$-invariant, and {\em linear} if it is
$\SL_d(\R)$-invariant but not $\ASL_d(\R)$-invariant. Here
$\ASL_d(\R)$ and $\SL_d(\R)$ are embedded in $M$ via the
top-left corner embeddings \equ{eq: embedding affine} and \equ{eq:
  embedding linear}.
\subsection{The homogeneous measures arising from the $F$-action on 
  $\ALN_n$}\name{subsec: homogeneous measures}
In this section we state a more precise version of Theorem
\ref{thm: classification}. Suppose $k_0$ is a subfield of $\CC$. We
say that a $k_0$-algebraic group 
$\mathbf{H}$ is {\em $k_0$-almost simple} if any normal
$k_0$-subgroup $\mathbf{H}'$ satisfies $\dim \mathbf{H}'= \dim
\mathbf{H}$ or $\dim \mathbf{H}' =0$. In this case we will also say that
a subgroup of finite index of $\mathbf{H}_{k_0}$ is $k_0$-almost
simple. 



\begin{thm}
  \name{thm: MS classification}
Let $\bar\mu$ be an $F$-invariant ergodic measure on $\ALN_n$, and let
$H$ and $\LL_1$ denote respectively the subgroup of $M$ and the point
in $\ALN_n$ involved in Theorem \ref{thm:
  Ratner measure classification}; i.e., $\bar  \mu$ is $H$-invariant
and supported on the
closed orbit $H  \cL_1$. Let $g_1 \in M$ such that $\LL_1
= g_1 \Z^n$ and let $H_1 \df g_1^{-1} H g_1$. Assume also that $\LL_1$
satisfies conditions {\bf (D)} and {\bf (I)}. Then
$H, \, H_1$ and $\cL_1$ are described as
follows: 
  \begin{itemize}
\item[(i)]
  In the  linear case, $H_1$ is semisimple and 
  $\Q$-almost simple. In this case write $H' \df H_1$. 
    In the affine case, we can write $H_1$ as a semidirect product $H' \ltimes
  \R^n$ where $H'$ is semisimple and 
  $\Q$-almost simple, and $\R^n$ denotes the full group of translations of
  $\R^n$.
\item[(ii)]
The group $H'$ in (i) is the connected component of the identity in
the group of $\R$-points of
$\Res_{\KK/\Q}(\mathbf{G})$, where $\KK$ is a real
number field and $\mathbf{G}$ is a
$\KK$-group which is $\KK$-isomorphic to either $\SL_k$ or
$\Sp_{2k}$, for some $k \geq d$. In the case $\mathbf{G} = \SL_k$ we have
   $n = k \, \deg (\KK/\Q)$, and there is a subspace $V$ of $\R^n$
   of dimension $k$ containing $g_1^{-1}\Vphys$ which is $H'$-invariant and
   such that the action of $H'$ on $V$ gives the group $\SL(V)$. The case
  $\mathbf{G} = \Symp_{2k}$ only arises when $d=2$, and in that case $n = 2k \,
  \deg (\KK/\Q)$, and there is a subspace $V$ of $\R^n$ of dimension
  $2k$ equipped with a symplectic form $\omega'$ such that $V$ is
  $H'$-invariant, the action of $H'$ on $V$ gives the symplectic group
  $\Symp(V, \omega')$, and $V$ contains $g_1^{-1}\Vphys$ as a symplectic subspace.
    \end{itemize}
  \end{thm}

The proof will involve a reduction to the space $\LLN_n$ of
lattices. We introduce some notation and give some preparatory
  statements. 

  As in \S \ref{subsec: strategy}, let  $M = \ASL_n(\R), \, \Gamma = \ASL_n(\Z), \,
  \ALN_n = M/\Gamma$, so that the closed 
 orbit $H  \LL_1$ is identified with $H g_1 \Gamma = g_1H_1
 \Gamma$. By Theorem \ref{thm: Ratner measure classification},
 $\Gamma_{H_1}\df H_1 \cap \Gamma$ is a lattice in $H_1$ and $\bar \mu$
 is the pushforward of the 
 unique $H_1$-invariant probability measure on $H_1/\Gamma_{H_1}$ under the
 map $h\Gamma_{H_1}\mapsto  g_1 h \Gamma$. By Theorem \ref{thm:
   properties of L}, $H_1$
 is the connected component of the identity in the group of  real
 points of a $\Q$-algebraic  group. In particular 
 there are at most countably many possibilities for $H_1$.

 Also let $\underline M = \SL_n(\R), \, \underline \Gamma
 = \SL_n(\Z), \, \LLN_n = \underline M/\underline \Gamma$ as above, and let
 $\pi, \underline \pi$ be the maps in \equ{eq: def projection}.
   The orbit $\underline \pi(Hg_1 \Gamma )=
 \underline{H} \underline g_1 \underline \Gamma=  \underline g_1
 \underline{H}_1 \underline \Gamma$ is
 closed, where $\underline H_1 = \pi(H_1), \underline g_1 = \pi(g_1)$
 and $\underline H = \pi(H)$
 contains $\pi(F) \simeq \SL_d(\R)$. 

 We say that property {\bf (irred)} holds if 
there is no proper $\Q$-rational
subspace of $\R^n$ that is $\underline H_1$-invariant (for the linear
action by matrix multiplication).
Note that by Theorem~\ref{thm: properties of L}, $\underline H_1$
is the connected component of the identity in the group of real points
of the smallest $\Q$-subgroup 
of $\SL_n$ containing
$\underline{g}_1^{-1}\SL_d(\R)\underline{g}_1$, and thus {\bf
(irred)} is equivalent to 
requiring that there is no proper $\Q$-rational subspace of
$\R^n$ that is $\underline{g}_1^{-1}\SL_d(\R)\underline{g}_1$-invariant.

We now state an analogue of Theorem \ref{thm: MS
classification} for the action on $\LLN_n$. 

\begin{lem}\name{lem:classification}
Assume ${\bf (irred)}$ holds. Then $\underline H_1$ is the connected
component of the identity of the group of real points of a semisimple
$\Q$-algebraic group $\mathbf{H}$, satisfying the properties listed in
statement (ii) of Theorem \ref{thm: MS 
  classification} (for the group $H'$).
\end{lem}
Lemma \ref{lem:classification} is the main result of this section, and
its proof will be given below in \S \ref{subsec: preparations Lemma}
and \S \ref{subsec: proof Lemma}.

\begin{proof}[Proof of Theorem \ref{thm: MS classification} assuming
  Lemma \ref{lem:classification}]
  Let $\mathbf{\tilde{H}}$ be the smallest $\Q$-subgroup of $\ASL_n$ containing
$g_1^{-1}\SL_d(\R)g_1$, so that by Theorem~\ref{thm: properties of L},
we have $H_1=(\mathbf{\tilde{H}}_\R)^\circ$. Similarly, let
$\mathbf{H}$ be the smallest $\Q$-subgroup of $\SL_n$ containing
$g_1^{-1}\SL_d(\R)g_1$. We extend $\pi$ to a projection map of
algebraic groups defined over $\Q$, mapping $\Q$-subgroups to
$\Q$-subgroups (\cite[Cor I.1.4]{Borel1}). Then it follows from
minimality of $\mathbf{H}$ and $\mathbf{\tilde{H}}$, that
$\pi(\mathbf{\tilde{H}})=\mathbf{H}$. 

As we will see in Lemma \ref{lem:assref}, under the assumptions of
Theorem \ref{thm: MS classification}, condition {\bf (irred)}
holds. In particular, the conclusion of Lemma~\ref{lem:classification}
applies. 
Hence $\mathbf{H}$ is semisimple.

Let $\mathbf{U}$ be the unipotent radical of $\mathbf{\tilde{H}}$.
Then
$\mathbf{U} \subset \ker \pi$,
and since $\ker \pi \cap \mathbf{\tilde{H}}$ is a unipotent normal subgroup,
$\mathbf{U} = \ker \pi \cap \mathbf{\tilde{H}}$.
This means that in the affine map determined by $h \in H_1$ on
$\R^n$, $\pi(h)$ is the linear part, and $U=\mathbf{U}_\R$ acts on $\R^n$ by
translations.
This implies the equality 
\eq{eq: def V0}{
 \spa \{u(x) -x : x \in \R^n, u \in U\} = \spa \{u(0): u \in
U\},}
and we denote the subspace of $\R^n$ appearing in \equ{eq: def V0} by $V_0$. 
Clearly, $V_0$ are the real points of a $\Q$-subspace of $\CC^n$ since $\mathbf{U}$ is defined over $\Q$.

Since $H_1$
normalizes $U$, $V_0$ is $H_1$-invariant, and since $\underline H_1 =
\pi(H_1)$ is the group of linear parts of elements of $H_1$, $\underline H_1$
also preserves $V_0$.
By {\bf (irred)} we must have  
$V_0 = \{0\}$ or $V_0 = \R^n$. If $V_0 = \{0\}$ then $U = \{0\}$. If $V_0
= \R^n$ then $U$ contains translations in $n$ linearly independent
directions and hence $U \simeq \R^n$ is the entire group of
translations of $\R^n$. This gives the description of the
translational part of $H_1$, in assertion (i). 
Assertion (ii) follows from Lemma \ref{lem:classification}.
\end{proof}

The next proposition shows that all the cases described in Theorem
\ref{thm: MS classification} do arise. Namely we have: 

\begin{prop}\name{prop: new report U}
For any $k \geq d \geq 2$ and any real number field $\KK$, there are
$\R$-algebraic groups $H$ and $H'$ in $ M $, and $\mathcal{L}_1 = g_1 \Z^n
\in \ALN_n$, where $n= k \deg(\KK/\Q)$ and $g_1 \in M$,
such that the 
following hold: 
\begin{itemize}
\item
  $H'$ is defined over $\Q$, and  is 
  $\Q$-isogenous to $\Res_{\KK/\Q}(\mathbf{G})$, where $\mathbf{G}$ is
  $\KK$-isomorphic to $\SL_k$.
\item
  $H$ is either equal to $H'$ (linear case) or to $H' \ltimes \R^n$
  (affine case). 
\item
The orbit $H\mathcal{L}_1$ is closed and supports an $H$-homogeneous
probability measure $\nu$. The pushforward $\Psi_* \nu$ is
an RMS measure. 
\end{itemize}
The same statement is true with $d=2, \, n = 2k \deg(\KK/\Q),$ and
with $\mathbf{G} $ being $\KK$-isomorphic to $ \Sp_{2k}$ for 
some $k \geq 2$. 
\end{prop}

\begin{proof}
  The proof amounts to reversing the steps in the preceding
  discussion. For concreteness, we give it for $\mathbf{G}= 
\SL_k$. Let $D \df \deg(\KK/\Q)$, $n \df Dk$ and $\underline{G}
\df \SL_n(\R)$. The standard
action $\varphi$ of $\mathbf{G}_{\KK}$ on $\KK^k$ gives 
rise to a $\Q$-embedding $\Res_{\KK/\Q}(\varphi): \Res_{\KK/\Q}(\mathbf{G} )\to
\SL_{n}$. Let $\underline{H}_1$ denote the connected component of the identity in the group of
$\R$-points in $\Res_{\KK/\Q}(\mathbf{G})$.
Similarly to \eqref{eq: embedding linear} and \eqref{eq: embedding
  affine}, we refer to 
\eq{eq: embedding linear2}{g \mapsto 
\left(\begin{matrix} g & \mathbf{0}_{d,n-d}\\ \mathbf{0}_{n-d,d} &
      \mathrm{Id}_{n-d} \end{matrix} 
  \right)
}
as the top-left corner embedding of $\SL_d(\R)$ in $\underline{M}$. By
the explicit description of  
restriction of scalars described in \S \ref{subsec: number
  fields}, there is $\underline{g}_1 \in \underline{M}$ such that
$\underline{H} \df \underline{g}_1 \underline{H}_1 \underline{g}_1^{-1}$ contains
the top-left corner embedding of $\SL_d(\R)$ in $\underline{M}$, and 
up to scaling, $\underline{g}_1 \Z^n$ is the
geometric embedding of $\mathcal{O}^k$ as in \eqref{eq: def of lattice
from order}, where $\mathcal{O}$  is the
ring of integers in $\KK$. In particular, the orbit $\underline{H}
\underline{g}_1 \Z^n$ is a closed orbit supporting an
$\underline{H}$-homogeneous measure in $\LLN_n$.

Recall that there is an embedding of $\underline{M}$ in $M$ and of
$\LLN_n$ in $\ALN_n$ (respectively as the stabilizer of
the origin in the standard action on $\R^n$, and as the set of
lattices in the space of grids). We let $H'$ denote
the image of $\underline{H}_1$ under this embedding, and in the linear
case we set $H \df H'$ and let $H \LL_0 $ be the image of  $\underline{H}
\underline{g}_1 \Z^n$ under this embedding, and let $\nu$ be the
$H$-homogeneous measure on $H\LL_1$. Because the action of
$\SL_d(\R)$ is ergodic with respect to $\nu$, we can find $g_1$ so
that for $\LL_1 = g_1 \Z^n$ we have 
$\overline{\SL_d(\R) \LL_1} = H\LL_1 = H\LL_0$. It is not hard to
check that with these choices, the desired conclusions hold.
The proof in the affine case is similar, taking $H = H' \ltimes \R^n$ and 
$\underline{\pi}^{-1}(\underline{H} \underline{g}_1 \Z^n)$. 
  \end{proof}

\subsection{Preparations for the proof of Lemma
  \ref{lem:classification}}\label{subsec: preparations Lemma}
Recall that $\underline{\LL}_1= \underline \pi(\LL_1)$.
A vector space $V \subset\bR^n$ is called {\em $\LL_1$-rational} if
$V\cap\underline \LL_1$
is a lattice in $V$. In other words a subspace $V$ is $\LL_1$-rational
if it is of the form $\underline g_1W$ for some rational subspace $W \subset
\R^n$, i.e., a subspace spanned by vectors with rational entries. 
\begin{lem}
\name{lem:assref}
The following implications hold. 
\begin{itemize}
        \item[(a)] ${\bf (D)} \Rightarrow $ $\Vphys$ is
          not contained in a proper $\LL_1$-rational subspace. 
        \item[(b)] ${\bf (I)} \Rightarrow$ $\Vint$
          contains no nontrivial $\LL_1$-rational subspace. 
        \item[(c)] ${\bf (I)} \text{ and }  {\bf (D)} \Rightarrow {\bf
            (irred)}$.
\end{itemize}
\end{lem}
 Variants of statements (a) and (b) are given in \cite{Pleasants},
  but we give a complete proof for the convenience of the reader.
\begin{proof}
  We will prove all three statements by contradiction. 
        Suppose that (a) fails, so that there is a proper
        $\LL_1$-rational 
        subspace $W$ containing $\Vphys$. Let $W^\perp$ be an
        $\LL_1$-rational complement of $W$.  
        Since $W^\perp$ is $\LL_1$-rational, $\underline{\LL}_1$ is
        mapped to a lattice in $W^\perp$ under the projection $\R^n\to
        W^\perp$, and hence the projection of $\LL_1$ to $W^\perp$ is
        discrete.  
        On the other hand, $\R^n\to W^\perp$ factors through $\Vint$
        since $\Vphys\subset W$, and by $\textbf{(D)}$ the
        image of $\LL_1$ is dense in $\Vint$. Thus, the projection of
        $\LL_1$ is dense in $W^\perp$, a contradiction. 

      Now suppose that (b) fails, and
      $\Vint$ contains a nontrivial $\LL_1$-rational subspace $W$.
        Then $\Vint$, which is the kernel of the map $\R^n\to \Vphys$,
        contains $W\cap \underline{\LL}_1$, which by
        assumption is nontrivial. This contradicts \textbf{(I)}. 

        Now suppose  {\bf (D)} and {\bf (I)} hold but {\bf (irred)} fails, so 
     that there is a proper $\underline H_1$-invariant $\Q$-rational
     subspace $W$. From (b) we know that $\underline g_1W$
  is not contained in $\Vint$. Hence some $u\in \underline g_1W$
  can be written as 
  $$
u=u_p+u_i, \ u_p \in \Vphys \sm \{0\}, \, u_i \in \Vint . 
  $$
  Since $\SL_d(\R) \subset \underline H = \underline g_1 \underline H_1 \underline
  g_1^{-1}$, $\underline g_1W$ is also 
  $\SL_d(\R)$-invariant. Since $\SL_d(\R)$ acts 
  trivially on $\Vint$, for any $g \in \SL_d(\R)$ we have
  $$gu-u = gu_p - u_p \in
  \Vphys.$$
  We can find $g \in \SL_d(\R)$ such that $g u_p \neq u_p$, and hence
   $\underline g_1W\cap \Vphys$ is nontrivial. Since $\SL_d(\R)$ acts irreducibly
  on $\Vphys$, $\Vphys\subset \underline g_1 W$. This contradicts the conclusion
  of (a).
\end{proof}

\begin{thm}[Morris]\name{thm: from morris}
Let $n \geq d \geq 2$, and let $S$ be a connected real algebraic group which is
$\R$-almost simple, and contains the
image of $\SL_d(\R)$ under the top-left corner embedding (see \eqref{eq:
embedding linear2}). Then there are $k \geq d, \ell \geq d$ and $g \in \SL_n(\R)$
such that $gSg^{-1}$ is the image of
either $\SL_k(\R)$ or 
$\Symp_{2\ell}(\R)$ under the top-left 
corner embedding, and the latter
can only occur when $d=2$. 
\end{thm}

In this statement, by the `top-left corner embedding of
$\Symp_{2k}(\R)$', we mean the image under \eqref{eq:
embedding linear2}, that is, the
elements of $\SL_{2k}(\R)$ stabilizing a non-degenerate alternating
bilinear form on $\R^{2k}$. As is well-known, such a form can be taken
to be defined by 
\begin{equation*}
  \omega(\vec{x}_i, \vec{y}_j) = - \omega(\vec{y}_j, \vec{x}_i)=\delta_{ij} , \ \
  \omega(\vec{x}_i, \vec{x}_j) = \omega(\vec{y}_i, \vec{y}_j)=0
\end{equation*}
for some basis $\vec{x}_1, \ldots, \vec{x}_k, \vec{y}_1, \ldots,
\vec{y}_k$ of $\R^{2k}$. 

This result was proved by Dave Morris in 2014, in connection
with prior work of one of the authors and
Solomon. Namely, the result appeared in an initial ArXiV version
\cite{Barak_Yaar1} (in a slightly different form) but eventually did
not appear in  
the published version \cite{Barak_Yaar}.

For any $k\geq d$, we will refer to the 
image of $\SL_k(\R)$ under the top-left corner embedding in \eqref{eq:
  embedding linear2} (replacing $d$ with $k$ in that embedding) as the {\em top-left copy of
  $\SL_k(\R)$}. Clearly, with respect to the 
decomposition
\begin{equation}\label{eq: decomposition}
  \R^n = \R^k \oplus 
  \R^{n-k},
\end{equation}
the top-left copy of $\SL_k(\R)$ acts via its standard action on the first
summand, and the second summand is the set of vectors fixed by the
action.

Let $k$ be maximal, such that $S$
contains a conjugate (over $\SL_n(\R)$) of the
top-left copy of $\SL_k(\R)$. To make the ideas more transparent we
separate the proof into cases according to whether $k \geq 3$ (the
easier case) or
$k=2$. The proofs in these  cases are not independent -- readers
interested in the case $k=2$ are encouraged to first read the proof for $k \geq
3$. 

\begin{proof}[Proof in case $k\geq 3$.]
We recall the following result of Mostow \cite{Mostow}:
If $G_1 \subset \cdots \subset G_r \subset \SL_n(\R)$ are connected
reductive real algebraic groups, then there is $x \in \SL_n(\R)$ such
that $x^{-1} G_i x$ is self-adjoint for every~$i$. That is,
if $g \in x^{-1} G_i x$, then the transpose of~$g$ is also
in $x^{-1} G_i x$. 

Replacing $S$ by a conjugate, we may
assume that $S$ contains the top-left  embedding of $\SL_k(\R)$, which
we denote by $F$. By
Mostow's theorem, 
there is $x \in \SL_n(\R)$, such that $x^{-1} F x$ and
$x^{-1}Sx$ are self-adjoint. Let $V$ be the $(n - k)$-dimensional
subspace of $\R^n$ which is pointwise fixed by $F$. Since
$\SO_n(\R)$ acts transitively on the set of subspaces of any given
dimension, there is some $h \in \SO_n(\R)$, such that $xh(V) =
V$. After replacing $x$ with $xh$, we may assume that
$x^{-1}Fx$ fixes pointwise the 
second summand in the splitting \eqref{eq: decomposition}, and
$x^{-1}Fx$ and $x^{-1}Sx$ are self-adjoint (because this
property is not affected by conjugation by an element of
$\SO_n(\R)$). We conclude that $x^{-1}Fx = F$. Thus, we may assume that
$S$ is self-adjoint and contains $F$. We will assume 
that $S \neq F$ and derive a contradiction to the maximality of
$k$. Since $F \varsubsetneq S$ are connected, their Lie algebras $\Liew
f, \Liew s$ satisfy $\dim \Liew f < \dim \Liew s$. 

For $1 \le i, j \le n$, let $e_{i,j}$ be the elementary matrix
with~$1$ in the $(i,j)$ entry, and all other entries~$0$. Write 
\begin{equation}\label{eq: in analogy with}
\Liew{SL}_n(\R) = \Liew{f} \oplus \Liew Z \oplus X_1
        \oplus \cdots \oplus X_k \oplus Y_1 \oplus  \cdots \oplus Y_k
        ,
        \end{equation}
where
	\begin{itemize}
	\item $\Liew{SL}_n(\R)$ and $\Liew{f}$ are the Lie
          algebras of $\SL_n(\R)$ and $F$, respectively, 
	\item $\Liew Z$ is the subspace of $\Liew{SL}_n(\R)$ fixed
          pointwise by $\Ad(F)$, where $\Ad: \SL_n(\RR) \rightarrow \rm{Aut}(\mathfrak{sl}_n(\RR))$ is the adjoint representation, 
	\item $X_i$ is the linear span of $\{\, e_{i,j} : k+1 \le j \le n \,\}$,
	and
	\item $Y_j$ is the linear span of $\{\, e_{i,j} : k+1 \le i \le n \,\}$.
	\end{itemize}
Now we denote by $A$ the group of diagonal matrices in $F$ with
positive entries. We write an element $a \in A$ as 
\begin{equation}\label{eq: def A Morris}
  a =  \mathrm{diag}\bigl( a_1,a_2, \ldots, a_{k-1},
        (a_1 a_2\cdots a_{k-1})^{-1}, 1, \ldots, 1 \bigr), 
      \end{equation}
      and denote by $\chi_i$ the characters $a \mapsto a_i$, where
      $a_k \df (a_1 \cdots a_{k-1})^{-1}$. 
Since $k \geq 3$, the characters $\chi_i, \chi_i^{-1}$ are distinct,
for $i=1, \ldots, k$, and the subspaces       
$X_1,X_2,\ldots,X_k$ and $Y_1,Y_2,\ldots,Y_k$ are the corresponding
weight spaces, that is, 
	\begin{itemize}
	\item $X_i = \{x \in \Liew{SL}_n(\R) : \Ad (a)(x) 
          =  \chi_i(a) x \text{ for all } a \in A\}$, and
        \item
$Y_j = \{x \in \Liew{SL}_n(\R) : \Ad(a)(x) 
          =  \chi_j^{-1}(a) x \text{ for all } a \in A\}$. 
	\end{itemize}
We will use repeatedly the fact that if 
$\Liew{l}$ is an $\Ad (A)$-invariant subspace of
$\Liew{sl}_n(\R)$, and 
$v \in \Liew{l}$ has a nontrivial projection onto some weight
space, then this projection is contained in $\Liew{l}$. 

 Since $A \subset S$, $\Liew s$ is invariant
under $\Ad (A)$. 
Since $S$ is $\R$-almost simple and $\dim
\Liew f < \dim \Liew s$,  
$\Liew S$ cannot be contained in $\Liew{f}\oplus
\Liew Z$, and hence $\Liew s$ projects
nontrivially to some $X_i$ or $Y_j$. In 
fact, since $S$ is self-adjoint, it must project nontrivially to
both $X_i$ and~$Y_i$, for some~$i$. Since 
$X_i$ is a weight space of $\Ad(A)$, 
we find that $X_i \cap \Liew s$ is nontrivial.
Conjugating by an
element of $I_k \times \SO_{n - k}(\R)$, we may assume that $\Liew s$
contains the matrix $e_{i,k+1}$. Applying an appropriate element of
$\Ad(\SO_k(\R))$ shows that $e_{k,k+1} \in \Liew s$. Then,
since $S$ is self-adjoint, $\Liew s$ also contains
$e_{k+1,k}$. Therefore, $\Liew s$ contains the Lie subalgebra generated by
$\Liew{f}$, $e_{k,k+1}$, and~$e_{k+1,k}$, which is the Lie subalgebra
of $F'$, the top-left copy of $\SL_{k+1}(\R)$. Thus $S$ contains $F'$,
contradicting the maximality
of~$k$, and completing the proof in case $k \geq 3$.
\end{proof}
\begin{proof}[Proof in case $k=2$.]
In this case we also have $d =2$. 
Arguing as in the case $k \geq 3$ we may
assume that $S$ properly contains $F$, the top-left copy of $\SL_2(\R)$, and is
self-adjoint. 
Let $\ell$ be the maximal number so that $S$
contains a copy of $H \df F_1 \times \cdots \times F_\ell$, where each $F_r$
is isomorphic to $\SL_2(\R)$ and there is an $H$-invariant direct sum
decomposition $\R^n = V_1 \oplus \cdots \oplus V_\ell \oplus V_0$,
where the spaces 
$V_1, \ldots, V_\ell$ are two dimensional, and each $F_r$ acts linearly on $V_r$
and trivially on $\bigoplus_{s \neq
  r} V_s$. By assumption $\ell \geq 
1$, and there is a conjugation taking $H$ into  a top-left copy of
$\SL_{t}(\R)$, where $t =2\ell \geq 2$. We replace $H$ and $S$ by
their images under this conjugacy (retaining the same names $H$ and
$S$). By Mostow's theorem we can assume that $H$ and $S$ are both 
self-adjoint.  

Our first goal is to show that
\begin{equation}\label{claim one} 
S \text{  is also contained in the top-left copy of }
\SL_{t}(\R).
\end{equation}
Indeed, in analogy with \eqref{eq:
  in analogy with}, 
consider the decomposition
$$
\Liew{SL}_n(\R) = \Liew{l} \oplus \Liew Z \oplus \Liew M, \text{ where
} \Liew M =  X_1
\oplus \cdots \oplus X_{t} \oplus Y_1 \oplus  \cdots \oplus Y_{t},
$$
and 
\begin{itemize}
\item  $\Liew{l}$ is the Lie
          algebra of the top-left $\SL_{t}(\R)$, 
	\item $\Liew Z$ is the Lie algebra of the centralizer of the
          top-left $\SL_t(\R)$,
	\item $X_i$ is the linear span of $\{\, e_{i,j} : t+1 \le j \le n \,\}$,
	and
	\item $Y_j$ is the linear span of $\{\, e_{i,j} : t+1 \le i \le n \,\}$.
        \end{itemize}
With this notation,
our claim \eqref{claim one} is that $\Liew{s} \subset
\Liew{l}$.

        We note that
        \begin{equation}\label{eq: no nonzero}
          \Liew{s}
          \text{ does not contain a nonzero element in some
          } X_i \text{ or some } Y_i.
        \end{equation}
        Indeed, if $v \in (\Liew{s} \cap X_i ) \sm \{0\}$, we could
re-index to assume $i=1$, and conjugate by an element of $I_{t}
\times \SO_{n-t}(\RR)$ and rescale to assume $v= e_{1, t+1}$. Since $\Liew
s$ is self-adjoint, we also have $e_{t+1, 1} \in \Liew s$. Since
$\Liew{f}_1, \, e_{t+1,1}$ and $e_{1, t+1}$ generate a Lie 
algebra isomorphic to $\Liew{sl}_3(\R)$, this gives a contradiction to
the choice of $k$ and proves \eqref{eq: no nonzero}. 

If $\Liew{s} \not \subset \Liew{l}$, using that
        $\Liew{s}$ is simple and the Lie algebras $\Liew l, \, \Liew
        z$ commute, we see that the projection of $\Liew{s}$ 
        onto $\Liew M$ is nontrivial; indeed, if $\Liew{s} \subset
        \Liew{l} \oplus \Liew{z}$ then the kernel of the projection of
        $\Liew{s}$ to $\Liew{z}$ contains $\Liew{f}$ and by simplicity
        is equal to $\Liew{s}$. 

Let $A'$ be the intersection of $H$ with
the diagonal subgroup and let $\Liew{a}'$ be its Lie algebra. For each
odd index $i <t$, the
spaces $X_i \oplus Y_{i+1}$ and $X_{i+1} \oplus Y_i$  are weight
spaces for $\Ad(A')$, 
and hence there is some $i$ such that $\Liew{s} \cap (X_i \oplus
Y_{i+1} \cup X_{i+1} \oplus
Y_{i})$
contains a nonzero element $u$. Re-indexing, conjugating and rescaling
as in the proof of \eqref{eq: no nonzero}, we can assume
$u= e_{1, t+1} + \sum_{j \geq
  t+1} a_j e_{j,2},$
where the $a_j$ are not all zero. By a further conjugation by an
element of $I_t \times \SO_{n-t}(\RR)$ that fixes $e_{1,t+1}$, we can also
assume that $a_j =0$ for $j>t+2$, that is, we can write
$$
u= e_{1, t+1} + a e_{t+1, 2} + b e_{t+2, 2},  \ \ \text{ with } (a,b)
\neq (0,0). 
$$
Using brackets to denote the commutator  $[x,y]=xy-yx$ , we compute
\[
  \begin{split}
w &\df  [u, [u, e_{2,1}]] =   \left[ e_{1, t+1} + a e_{t+1,2} +
  be_{t+2,2}, -e_{2, t+1} +a e_{t+1,1} +be_{t+2,1} \right]
\\ & = a(e_{1,1} + e_{2,2} -2e_{t+1, t+1}) -2be_{t+2, t+1}
 \end{split}
\]
and
$$
[w,u] = 3a e_{1,t+1} - 3a^2e_{t+1,2}
    -3ab e_{t+2,2}, 
    $$
    so that
    $$
 6ae_{t+1,1}= [w,u] + 3au \in \Liew{s}.
$$
It follows from \eqref{eq: no nonzero} that $a=0$, and thus $\Liew{s}$
contains $\frac{-1}{2b} w = e_{t+2, t+1}$. Since $\Liew{s}$ is self-adjoint
it also contains $e_{t+1, t+2}$, and since these two vectors generate
a copy of $\Liew{sl}_2(\R)$ which is contained in $\Liew{z}$, and acts
on $\R^n$ 
by the standard two-dimensional representation, we have a
contradiction to the 
definition of $\ell$. This proves
\eqref{claim one}. 

Since $S$ properly contains $F$ we have $\ell>1$. We will now show
that $\Liew{s}$ is the Lie algebra $\Liew{sp}(2\ell, \R)$ of the
top-left  corner embedding 
of $\Symp_{2\ell}(\R)$. We will
begin with the case $\ell=2$ as it will make the argument more
transparent. That is, up to a conjugation in $\SL_n(\R)$, we want to
show that
\begin{equation}\label{eq: given by 1}
\Liew{s} = \Liew{h} \oplus \Liew{s}_{1,3} \oplus \Liew{s}_{1,4}\oplus
\Liew{s}_{2,3}\oplus \Liew{s}_{2,4},
\end{equation}
where $\Liew{h} \cong \Liew{sl}(2, \R) \oplus \Liew{sl}(2,\R)  \subset
\Liew{s}$ is the Lie algebra of $H$, and  
\begin{equation}\label{eq: given by 2}
  \begin{split}
\Liew{s}_{1,3} \df \spa (e_{1,3}- e_{4,2}) \ \ \  \ \ & \Liew{s}_{1,4} \df \spa
(e_{1,4}+e_{3,2}) \\ \Liew{s}_{2,3} \df \spa (e_{2,3}+ e_{4,1}) \ \ \
\ \ 
& \Liew{s}_{2,4} \df \spa (e_{2,4}- e_{3,1}).
\end{split}\end{equation}
To this end, let 
\begin{equation}\label{eq: contained in}
  \begin{split} \Liew{l}_{1,3} \df & \spa (e_{1,3}, e_{4,2}), \ \ \ \ 
\Liew{l}_{1, 4} \df \spa (e_{1,4}, e_{3,2}),  \\ 
\Liew{l}_{2,3} \df & \spa (e_{2,3}, e_{4,1}),  \ \ \ \ \Liew{l}_{2,4} \df
\spa (e_{2,4}, e_{3,1})
\end{split}\end{equation}
be the weight spaces for the action of $\Ad(A')$, which are
not in $\Liew{h}$. Let
$$
\Liew{s}'_{i,j} \df \Liew{l}_{i,j} \cap \Liew{s},
$$
where the indices $(i,j)$ range over $\{1,2\} \times \{3,4\}$. 
Our goal is to show that
\begin{equation}\label{eq: every pair}
  \text{ for every  } i , j,
  \ \  \ \Liew{s}'_{i,j}  = \Liew{s}_{i,j} 
 .
\end{equation}
We first show that
\begin{equation} \label{eq: every pair 2}
  \text{ for every } i,j, \ \ \  \dim ( \Liew{s}'_{i,j})=1.
\end{equation}
To this end, note that the $\ad$-action of the off-diagonal
elements of $\Liew{h}$ permutes the spaces $ \Liew{l}_{i,j}$ transitively. For
example,
$$
\Liew{l}_{1,3} = \left[e_{1,2}, \Liew{l}_{2,3} \right], \
\Liew{l}_{1,3} = \left[e_{4,3}, \Liew{l}_{1,4} \right],  
$$
and so on. 
Since $e_{1,2}, e_{2,1}, e_{3,4}, e_{4,3} \in \Liew{s}$, this
$\ad$-action also permutes the intersections $\Liew{s}'_{i,j} 
$, and thus they all have the same dimension. If this
dimension is 0 then $\Liew{s} = \Liew{h}$, contradicting the fact that
$\Liew{s}$ is simple, and if this dimension is $2$, then $\Liew{s} =
\Liew{sl}(4,\R)$, contradicting the definition of $k$. We have shown
\eqref{eq: every pair 2}. 

We now claim that 
\begin{equation}\label{eq: also we have}
  \Liew{s}'_{1,3}  \text{ is equal to  either } \Liew{s}_{1,3} =\spa
  (e_{1,3} -e_{4,2}) \text{ or } 
  \spa (e_{1,3} +e_{4,2}).
\end{equation}
To see this,  let $u = a e_{1,3} + b
e_{4,2} \in \Liew{s}'_{1,3} \sm \{0\}$. By \eqref{eq: no nonzero},
$a,b$ are both nonzero. Since $\Liew{s}$ 
is self-adjoint $v \df a 
e_{3,1} + b e_{2,4} \in \Liew{s}$ and hence $v \in
\Liew{s}'_{2,4}$. Also we have
\[
 w \df  [e_{2,1}, [e_{3,4}, u] ] 
  =[e_{2,1} ,
    -ae_{1,4}+be_{3,2}] = -ae_{2,4}-be_{3,1} 
    \in \Liew{s}.
  \]
  Since $w$ and $v$ are both nonzero elements of $\Liew{s}'_{2,4}$, by
  \eqref{eq: every pair 2} 
  they are scalar multiples of each 
  other and thus there is $c\neq 0$ so that $w = cv$. This forces $-a
  = cb$ and $-b = ca$ and so $c=\pm 1$, proving \eqref{eq: also we
    have}. 

  Using the
$\ad$-action as before we see that in order to obtain \eqref{eq: every
  pair}, it suffices to show that after a conjugation, we have
$\Liew{s}'_{1,3} = \Liew{s}_{1,3}$.  
  Suppose that $\Liew{s}'_{1,3}
  = \spa (e_{1,3} + e_{4,2})$. Then
  $$\Liew{s}'_{1,4} = \spa (
[e_{3,4},e_{1,3}+e_{4,2} ]) = \spa (e_{1,4}-e_{3,2}),$$
and we can apply a permutation
matrix swapping the indices 3,4 to obtain 
$$\Liew{s}'_{1,3}
= \spa (e_{1,3} - e_{4,2})  =  \Liew{s}_{1,3}.
$$
We have shown \eqref{eq: every
  pair}, completing the proof in case $\ell =2$. 

Note that for the case $\ell =2$
we only applied one conjugation, namely the conjugation swapping
the indices 3,4. Thus, by induction on $\ell$, we see that
after a conjugation, we have the following. For $i \in 
\{1, \ldots, \ell-1\}$, let $\SL_4^{(i)}(\R)$ be the copy of
$\SL_4(\R)$ embedded in $\SL_n(\R)$ in a $4 \times 4$ block
corresponding to indices $2i-1, 2i, 2i+1, 2i+2$. Let $H^{(i)} = F_i
\times F_{i+1} \subset \SL_4^{(i)}(\R)$ be the corresponding diagonal
copies of $\SL_2(\R)$, and let $\Liew{s}^{(i)}$ be the intersection of
$\Liew{s}$ with the Lie algebra of $\SL_4^{(i)}(\R)$. Then
$\Liew{s}^{(i)}$ is the obvious embedding of $\Liew{sp}(4,\R)$
(namely, the embedding given for $i=1$ by \eqref{eq: given by 1} and
\eqref{eq: given by 2}). 
The Lie algebras $\Liew{s}^{(i)}$ generate $\Liew{sp}(2\ell,\R)$
(namely, the Lie algebra of the top-left $\Symp_{2\ell}(\R)$). This
implies that $H$ contains
$\Symp_{2\ell}(\R)$. Since  $\Symp_{2\ell}(\R)$ is a maximal subgroup
among the connected Lie subgroups  
of $\SL_{2\ell}(\R)$ (see \cite{Karpelevich}), we must have that $S=\Symp_{2\ell}(\R)$. 
\end{proof}

\subsection{Proof of Lemma \ref{lem:classification}}\label{subsec:
  proof Lemma}
Since $\underline \pi$ is proper, we have 
$$\underline{H}
\underline{g}_1 \underline \Gamma = \underline \pi(H g_1 \Gamma) =\underline \pi \left(
  \overline{F g_1 \Gamma} \right) 
= \overline{ \SL_d(\R) \underline{g}_1 
  \underline \Gamma} .$$
Since $\underline{H}_1  = \underline{g}_1^{-1} \underline{H} \underline{g}_1$, by 
Theorem \ref{thm:  
  properties of L}, $\underline H_1$ is the connected component of the
identity in the
group of real points of a  $\bQ$-algebraic group $\mathbf{H}$. From
now on we replace $F$ with its image under $\pi$, i.e., denote $F
=\SL_d(\R)$. We also write 
$$
F' \df \underline g_1^{-1} F \underline g_1, \ \text{ so that} \ 
\overline{F' \underline \Gamma} = \underline H_1 \underline \Gamma.
$$

We need
to show that $\mathbf{H}$ admits the description given in the
statement. We divide the proof into steps.

\medskip

\textbf{Step 1: $\mathbf{H}$ is semisimple.} 
Let $\mathbf{U}$ be the radical of $\mathbf{H}$. 
By Theorem~\ref{thm:  
  properties of L}, it is defined over $\Q$ and unipotent,
$U=\mathbf{U}_\R^\circ$ is the unipotent radical of $\underline
H_1$, and $\mathbf{U}$ is connected (\cite[11.21]{Borel1}). 
Let $V^U$ be the subspace of $\R^n$ fixed by $U$. 
Since $\mathbf{U}_{\Q} \subset U$ is 
Zariski dense in $\mathbf{U}$ (see \cite[Cor. 18.3]{Borel1}), we have 
$$V^{U}=\{z\in \R^n: uz=z \text{ for all } u\in \mathbf{U}_{\Q} \}.$$
Thus $V^U$ is
defined over $\Q$.

Furthermore, since every unipotent
subgroup can be put in an upper triangular form, $V^U \neq \{0\},$ and
is a proper 
subspace of $\R^n$ unless $U$ is trivial. Since $U$ is normal in $\underline
H_1$, the space
$V^U$ is $\underline H_1$-invariant, and thus by assumption {\bf (irred)}, $V^U$
is not a proper subspace of $\R^n$. It follows that $U$ is trivial,
and hence $\underline H_1$ is semisimple. Therefore so is $\mathbf{H}$.

\medskip

For a group $M$ and normal subgroups $M_1, \ldots, M_k$, the {\em
  product} is the subgroup
$$ \prod M_i \df  \{m_1 \cdots m_k: m_i \in M_i, \ i=1,
\ldots, k\}.$$
Note that $\prod M_i$ is also
normal and does
not depend on the ordering of the $M_i$. 
Let $k_0$ be one of the fields $\Q$ or $\R$. 
Recall that an {\em almost direct product} is the image of a
direct product under a homomorphism with finite kernel (that is,
isogenous to a direct product).
A semisimple
$k_0$-group is an almost direct product of its $k_0$-almost simple
normal subgroups, 
and such a decomposition is unique up to permuting the $k_0$-almost
simple factors.

We write $\mathbf{H}$ in two ways: as an almost direct product of its
$\R$-almost simple factors 
$\mathbf{S}_i$, and as an almost direct product
of its $\Q$-almost simple factors $\mathbf{T}_j$, and let $S_i$ and $T_j$ denote
respectively the connected component of the identity in the group of 
$\R$-points of $\mathbf{S}_i$ and $\mathbf{T}_j$. Since every
$\mathbf{T}_j$ can be further decomposed into $\R$-almost simple
factors, and since these
decompositions are unique, the decomposition of
$\mathbf{H}$ into the $\mathbf{S}_i$ refines the decomposition of $\mathbf{H}$
into the $\mathbf{T}_j$. In other words, there is a partition of the
$\mathbf{S}_i$ into subsets such that each $\mathbf{T}_j$ is a
product of the $\mathbf{S}_i$ in one subset of the
partition. Then $\underline H_1$ is the product of
the $S_i$. 
For $h \in \underline{H}_1$,
we can write $h = h_1 \cdots 
h_t$, where $h_i \in S_i$, and if $h=h'_1 \cdots h'_t$ is another such
presentation, then for each $i$, $h'_i h_i^{-1}$ belongs to the finite
center of $\underline H_1$.

\medskip

{\bf Step 2: $F'$ is contained in one of the $S_i$, and $\mathbf{H}$
  is $\Q$-almost simple.} 
The second assertion follows from the first one. Indeed, by
re-indexing, let $S_1$ 
and $T_1$ denote respectively the connected
component of the identity in the real points of the $\R$- and
$\Q$-simple factors containing $F'$. 
Then $S_1 \subset T_1$ and $T_1$ does not 
properly contain the real points of any $\Q$-subgroup containing
$S_1$, and by the last assertion of 
Theorem \ref{thm: properties of L} we have that $\underline H_1  = T_1.$

Turning to the first assertion, let ${\rm Z}(\underline H_1)$ denote the center of $\underline H_1$, 
for each $i$ let $S'_i$ be the 
quotient group $\underline H_1 /\left( {\rm Z}(\underline H_1 ) \cdot \prod_{j \neq i} S_{j} \right)$, and
let $F'_i$ denote the image of the projection of $F'$ to $S'_i$. 
Let
$$
\underline{H}_2 \df \prod_{i \in \mathcal{I}} S_i, \ \ \text{ where }
\ \mathcal{I} \df \{ i : F'_i \text{
  is nontrivial} \}. 
$$
Note that $i_0 \in \mathcal{I}$ if and only if for any subset
$\mathcal{F}' \subset F'$ which generates a dense subgroup, there is
$f' \in \mathcal{F}'$ which can 
be written as a product of elements $f'_i$ in $S_i$, where $f'_{i_0}$ is not
central in $S_{i_0}$. 
Clearly $F' \subset \underline{H}_2$,  
  and our goal is to show that $\underline{H}_2$ is
equal to one of the $S_i$, or in other words that $\# \,
\mathcal{I}=1$. Also, for $i \in \mathcal{I}$, $F'_i$ is isogenous to
$\SL_d(\R)$. 

Recall that a representation of
a group $H$ on a vector space $V$ is
{\em isotypic} if $V$ is the direct sum of $k \in \N$ isomorphic
irreducible representations 
for $H$, where $k$ is referred to as the {\em
  multiplicity}. We will also use the term {\em $H$-isotypic}, if we want to make the
dependence on $H$ explicit. A linear representation of a
semisimple group has a 
unique presentation as a direct sum of isotypic representations (up to
permuting factors). Let $\Vphyss \df \underline
g_1^{-1}(\Vphys)$ and $\Vintt \df \underline g_1^{-1} (\Vint)$. Then
the decomposition $\R^n = \Vphyss
\oplus \Vintt$, is the decomposition of $\R^n$ into $F'$-isotypic
representations, and the action of $F'$ on $\Vphyss$ is
irreducible. In particular, the multiplicity of the  representation on
$\Vphyss$ is equal to one. 

Let $ V_1 \oplus \cdots \oplus V_t$ be a
decomposition of $\R^n$ into $\underline H_2$-isotypic
representations. Since $F' \subset \underline H_2$, each $V_\ell$ is
$F'$-invariant, and decomposes further into isotypic representations
for $F'$. 
Since $\Vphyss$ is an isotypical component of $F'$ of multiplicity one, 
$\Vphyss$ is contained in one of the $V_\ell$.
By renumbering we can assume
  $\Vphyss \subset V_1$.
  Since $F'$ acts on $\Vphyss$ irreducibly, the
  action of $\underline{H}_2$ on $V_1$ is irreducible, and the
  $\underline{H}_2$-isotypic component 
    associated to $V_1$ has multiplicity one. Since $F'$ acts 
  trivially on $\Vintt$, which 
  is a complementary subspace to $\Vphyss$, the action 
  of $F'$ on each $V_\ell$ is trivial for $\ell =  2, \ldots, t$, that is,
  \eq{eq: F contained}{
F' \subset \bigcap_{\ell = 2}^t \ker \left(\underline H_2|_{V_\ell}\right).
}
The right-hand side of \equ{eq: F contained} is a  normal subgroup of
  $\underline{H}_2$, and thus a product $\prod_{ i \in \mathcal{J}}
  S_i$ for some $\mathcal{J} \subset \mathcal{I}$.
  By the 
  assumption that $F'_i$ is nontrivial for each $i \in \mathcal{I}$,
  we must have that $\mathcal{J} = \mathcal{I}$, that is, the group on
  the right-hand side of \equ{eq: F contained} must coincide with
  $\underline{H}_2$. This means that for 
  $\ell \geq 2$, the $V_\ell$ are trivial representations for
  $\underline{H}_2$, and hence of $S_i$ for each $ i \in
  \mathcal{I}$.

  Let $\mathcal{F}'$ denote the elements of $F'$ whose
eigenvalues on $\Vphyss$ are all real, distinct from each other, and
not equal to 1. Since these conditions are invariant under
conjugation and $F'$ is simple, $\mathcal{F}'$ generates a dense subgroup of $F'$.
Write $f'$ as a product of elements $f'_i$, where $f'_i \in S_i$.
Then the elements $f'_i$ commute with each other and with $f'$. Thus
each $f'_i$ fixes the eigenspaces for  $f'$ and hence each $f'_i$
preserves the eigenspace decomposition of the action of $f'$ on
$\R^n$. In particular, $f'_i$ preserves
$\Vphyss$ for each $i \in 
\mathcal{I}$.

Re-indexing if
necessary we can assume that $1 \in \mathcal{I}$, and suppose by
contradiction that there is $i_0 \in \mathcal{I} \sm \{1\}$. There is $f'
\in \mathcal{F}'$ such that, when writing $f'$ as a product of elements $f'_i \in
S_i$, $f'_1$  acts on $\Vphyss$ with infinite order (this property
does not depend on the presentation 
of $f$ as a product of the $f'_i$). 
Then the action of $f'_1$ on 
$\Vphyss$ preserves an eigenspace $V'$, with $d' \df \dim V' < d = \dim
\Vphyss$. 
Since the action of $S_{i_0}$ commutes with the action of $f'_1$, the space
$V'$ is preserved by $S_{i_0}$, and hence by $f'_{i_0}$. The group
generated by all such elements $f'_{i_0}$ is isogenous to $F'_{i_0}$ and hence
to $\SL_d(\R)$. Thus, it has no nontrivial representations on any
$d'$-dimensional real vector space, for 
$d'<d$. This implies that the action of $S_{i_0}$ on $V'$ has an infinite
kernel, but since $S_{i_0}$ is simple, the action of $S_{i_0}$ on 
$V'$ must also be trivial.

So the
space
$$V'' \df \spa \, S_1 (V') \subset \spa \, S_1 (\Vphyss) \subset V_1$$
is acted on trivially by $S_{i_0}$ for
any ${i_0} \in \mathcal{I} \sm \{1\}$. In particular, $V''$ is $\underline{H}_2$-invariant. By
the irreducibility of the $\underline{H}_2$-action on $V_1$, this means that $V_1 =
V''$, and therefore $S_{i_0}$ acts trivially on $V_1$.
It follows that 
$F'_{i_0}$ acts trivially on $V_1$ for each $i_0  \in \mathcal{I} \sm \{1\}$. Since 
$S_{i_0}$ acts trivially on $V_\ell$ for all $ i_0 \in \mathcal{I}$ and all
$\ell \geq 2$, we get that in any decomposition of $f' \in F'$, all
the elements $f'_i$ for $i \geq 2$ act trivially on $\R^n$. That is,
$\mathcal{I} = \{1\}$. 

\medskip

\textbf{Step 3: Restriction of scalars, in explicit form.}
Since $\mathbf{H}$ is $\Q$-almost simple,  it is
obtained by restriction of scalars from an absolutely almost simple algebraic group defined
over a number field $\KK$ -- see
\cite[6.21]{BT_reductif} for a proof. We will reprove this result in our
setup, obtaining more information about the embedding of $\underline H_1$ in
$\SL_n(\R)$.


Using Step 2 and re-indexing, let $S_1 = (\mathbf{S}_1)^{\circ}_{\R}$ be the
connected component of the identity in the
$\R$-almost simple group containing $F'$, and set $\mathbf{G}\df
\mathbf{S}_1, \, G \df S_1$. It follows from \cite[\S
2.15b]{BT_reductif} that $\mathbf{G}$ is Zariski connected, which
implies via \cite[Cor. 18.3]{Borel1} that $G$ is Zariski dense in
$\mathbf{G}$. 
From Theorem 
\ref{thm: from morris}, we only have two possibilities for
$G$, and its Zariski closure is a conjugate
of either $\SL_k$ or $\Sp_{2\ell}$. Hence $\mathbf{G}_\R$ is a
conjugate of $\SL_k(\R)$ or $\Sp_{2\ell}(\R)$. 
In particular, we have that $\mathbf{G}$ is actually
$\CC$-almost simple. 
Since $\mathbf{H}$
is defined over $\Q$, the
$\CC$-almost simple factors of $\mathbf{H}$ are defined over a finite
extension of $\Q$; this is well-known (see e.g. \cite[\S
2.15b]{BT_reductif}) but we were unable to find a suitable reference,
so  we sketch the argument.
The group $\mathbf{H}$ has a maximal torus which is defined over $\Q$
and split over a finite extension $\mathbb{L}$ of $\Q$ by \cite[\S8, \S18]{Borel1}.
For each root $\alpha$, the group $G_\alpha$, which is the centralizer
of the connected component of the identity in $\ker \alpha$, is defined
over  $\mathbb{L}$ (see \cite[Proof of Thm. 18.7]{Borel1}). The groups
$G_\alpha$ generate $\mathbf{H}$ \cite[\S 14]{Borel1} and each
$\CC$-almost simple factor either contains $G_\alpha$, or intersects
it trivially. Thus, any $\CC$-almost simple factor $\mathbf{S}$ can be described as
the elements commuting with all the $G_\alpha$ not contained in
$\mathbf{S}$. In particular, the $\CC$-almost simple factors are
defined over $\mathbb{L}$. 


Replacing $\mathbb{L}$ if necessary with its Galois extension, suppose that
$\mathbb{L}$ is the smallest Galois extension of $\Q$ such that all
$\CC$-almost simple factors of $\mathbf{H}$ are defined over $\mathbb{L}$. 
Let $\on{Gal}(\mathbb{L}/\Q)$ denote the Galois group of
$\mathbb{L}$, which we can think of explicitly as the group of field
automorphisms of $\mathbb{L}$. 
If $\mathbf{V} \subset \mathbb{A}^n$ is
an affine variety defined 
over $\mathbb{L}$ then for any $\sigma \in \on{Gal}(\mathbb{L}/\Q)$ there
is a new affine variety, which we will denote by ${}^\sigma\mathbf{V}$, 
obtained by acting on the coefficients of the defining polynomial
equations, and $\sigma$ acts on the points of  $\mathbb{L}^n$ by acting
 separately on each component. 
The assignments $\mathbf{V} \mapsto {}^\sigma\mathbf{V}$ and $\sigma:
\mathbb{L} \to \mathbb{L}$ are compatible in the
sense that for $x \in \mathbb{L}^n$, $x \in \mathbf{V}_{\mathbb{L}}$ if
and only if $\sigma(x) \in {}^\sigma\mathbf{V}_{\mathbb{L}}$. Moreover,
if $\mathbf{V}$ is defined over $\mathbb{L}$, then it is defined 
over $\Q$ if and only if ${}^\sigma\mathbf{V}=
\mathbf{V}$ for every $\sigma \in \on{Gal}(\mathbb{L}/\Q)$; this
follows from the more general fact (see
\cite[\S AG12-\S AG14]{Borel1}), that
if $\mathbb{L}'$ is
a number field then $\mathbf{V}$ is 
defined over $\mathbb{L}'$ if and only if for any $\sigma \in
\on{Gal}(\bar{\Q}/\Q)$ such that $\sigma|_{\mathbb{L}'}=\mathrm{Id}$ we
have ${}^\sigma\mathbf{V}= \mathbf{V}$, where $\bar \Q$ denotes the
algebraic closure of $\Q$.

Let $D$ denote the number of $\CC$-almost simple factors of
$\mathbf{H}$, or equivalently, the number of $\mathbb{L}$-almost simple factors of
$\mathbf{H}$. 
The action of $\on{Gal}(\mathbb{L}/\Q)$ permutes these 
factors, and this permutation action is transitive since $\mathbf{H}$
is $\Q$-almost simple. 
Thus, the subgroup
$$\Delta \df \{\sigma \in \on{Gal}(\mathbb{L}/\Q):{}^{\sigma} \mathbf{G}=
\mathbf{G} \}$$
is of index $D$ in $\on{Gal}(\mathbb{L}/\Q)$,
and the 
$\CC$-almost simple factors are the (distinct) images of $\mathbf{G}$ by elements
$\sigma_1, \ldots, \sigma_D \in \on{Gal}(\mathbb{L}/\Q)$, where the
$\sigma_i$ are coset representatives of
$\on{Gal}(\mathbb{L}/\Q)/\Delta$.

Let $$\KK \df \{x \in \mathbb{L}: \forall \sigma
\in \Delta, \, \sigma(x) =x\}.$$
Complex conjugation $z \mapsto \bar z $ induces an
automorphism of $\mathbb{L}$ belonging to $\Delta$ since $\mathbf{G}$
is defined over $\R$, hence we see that $\KK \subset \R$. 
By the Galois
correspondence,  $\deg(\KK/\Q)=D$ and 
$$
\Delta = \{\sigma \in \on{Gal}(\mathbb{L}/\Q) : \text{ for all } x \in
\KK, \ \sigma (x)=x \}. 
$$
We claim that 
$\mathbf{G}$ is defined over $\KK$, and 
$\mathbf{G}$ is not defined over any proper subfield
of $\KK$.
 Indeed, if $\sigma \in \on{Gal}(\bar \Q/\Q)$ satisfies $\sigma|_{\KK}
 = \mathrm{Id}$, then $\sigma|_{\mathbb{L}} \in \Delta$ and hence
 ${}^{\sigma}\mathbf{G} = \mathbf{G}$. Furthermore, if $\mathbf{G}$ were
 defined over a proper subfield $\KK' \varsubsetneq \KK$,  then its
 stability group $\Delta'$ would be 
 of index $D' < D$ and therefore the collection $\{{}^{\sigma}\mathbf{G} :
 \sigma \in \on{Gal}(\mathbb{L}/\Q) \}$ would have cardinality $D'$.


We will show that $\mathbf{H}$ is isomorphic (as a $\Q$-algebraic
group) to $\Res_{\KK/\Q}(\mathbf{G})$. Moreover, we will show that the
given inclusion $\mathbf{H} \hookrightarrow \SL_n$ is, up to a
conjugation over $\SL_n(\R \cap \bar \Q)$, the matrix presentation described
in \S \ref{subsec: number fields}.  
By Theorem~\ref{thm: from morris} $G$ is, up to a conjugation in
$\SL_n(\R)$, either the top-left copy of 
$\SL_k(\R)$ or the top-left copy of $\Symp_{2k}(\R)$ for some $k \geq
2$ (and the latter can only arise when $d=2$). 
In the remainder of the proof we
will refer to these two cases as the {\em $\SL_k$ case} and the {\em $\Sp_{2k}$ case.}

We know that $G$ is conjugate over $\SL_n(\R)$ to the top-left
copy of $\SL_k(\R)$ (in the $\SL_k$ case) or $\Sp_{2k}(\R)$ (in the
$\Sp_{2k}$ case). Therefore there is
a $G$-invariant subspace
$V \subset \R^n$, of dimension $k$ (in the $\SL_k$ case) and $2k$ (in
the $\Sp_{2k}$ case) and a
complementary subspace $V_0$ such that $\R^n = V \oplus V_0$, the action
of $G$ on $V$ is irreducible, and $V_0$ is the subspace of $G$-fixed
vectors in $\R^n$. We claim that we can recover $V$ explicitly as
\begin{equation}\label{eq: can recover}
V = \spa \left \{ gx -x : g \in G, \, x \in \R^n \right \}.
\end{equation}
Indeed, denote the RHS of \eqref{eq: can recover} by $W$. We clearly have $W
\subset V$, and for the reverse inclusion, it is enough to show that $W$
is $G$-invariant. To see this, let $g_0, g \in G$ and $x \in
\R^n$. Then
$$g_0(gx-x) = g_0 g g_0^{-1}g_0 x - g_0x = g'x' - x',$$
where $ g' \df g_0 g g_0^{-1} $ and $ x' \df g_0x.$
This shows that the generators of $W$ are mapped to $W$ by any $g_0
\in G$. 

From
\eqref{eq: 
  can recover} and since $\mathbf{G}$ is defined over $\KK \subset
\R$, we deduce that $V = \mathbf{V}_{\R}$ for a subspace $\mathbf{V}
\subset \mathbb{A}^n$
defined over $\KK$. 
Clearly $V_0=(\mathbf{V}_0)_\R$ for a subapce $\mathbf{V}_0$ which is also defined over $\KK$. Arguing as in \eqref{eq: can
  recover}, but using $F'$ in place of $G$ 
and $\Vphyss$ in place of $V,$ we have $\Vphyss = \spa \{f'x-x: f' \in F',\
x \in \R^n\}$, and therefore $\Vphyss \subset V$.

We can think of
$\mathbf{V}_{\bar \Q}$ as a
$\bar{\Q}$-linear subspace of $\bar \Q^n$, and can discuss the
action of $\on{Gal}(\bar{\Q}/\Q)$ as before. We have
that $(\mathbf{G}_i)_{\bar \Q}$ 
preserves the decomposition $\bar \Q^n ={}^{\sigma_i}\mathbf{V}_{\bar \Q}
\oplus \left({}^{\sigma_i}\mathbf{V}_0\right)_{\bar \Q}$.  We claim
that
\begin{equation}\label{eq: a direct sum}
  \bar \Q^n = \bigoplus_{i=1}^D {}^{\sigma_i}\mathbf{V}_{\bar \Q}. 
\end{equation}
To see this, let $\mathbf{W}$ denote the vector
subspace of $\mathbb{A}^n$ spanned by $\bigcup_i {}^{\sigma_i}\mathbf{V}
$. Since it is 
$\on{Gal}(\bar{\Q}/\Q)$-invariant, it is  defined over $\Q$. Since $\Vphyss =
\underline{g}_1^{-1}\Vphys$ and $\Z^n =\underline{g}_1^{-1}  \LL_1$,
Lemma \ref{lem:assref} implies that $\Vphyss$ is not contained in any proper
rational subspace of $\R^n$. This implies that $\mathbf{W}_{\R} =
\R^n$ and thus $\mathbf{W} = \mathbb{A}^n$. The groups $\mathbf{G}_i$
commute, and  ${}^{\sigma_i}\mathbf{V}$ is a
$\mathbf{G}_i$-isotypic component of multiplicity one. 
For each pair of distinct $i, j$, each $g \in \mathbf{G}_i$ defines an intertwining
operator for the action of $\mathbf{G}_j$,
and thus by Schur's lemma
(see e.g. \cite[Cor. 4.9]{Knapp}),
the action of $\mathbf{G}_i$ on
${}^{\sigma_j}\mathbf{V}$ factors through an abelian group. Since
$\mathbf{G}_i$ is simple, this means that 
each $\mathbf{G}_i$ acts trivially on
${}^{\sigma_j}\mathbf{V}$ for $j \neq i$.
In particular,
${}^{\sigma_i}\mathbf{V} \cap \sum_{j\neq i} {}^{\sigma_j}\mathbf{V} = \{0\}$,
and we have shown \eqref{eq: a direct sum}.

It follows from \eqref{eq: a direct sum} that $\R^n$ is the
space of $\R$-points of $\Res_{\KK/\Q}(\mathbf{V})$. Write $D = r+2s$
as in \S \ref{subsec: number fields}. Since $\dim
{}^{\sigma_i}\mathbf{V} = \dim {}^{\sigma_j}\mathbf{V}$ for every $i \neq
j$, we have that $\underline{H}_1$ is realized explicitly in $r+s$
blocks. For real embeddings 
$\sigma_i, \, i=1, \ldots, r$ we have that the dimension (over $\R$)
of ${}^{\sigma_i}\mathbf{V}_{\R} $ is $k $
(in the $\SL_k$ case) and $2k$ (in the $\Sp_{2k}$ case), and
for $\sigma_{r+j}, \ j=1, \ldots, 
s$ which are non-conjugate complex embeddings of $\KK$ we have that
the dimension (over $\R$) 
of ${}^{\sigma_{r+j}}\mathbf{V}_{\CC} $ is $2k $ (in the $\SL_k$ case) and
$4k$ (in the $\Sp_{2k}$ case). Putting
this together we get that $n = Dk$ (in the $\SL_k$ case) and  $n=2Dk$ (in
the $\Sp_{2k}$ case), and the embedding of $\underline{H}_1$ in
$\SL_n(\R)$ is the one 
given in \eqref{eq: shape of representation}, where $\varphi: \SL_k
\to \SL_k$ is the identity map (in the $\SL_k$ case), and
$\varphi: \Sp_{2k} \to \SL_{2k}$ is the natural embedding (in the 
$\Sp_{2k}$ case).  In particular, we have proved that $\mathbf{H} =
\Res_{\KK/\Q}(\mathbf{G})$, with the explicit form of restriction of
scalars given in \S \ref{subsec: number fields}. 

\medskip

{\bf Step 4: $\mathbf{G}$ as a $\KK$-group.} 
It remains to identify the $\KK$-isomorphism type of $\mathbf{G}$.
We
proved in Step 3 that $\KK \subset \R$, the
decomposition $\R^n = V 
\oplus V_0$ into $G$-invariant subspaces is defined
over $\KK$, and there is a conjugacy over $\SL_n(\R)$
sending $G$ 
%
to the top-left corner embedding of $\SL_k(\R)$ or of 
$\Sp_{2k}(\R)$ (as defined after the statement of Theorem \ref{thm: from
  morris}).
We now show that as a $\KK$-group, $\mathbf{G}$
is $\KK$-isomorphic to either $\SL_k$ or $\Symp_{2k}$.

Consider first the $\SL_k$-case.
Let $\mathbf{W}\oplus\mathbf{W}_0=\CC^{k}\oplus\CC^{n-k}$ (whose real
points we used in equation~\eqref{eq: decomposition}), 
and note that both subspaces are defined over $\KK$.
Since $\mathbf{V}, \mathbf{V}_0$ are $\KK$-subspaces,
we can find $g\in \SL_n(\KK)$,
such that $g\mathbf{V}= \mathbf{W}$, $g\mathbf{V_0}=\mathbf{W_0}$,
and hence, $\mathbf{G}' = g\mathbf{G} g^{-1}$ is contained in the top-left corner embedding of $\SL_k$.
In particular, 
the groups $\mathbf{G}$ and 
$\mathbf{G}'$ are $\KK$-isomorphic, and $\mathbf{G}'_{\R}$ is
$\R$-isomorphic to the top-left $\SL_k(\R)$. Let $\mathbf{G}''
= \SL(\mathbf{W})=\SL_k$ (top-left corner embedding) considered as a $\KK$-group. Then $\mathbf{G}''_{\R}$ is
also $\R$-isomorphic to $\SL_k(\R)$, and thus $\mathbf{G}'$ and $ 
\mathbf{G}''$ have the same dimension (as algebraic varieties). Since
$\mathbf{G}'_{\KK} = g \mathbf{G}_{\KK} g^{-1} 
\subset \mathbf{G}''$, there is a $\KK$-embedding $\mathbf{G} \hookrightarrow
\mathbf{G}''$, and since these groups have the same dimension and are
Zariski connected,
$\mathbf{G}$ and $\mathbf{G}''$ are $\KK$-isomorphic. 

Now consider the $\Symp_{2k}$ case.
We have shown that $\dim V = 2k$ is even, and we adjust the definitions $\mathbf{W}\oplus\mathbf{W}_0=\CC^{2k}\oplus\CC^{n-2k}$.
We let again $g \in \SL_n(\KK)$ be
the conjugating element sending $\mathbf{G}$ to $\mathbf{G}' =
g\mathbf{G} g^{-1}\subset \SL(\mathbf{W})$.
$\mathbf{G}'_\R$ is $\R$-isomorphic to $\Symp_{2k}(\R)$, that is,
there is a nondegenerate alternating bilinear form $\omega$ on
$\mathbf{W}_\R$
such that $\mathbf{G}'_{\R}$ is the group of all $\R$-linear transformations of
$\mathbf{W}$ preserving $\omega$.
Note that $\omega$ is $\R$-bilinear and takes
values in $\R$.
We claim that there is a form
$\omega'$ which is defined over $\KK$ on $\mathbf{W}$ (that is, takes values in
$\KK$ when evaluated on elements of 
$\mathbf{W}_{\KK}$),
so that $\mathbf{G}'_{\R}$ is contained in the group of
$\R$-linear transformations of $\mathbf{W}$ preserving $\omega'$.
Once the claim is proved, we will have that there is a
$\KK$-embedding $\mathbf{G} \hookrightarrow \Symp(\mathbf{W}, \omega')$ (the
group of linear transformations of $\mathbf{W}$ preserving $\omega'$) which
will be an isomorphism by dimension considerations as in the preceding
case, thus proving that $\mathbf{G}$ is $\KK$-isomorphic to
$\Symp(\mathbf{W}, \omega') \cong \Symp_{2k}$. 

To prove the claim, consider the collection $\bigwedge^2(\mathbf{W}^*)$ of
alternating bilinear forms on $\mathbf{W}$.
This collection is a linear space, and the
nondegenerate forms form a Zariski open subset (since nondegeneracy
is equivalent to the non-vanishing of the determinant of the Gram
matrix of the form).
Since $\mathbf{G}'$ is a $\KK$-group, the
subspace $\bigwedge^2(\mathbf{W}^*)^{\mathbf{G}'}$ of
$\mathbf{G}'$-invariant forms is a $\KK$-subspace, which is nonempty
since its collection of $\R$-points contains $\omega$. Since
$\KK$-points are Zariski dense in 
$\KK$-subspaces, we find that there are nondegenerate symplectic
$\KK$-forms which are $\mathbf{G}'$-invariant.

Finally, the
proof of Theorem~\ref{thm: from morris} shows that in the symplectic
case, the space $g V'_{\text{phys}}\cong \R^2$ is spanned by two vectors
$\vec{x}, \vec{y}$ satisfying $\omega(\vec{x}, \vec{y})=1$; that is,
$g V'_{\text{phys}}$ is a symplectic subspace for $\omega$.
(We recall at this point that $V'_{\text{phys}}=g_1^{-1}\Vphys \subset V$,
$g$ is the conjugation mapping $V$ to $W$,
and $\omega$ is the real symplectic form on $W$ induced by the isomorphism of
$G'_\R\simeq \Sp_{2k}(\R)$ from Theorem~\ref{thm: from morris}.)

Write $\omega$ as a
linear combination of forms $\omega'$ which are defined over $\KK$ and 
$\mathbf{G}'$-invariant.
Since $\omega(\vec{x}, \vec{y}) \neq 0$,
there has to be some $\omega' \in \left(\bigwedge^2(\mathbf{W}^*)^{\mathbf{G}'}\right)_\KK$ for
which $\omega'(\vec{x}, \vec{y}) \neq 0$.
This shows that $V'_{\text{phys}}$ is a symplectic subspace of $V$ under the form induced by $\omega'$.
\qed

\begin{remark}
In the symplectic case, Step 4 also shows that there is a symplectic
form on the entire space $\R^n$ that is preserved by the entire group
$\underline{H}_1$. Indeed, the form $\omega'$, which is symplectic and
defined over $\KK$, can be `pushed' using the field embeddings
$\sigma_i$ to induce symplectic forms on the spaces
${}^{\sigma_i}\mathbf{V}$. We will not be using this fact and we leave
the details to the reader. 
  \end{remark}

\section{An intrinsic description of the measures arising via
  $\Psi_*$}
The following result shows that all RMS measures arise via the map
$\Psi_*$. For a given constant $c>0$, we denote by $\rho_c: \Cl(\R^d) \to
\Cl(\R^d)$  the map induced by the dilation by $c$, that is,
$
\rho_c(F) =\{ cx: x \in F\}.
$

\begin{thm}\name{thm: MS surjective}
Let $F$ be as in \equ{eq: def G0} and
embedded in $G$ via the top-left corner embedding. For any ergodic $F$-invariant Borel probability measure $\mu$ on
$\Cl(\R^d)$ which assigns full measure to irreducible cut-and-project
sets, there is an irreducible cut-and-project construction with $\R^n =
\Vphys \oplus \Vint, \piphys, \piint,  W$ and with $\Psi$ as
in \equ{eq: def Psi}, a constant $c>0$, and an
$F$-invariant ergodic homogeneous 
measure $\bar \mu$ on $\ALN_n$, such that
$\mu = \rho_{c*}  \Psi_* \bar \mu$. For $\mu$-a.e.\  $\Lambda$ we have 
\eq{eq: for scaling}{
c = \left( \frac{\vol(W)}{D(\Lambda)} \right)^{\frac{1}{n}}, 
}
where $D(\Lambda)$ is the density of $\Lambda$ as defined in \equ{eq: density exists}. 
\end{thm}

We will split the proof into the linear and affine case. 

  \begin{proof}[Proof of Theorem \ref{thm: MS surjective}, affine case]
  Suppose $\mu$ is $\ASL_d(\R)$-invariant and $F = \ASL_d(\R)$, and
  let $\{g_t\}$ be a 
one-parameter diagonalizable subgroup of $\SL_d(\R) \subset F$. By the Mautner
phenomenon (see \cite{ew}), the action of $\{g_t\}$ on $\left(
  \Cl(\R^d), \mu \right)$ is ergodic. Thus, by the Birkhoff pointwise
ergodic theorem, there is a subset  $X_0 \subset \Cl(\R^d)$ of full
$\mu$-measure such that for all $\Lambda \in X_0$ we have
$$
\frac{1}{T} \int_0^T (g_t)_* \delta_\Lambda \, dt \to_{T \to \infty} \mu.
$$
Since the function $\Lambda \mapsto D(\Lambda)$ is measurable and
invariant, we can further assume that the value 
of $D(\Lambda)$ is the same for each $\Lambda \in X_0$.

Let $U, \Omega, m_U$ be as in Theorem \ref{thm: Shah}. Then by
Fubini's theorem, and since $\mu$ is $U$-invariant, we have
\[\begin{split}
    1 & = \mu(X_0) = \frac{1}{m_U(\Omega)}\int_{\Omega}
    \mu(u^{-1}X_0) \, dm_U(u)\\
& =\int \left[ \frac{1}{m_U(\Omega)}\int_{\Omega}
\mathbf{1}_{X_0}(u \Lambda) \, dm_U(u)  \right]\, d\mu(\Lambda),
\end{split}
\]
where $\mathbf{1}_{X_0}$ is the indicator function of $X_0$. 
Thus the inner integral on the RHS is equal to one on a subset of full
measure; i.e., there is 
$X_1 \subset \Cl(\R^d)$ of full measure such that 
for every $\Lambda \in X_1$ we have
$
u\Lambda \in X_0 \text{ for $m_U$-a.e.\  } u \in \Omega.  
$
This implies that for $\Lambda \in X_1 $ we have 
\eq{eq: nice convergence}{
\frac{1}{T} \int_0^T \int_{\Omega} (g_t u)_* \delta_{\Lambda} \,
dm_U(u) \, dt \to_{T
  \to \infty} \mu. 
  }
Let $\Lambda
\in X_1$ be an irreducible cut-and-project set, that is, $\Lambda =
\Psi(\LL)$, where $\LL$ is a grid and $\Psi$ is defined using 
data $d,m,n, \Vphys, \Vint, W$ satisfying {\bf (D), (I),
  (Reg)}. 
We can simultaneously rescale $\LL$, the window $W$, and the 
metric on $\Vphys$ by the same positive scalar, in order to assume that $\LL \in
\ALN_n$. Namely, set $c_1 \df
\covol(\LL)^{-\frac{1}{n}}$, so that $\LL_1 \df c_1 \LL \in \ALN_n$
 satisfies
$$
\Lambda = \Lambda(\LL, W) = \frac{1}{c_1} \Lambda(\LL_1, c_1W).
$$
Now 
solving for $c = \frac{1}{c_1}$ in \equ{eq: density exists} gives \equ{eq: for scaling}.

Define a 
sequence of measures $\eta_T$ on  
$\ALN_n$ by
$$
\eta_T \df \frac{1}{T} \int_0^T \int_{\Omega} (g_t u)_* \delta_{\LL}
\, dm_U(u) \, dt. 
$$
That is, the measures $\eta_T$ are defined by the same averaging as in
\equ{eq: nice convergence}, but for the action on $\ALN_n$ rather than
on $\Cl(\R^d)$. By \equ{eq: equivariance Psi}, their pushforward
under $\Psi$ are the measures appearing on the LHS of \equ{eq: nice 
  convergence}. 
By Theorem \ref{thm: Shah} we have $\eta_T \to_{T \to \infty} \bar \mu$ for some
homogeneous measure $\bar \mu$ on $\ALN_n$. By assertion (i) of
Theorem \ref{thm: MS 
  classification}, $\bar \mu$ is invariant under translation by any 
element of $\R^n$, and in particular any element of $\Vint$. Hence, by
Corollary \ref{cor: continuity}, $\bar \mu$ is a 
continuity point of the map $\Psi_*$. By \equ{eq: nice
  convergence}, $\Psi_* \eta_T \to \mu$ and by continuity, $\mu =
\Psi_* \bar \mu$.
\end{proof}

For the case in which $\mu$ is $\SL_d(\R)$-invariant but
not $\ASL_d(\R)$-invariant, we will need the following result:
\begin{lem}\name{lem: dense orbits H'}
  With the notation of Theorem \ref{thm: MS classification}, let
  $$
H'_1 \df g_1 H' g_1^{-1}
$$
(so that $H'_1 = H$ in the linear case and $H'_1$ is a Levi subgroup
of $H$ in the affine case), and let $v$ be
a nonzero vector in $\underline{\LL}_1$. Then the orbit of $v$ under the linear
action of $H'_1$ is an open dense subset of $\R^n$. 
  \end{lem}
  \begin{proof}
Write $v = g_1u$ for $u \in \Z^n \sm \{0\}$. It suffices to show that
the orbit $H'u$ is open and dense in $\R^n$. The linear action of $H'$
on $\R^n$ factors through the group $\underline{H}_1$ so we may
replace $H'$ with $\underline{H}_1$. 

The action of $\SL_{k}(\R)$ on $\R^k$ has the property that the orbit
of every nonzero vector is dense. The same is true for the action of
$\Symp_{2k}(\R)$ on $\R^{2k}$ (since any vector can be completed to a
symplectic basis), for the action of $\SL_k(\CC)$ on
$\CC^k \simeq \R^{2k}$, and for the action of $\Symp_{2k}(\CC)$ on
$\CC^{2k} \simeq \R^{4k}$. By Step 3 of the proof of Lemma \ref{lem:classification},
$\underline{H}_1$ is the product of groups $G_i$, and we have a direct product $\R^n =
\oplus_{i=1}^{r+s} V_i$, with the following properties:
\begin{itemize}
  \item
  For $i=1, \ldots, r$ we have a real field embedding $\sigma_i$, and
  $V_i = \sigma(\mathbf{V})_{\R}$; for $i=r+1, \ldots, r+s$ we have
  representatives $\sigma_i$ of pairs of complex embeddings, and $V_i
  = \sigma(\mathbf{V})_{\CC} $.
\item
  For $i=1, \ldots, r$ we have $G_i = \sigma_i(\mathbf{G})_{\R}$ and
  for $i=r+1, \ldots, s$ we have $G_i = \sigma_i(\mathbf{G})_{\CC}$. 
\item
  In the $\SL_k$-case (resp., the $\Sp_{2k}$ case), $V_1$ is
  isomorphic to $\R^k$ (resp., $\R^{2k}$), with the standard action.

\item
  The action of $G_i $ on $V_i$ is the
  obtained from the action of $G_1$ on $V_1$ by applying
  $\sigma_i$. In particular, for real embeddings it is isomorphic to
  the standard action of $\SL_k(\R)$ or $\Sp_{2k}(\R)$, and for
  complex embeddings it is isomorphic to the standard action of
  $\SL_k(\CC)$ or $\Sp_{2k}(\CC)$. 

  \end{itemize}
Thus, it is enough to show that for any $u \in \Z^n \sm \{0\}$, and for
any field embedding 
$\sigma_j$ of $\KK$, the projection $u_j $ of $u$ to the factor corresponding
to $\sigma_j$ is nonzero.

Suppose to the contrary that $u_j=0$ for some $j$, and let $a \in \SL_n(\R)$ be
a diagonalizable matrix, such that $a$ acts on the $\ell$-th factor of
$\R^n$ corresponding to the field embedding $\sigma_\ell$ as a scalar
matrix $\lambda_\ell\cdot \mathrm{Id}$, where the 
$\lambda_\ell$ are positive real scalars satisfying
$$
\lambda_j>1, \ \ \lambda_i<1 \text{ for } i \neq j, \ \ \text{ and} \
\ \prod_\ell
\lambda_\ell =1.
$$
That is, $a$ belongs to the centralizer of $H'$ in $\SL_n(\R)$, and $a^i
u \to_{i \to \infty} 0$. This implies by Mahler's compactness criterion that the
sequence $a^i \Z^n$ is 
divergent (eventually escapes every compact subset of $\LLN_n$). In
particular, the orbit of the identity coset $\SL_n(\Z)$ under the centralizer of
$H'$ is not compact. From this, via the implication $3 \implies 2$ in
\cite[Lemma 5.1]{EMS_gafa}, we see that 
$H'$ is contained in a proper $\Q$-parabolic subgroup of $\SL_n(\R)$,
and hence (see e.g. \cite[\S 11.14]{Borel_arithmetiques}) leaves
invariant a proper $\Q$-subspace of $\R^n$. This is a contradiction to {\bf (irred)}. 
    \end{proof}

  \begin{proof}[Proof of Theorem \ref{thm: MS surjective}, linear
    case] 
    We repeat the argument given for the affine case. The only
    complication is in establishing $\eta_T \to \bar \mu$ implies
    $\Psi_* \eta_T \to \Psi_* \bar \mu$, as  in the last
    paragraph of the proof. In the proof for the affine case, this was
    obtained from Corollary \ref{cor: continuity}, which shows that
    $\bar \mu$ is a 
    continuity point for the map $\Psi_*$, using the fact that
        $\bar \mu$ is invariant under translations by elements of
        $\Vint.$ In the linear situation $\bar \mu$ no longer has this
        continuity property. 

        To overcome this difficulty we argue as follows. We note that
        if 
    \eq{eq: can still show}{\bar \mu \left( \left\{ \LL \in \ALN_n:
          \piint(\LL) \cap \partial W \neq
          \varnothing\right\}\right)=0
    }
    then Corollary \ref{cor: continuity} can still be applied to show
    that $\bar \mu$ is a continuity point for $\Psi_*$.
Thus, we can assume from now on that \equ{eq: can still show} fails. 
Since $\supp \, \bar \mu
    = H\LL_1, $ this implies that the Haar measure $m_H$  of $H$
    satisfies

    \eq{NotSiegelMS}{
    m_H \left( \left
\{h \in H: \piint(  h \LL_1) \cap \partial W \neq \varnothing \right \} \right) >0.
}
Since $\LL_1$ is countable, 
there must be some $v \in \LL_1$ such that
\eq{eq: positive measure in boundary}{
m_H\left(\left\{h \in H: \piint( hv) \in \partial W \right\} \right)>0.
}
By Lemma \ref{lem: dense orbits H'}, there is a unique element $v_1\in
\R^n$ which is fixed by $H$ (namely $v_1 = g_1(0)$), and for any $v
\neq v_1$,  the orbit of $v$ under the action of $H$ is an open dense
subset of $\R^n$. In particular, if $v \neq v_1$ then the map
$h \mapsto hv$ sends $m_H$ to an absolutely continuous measure on
$\R^n$, and for such $v$ \equ{eq: positive measure in boundary} cannot
hold by {\bf   (Reg)}.

Thus, we must have $v =v_1$. In this case $hv=v$ and
$\piint(hv) \in \partial W$ for all $h \in H$. By 
examining the proof of Proposition \ref{prop: continuity
  Chabauty-Fell}, we see that the map
$$H \to  \Cl(\R^d), \ \ \ \ \ h \mapsto \Psi(h \LL_1)$$
is still continuous at any point outside a set of zero measure;
namely, the set of $h$ for which there is $v \neq v_1$ such that  $\piint(hv) \in\partial W $. Furthermore, the measure $\bar
\mu$ and the measures $\eta_T$ are all supported on the orbit
$H\LL_1$. Thus, we can apply the argument 
proving Corollary \ref{cor: continuity}, to see that the restriction
of $\Psi_*$ to measures supported on the orbit $H\LL_1$ is
continuous. This is sufficient to conclude that 
$\Psi_* \eta_T \to \Psi_* \bar \mu$ as $T \to \infty$. 
\end{proof}

\begin{remark} 
Theorem~\ref{thm: MS surjective} remains valid when one considers
other topologies (and
potentially, Borel structures) on
$\Cl(\R^d)$, as is done for example in \cite{Veech_siegel_measures, MSnew}. Thus, in
the terminology of \cite{Veech_siegel_measures}, the theorem is valid if
$\bar \mu$ is a Siegel measure giving full measure to cut-and-project sets.
Indeed, the only properties of the topology on $\Cl(\R^d)$ used in the proof
are the validity of Corollary~\ref{cor: continuity} (in the affine case)
and Proposition~\ref{prop: continuity Chabauty-Fell}, and the arguments deriving
Corollary~\ref{cor: continuity} (in the linear case).
These topological ingredients are easily seen to hold for the vague topology used in \cite{Veech_siegel_measures} and \cite{MSnew}. For example, for the analogue of Proposition~\ref{prop: continuity Chabauty-Fell}, see \cite[Lemma 5.14]{MSnew}. 
  \end{remark}

  \section{Some consequences of the classification}\label{sec:consequences}
   With Theorem \ref{thm: MS classification} in hand it is easy to obtain
  explicit descriptions 
  of RMS measures in low dimensions. Recall that we refer to the
  unique $\ASL_n(\R)$-invariant probability measure on $\ALN_n$ and the unique
  $\SL_n(\R)$-invariant probability measure on $\LLN_n$ as the {\em Haar-Siegel
  measures}. 

\begin{cor}\name{cor: MS Prop 2.1} With the notation above, suppose
  that 
  $\dim \Vphys > \dim \Vint$. Then the only affine RMS
  measure is the one for which $\bar \mu$ is the Haar-Siegel 
  measure on $\ALN_n$, and the only linear RMS
  measure is the one for which $\bar \mu$ is the Haar-Siegel 
  measure on $\LLN_n$. 
\end{cor}

This reproves a result stated without proof in 
\cite[Prop. 2.1]{MS}.

\begin{proof} 
In our classification result, there is $k \in \{d,
\ldots, n\}$ and $D = \deg(\KK/\Q) $ such that $n = D k$ in the
$\SL_k$-case and $n=2Dk $ in the $\Symp_{2k}$-case. Since
\eq{eq: strict inequality}{k \geq d = \dim \Vphys >\dim
  \Vint = n-d \geq n-k,}
we obtain $k > (D-1)k$ in the $\SL_k$-case and $k>(2D-1)k$ in the
$\Symp_{2k}$-case. This is only possible if $D=1$ and we are in the
$\SL_k$-case. 
That is, the only possible case is $H' = \SL_n(\R)$, and this gives
the  required result. 
  \end{proof}
We extend Corollary \ref{cor: MS Prop 2.1} to the case of equality:

\begin{cor}\name{cor: MS equality}
With the above notation, suppose that $\mu$ is not one of the
Haar-Siegel measures mentioned in Corollary 
\ref{cor: MS Prop 2.1}, and suppose $\dim \Vphys = \dim \Vint$. Then
either $d =2$ and $H' = \Symp_{4}(\R)$, or $d \geq 2$ and there is a
real quadratic field $\KK$ 
such that
$H'$ is (the group of real
points of) $\Res_{\KK/\Q}(\SL_d)$. 
\end{cor}

\begin{proof}
If the strict inequality in \equ{eq: strict inequality} becomes
non-strict, it is also possible that $H' =
\Res_{\KK/\QQ}(\SL_d)$ and $\KK$ is a real quadratic field, or $\KK = \Q,
\, d=2$ and $H' = \Symp_4(\R)$.  
  \end{proof}
As shown by Pleasants \cite{Pleasants}, an example of a cut-and-project
set  associated with a  real 
quadratic field as in Corollary \ref{cor: MS equality} is the vertex
set of an Ammann-Beenker tiling, where in this case the associated field is $\KK=\mathbb{Q}(\sqrt{2})$. Similarly, as discussed in \cite[\S
2.2]{MS}, the Penrose tiling vertex set can be described as a finite
union of cut-and-project 
sets associated with the real 
quadratic field $\mathbb{Q}(\sqrt{5})$.

\medskip

We record the following trivial but useful fact.

\begin{prop}\name{prop: origin in window}
For any affine RMS measure $\mu$, one can assume the window $W$
contains the origin in its interior. 
\end{prop}

\begin{proof}
Let $W$ be the window in the construction of the RMS measure $\mu$. By
{\bf (Reg)}, let $x_0 \in \Vint$ be a point in the interior of $W$. By assertion
(i) of Theorem \ref{thm: MS 
  classification}, the measure $\bar \mu$ is invariant under
translations by the full group $\R^n$ of translations, and in
particular by the translation by $x_0$. So we can replace any $\LL \in
\ALN_n$ by $ \LL-x_0$ without affecting the measure
$\bar \mu$. But clearly for $x_0 \in \Vint$ we have
$$\Lambda (\LL, W) = \Lambda(\LL-x_0, W-x_0).$$
So the measure $\mu$ can be obtained from $\bar \mu$ by using the
window $W-x_0$, which contains the origin in its interior. 
  \end{proof}
\medskip

Recall that we have an  inclusion
$$\iota: \SL_n(\R) \to \ASL_n(\R), \ \iota( g) = (g, \mathbf{0}_n),$$
i.e., $\iota(\SL_n(\R))$ is the stabilizer of the origin in the affine
action of $\ASL_n(\R)$ on $\R^n$.  
This induces an inclusion $\bar \iota : \LLN_n \to \ALN_n$, and these maps form
right inverses to the maps appearing in \equ{eq: def projection}:
$$
\pi \circ \iota = \mathrm{Id}_{\SL_n(\R)}, \ \ \underline \pi \circ
\bar \iota =
\mathrm{Id}_{\LLN_n}. 
$$
In the linear case, we can use these maps to understand the measures $\bar
\mu$ on $\ALN_n$ appearing in Theorem \ref{thm: MS classification} in
terms of measures on $\LLN_n$. Namely we have:

\begin{prop}\name{prop: restrictions linear}
  Let $F = \SL_d(\R)$, embedded in $\ASL_n(\R)$ via \equ{eq: embedding
    linear}, and let $\bar \mu$ be a measure on
  $\ALN_n$ projecting to a linear RMS measure on $\Cl(\R^d)$;
  i.e., $\bar \mu$ is $F$-invariant and ergodic, and not invariant
  under $\ASL_d(\R)$. Let $H, \, \LL_1$ be as in Theorem \ref{thm: MS
    classification}. Let $\underline F \df  \pi(F)$.  
  Then one of the following holds:
\begin{itemize}
\item[(i)] We have $\supp \,
  \bar \mu \subset \bar \iota(\LLN_n)$ and
  $\underline{\pi}|_{\supp \, \bar \mu}$ is a homeomorphism which maps $\bar \mu$ to
  an $\underline F$-invariant ergodic measure 
  on $\LLN_n$. In this case $H$ is contained in $\underline G \df \iota(\SL_n(\R))
  $, i.e., $H = \iota \circ \pi (H)$.  
\item[(ii)] We have  $\bar \mu (\bar \iota  (\LLN_n))=0$, and there are
  $D_1, D_2 \in \N$ such that $\underline{\pi}|_{\supp \bar \mu}$ is a
  closed map of degree $D_1$, and for every $\LL \in \supp \, \bar
  \mu$ there is a lattice $\LL' \in \LLN_n$,
depending only on $\underline \pi (\LL)$, such that $\LL'$ 
 contains $\underline{\pi}(\LL)$ with index 
$\left[\LL': \underline{\pi}(\LL)\right] = D_2,$ and such that $\LL$ is a
translate of $\underline \pi(\LL)$ by an element of $\LL'$. 
  
  \end{itemize}
  \end{prop}
  \begin{proof}
The set of lattices $\bar \iota(\LLN_n) \subset \ALN_n$ is
clearly $F$-invariant, so by ergodicity is either null or conull for
the measure $\bar \mu$. If it is conull then $\bar \iota(\LLN_n)$ is a
closed subset of full measure, i.e., $\supp \, \bar \mu \subset
\bar \iota(\LLN_n).$ Since $\bar \iota$ is a right inverse for
$\underline \pi$ we have that $\underline{\pi}|_{\supp \, \bar \mu}$
is a homeomorphism. Furthermore, since we have a containment of orbits
$$H  \LL_1 = \supp \, \bar
\mu \subset \bar \iota (\LLN_n) = \underline G  \Z^n = \underline
G  \LL_1,$$
and the groups $H, \underline G$ are connected analytic submanifolds of $G$, 
we have a containment of groups $H \subset \underline G$. This proves
(i). 

Now suppose  $\bar \mu \left(\bar \iota(\LLN_n)  \right)=0,$
and let $H, \, \LL_1$ be as in the statement of Theorem \ref{thm: MS
  classification}, so that $\supp \, \bar \mu = H 
\LL_1$.
Let $\mathbb{T}^n \df \underline \pi^{-1} (\underline \pi
(\LL_1)) 
$ be the orbit of $\LL_1$ under
translations. Since we are in the linear case, $H$ is transverse to
the group of translations $\R^n$ which moves along the fibers of
$\underline \pi$, and since $H \LL_1$ does not
accumulate on itself and 
$\mathbb{T}^n$ is 
compact, the intersection $\Omega \df \mathbb{T}^n \cap H \LL_1$ is a
finite set. Then by \equ{eq: equivariance property}, for any $\LL = h \LL_1
\in \supp \bar \mu$ we have
$$
h \Omega = \underline \pi^{-1}(\underline \pi(\LL)) \cap H\LL_1,
$$
and thus the map $\underline \pi|_{\supp \, \bar \mu}$ has fibers of a
constant cardinality $D_1 \df |\Omega|$.

Now denote 
$$
\Gamma_1 \df \{h \in H: h \LL_1 = \LL_1\}, \ \ \ \ \Gamma_2 \df \{h \in H:
h \Omega = \Omega\}.
$$
By equivariance we have $\Gamma_1 \subset \Gamma_2$ and the index of
the inclusion is $D_1$ since $\Gamma_2$ acts transitively on
$\Omega$. The bijection
$$ \R^n/\underline \pi(\LL_1 )\to \mathbb{T}^n, \ x \ \mathrm{mod} \,
\underline \pi(\LL_1) \mapsto x +\LL_1$$ 
endows $\mathbb{T}^n$ with the structure of a real torus, whose identity
element corresponds to $\LL_1$. In these coordinates $\Gamma_2$ acts
by affine maps of $\mathbb{T}^n$ but $\Gamma_1$ acts
by toral automorphisms, since it preserves $\LL_1$. Thus, $\Omega$ is a
finite invariant set for the action of an irreducible lattice in a
group acting $\LL_1$-irreducibly on $\R^n$, and thus
by \cite{guivarch_starkov} consists of torsion points in
$\mathbb{T}^n$. That is, there is $q 
\in \N$ so that they belong to the image of $\frac{1}{q} \cdot \LL_1 $ 
in $\mathbb{T}^n$. By equivariance the same statement holds, with the
same $q$, for $h \LL_1$ in place of $\LL_1$. Thus, the second
assertion holds if we let $\LL' = 
\frac{1}{q} \cdot \LL, \ D_2 = q^n$.  
    \end{proof}

    \begin{example}
It is possible that in case (ii) we have $\supp \bar \mu \cap
\bar \iota (\LLN_n) \neq \varnothing.$ For example, take $n=3, d=2, $
let $f$ be the translation $f(x) \df x + \frac{1}{2} \mathbf{e}_3, $
where $\mathbf{e}_3$ is the unit vector in the third axis. Let $H$ be
the conjugate of $\SL_3(\R)$ by $f$ and let $\LL_1 = f(\Z^3).$ Then $F
\subset H$ and  $H\LL_1$ is a closed homogeneous orbit.  Since
$\LL_1 \notin \bar \iota(\LLN_3)$, the corresponding homogeneous
measure does not satisfy  (i). But one can check
that the lattice $\spa_{\Z}(\mathbf{e}_1, 2 \mathbf{e}_2,
\frac12 \mathbf{e}_3)$ is contained in $
H\LL_1$, that is, $H\LL_1 \cap \bar \iota (\LLN_3) \neq \varnothing.$ 
      \end{example}
     \ignore{
Let $v \in \R^3$ be some vector not in the direction of $\mathbf{e}_3$
and let $h \in \SL_3(\R)$ such that $h (\mathbf{e}_3)  = v$ and $h(1,
0, -\frac{1}{2}) = \frac{1}{2}v-\frac{1}{2} \mathbf{e}_3$. Such an $h$ exists because we have
imposed two conditions on two independent vectors in $\R^3$. Let $h^f
\in H$ be the conjugate of 
$h$ by the translation $f$. Then
$h^f (1,0,-\frac{1}{2}) = \frac{1}{2} \mathbf{e}_3 + h ( f^{-1}(1, 0, \frac{1}{2})) =
\frac{1}{2} \mathbf{e}_3 + h(1,0,-\frac{1}{2}) - \frac{1}{2} h (\mathbf{e}_3) =\frac{1}{2}
\mathbf{e}_3 + \frac{1}{2} v -\frac{1}{2} \mathbf{e}_3 -\frac{1}{2} h (\mathbf{e}_3)=0,
$
so $h^f(\LL_1)$ is a grid containing $0$, i.e., a lattice. 
      }

\section{Integrability of the Siegel-Veech transform}\name{sec:
  reduction}
In this section we prove Theorem \ref{thm: reduction theory}.  Let $\mu$ be an RMS
measure and let $\bar \mu$, $H_1$, $\LL_1 = 
 g_1\Z^n$ be as in Theorem \ref{thm: MS classification}.  Recall
that the function $\hat f$ defined in \equ{eq: Siegel Veech transform}
is defined on $\supp \, \mu$. Also let
 $\pi: \ASL_n(\R) \to \SL_n(\R), \, \underline \pi: \ALN_n \to \LLN_n,
 \,\underline H_1 = \pi(H_1)$ be as in \S \ref{subsec: homogeneous
 measures}.  Let
$\Gamma_{1} \df H_1 \cap \ASL_n(\Z),\,
 \underline{\Gamma_1} \df \underline H_1 \cap \SL_n(\Z)$ be the $\Z$-points
 of $H_1$ and $\underline{H}_1$, and let $\mathbf{X}_1 \df
  H_1/\Gamma_1, \, \underline{\mathbf{X}}_1 \df \underline
 H_1/\underline \Gamma_1$.  We will use the results of
 \S \ref{subsec: homogeneous measures} to lift
$\hat f $ to a function on $\mathbf{X}_1$, and show that it is dominated by the
pullback of a 
function on $\underline{\mathbf{X}}_1$. For the arithmetic
homogeneous space
$\underline{\mathbf{X}}_1$ we will develop the 
 analogue of the Siegel summation formula and its
 properties. Specifically, we will describe a 
Siegel set $\sset \subset \underline H_1$, which is an easily described subset 
projecting onto $\underline{\mathbf{X}}_1$, and estimate the rate of decay of the
Haar measure of the subset of $\sset$ covering the 
`thin part' of $\underline{\mathbf{X}}_1$.

\subsection{Reduction theory for some arithmetic homogeneous spaces}
\label{subsec: reduction theory}
We begin our discussion of Siegel sets. For more 
details on the terminology and statements given below, see
\cite[Chaps. 11-13]{Borel_arithmetiques}. 

Let $\mathbf{H}$ be a semisimple $\Q$-algebraic group,
let $\mathbf{P}$ be a minimal $\Q$-parabolic subgroup,
and let $H = \mathbf{H}_\R$. Then $P = \mathbf{P}_{\R}$ has a decomposition $P 
= MAN$ (almost direct product), where:
\begin{itemize}
\item
 $A$ is the group of $\R$-points of a maximal
$\Q$-split torus $\mathbf{A}$ of $\mathbf{P}$;
\item
$N$ is the unipotent radical of $P$;
\item
 and $M$ is the connected component of the identity in the group of $\R$-points
of $\mathbf{M}$, a maximal
$\Q$-anisotropic $\Q$-subgroup of the centralizer of $\mathbf{A}$ in
$\mathbf{P}$. 
\end{itemize} 
Furthermore, $H = KP$ for a maximal compact subgroup $K$ of $H$.

As in \S \ref{subsec: number fields}, we think of $\mathbf{H}$ as concretely embedded in
$\SL_{n_0}(\R)$ for some $n_0 \in \N$, where we take this embedding to
be defined over $\Q$ for the standard $\Q$-structure on
$\SL_{n_0}(\R)$. 
Let $\mathfrak{a}$ and $\mathfrak{n}$ denote respectively the Lie
algebras of $A$ and $N$, let $\Phi \subset
\mathfrak{a}^*$ denote the $\Q$-roots of $H$ 
and choose an order on $\Phi$ for which $\mathfrak{n}$ is
generated by the positive root-spaces.

Every element of $H$
can be written in the form
\eq{eq: decompose h}{h = kman \ \ \ (k \in K, \, m\in M, \, a \in A, \, n \in N),}
and one can express the Haar
volume element $dh$ of $H$ in these coordinates in the form
\eq{eq: form of haar measure}{
dh = dk \, dm \, dn \, \rho_0(a) da,
}
where $dk, dm,dn, da$ denote respectively the volume elements
corresponding to the Haar measures on the (unimodular) groups $K, M, N,
A,$ and
\eq{eq: def rho}{
  \rho_0(a) =
  \left 
  |\det \left(\Ad(a)|_{\mathfrak{n}}\right)\right| 
=
\exp (2\rho(X)),
}
where $ a=\exp(X)$ and $\rho$ is the character on $\mathfrak{a}$ given by 
$\displaystyle{\rho= \frac{1}{2} \sum_{\alpha \in \Phi^+} c_\alpha \alpha,
}$
for $\Phi^+$ the positive roots in $\Phi$, and $c_\alpha = \dim \mathfrak{h}_{\alpha}$. 
We note that this formula for Haar measure is 
well-defined despite the fact that the decomposition \equ{eq:
  decompose h} is not unique.

Let $\Delta \subset \Phi^+$ be a basis of simple $\Q$-roots.
For fixed $t \in \R$, let
\eq{eq: def At}{
A_t \df \{\exp(X) : X \in \mathfrak{a}, \, \forall \chi \in \Delta, \,
\chi(X) \leq t\}
}
and for a compact neighborhood of
the identity $\omega \subset MN$, let 
$$
\sset_{t, \omega} \df K A_t \omega. 
$$
These sets are referred to as {\em Siegel sets}, and by a fundamental
result, a 
finite union of translates of Siegel sets contains a fundamental
domain for the action of an arithmetic group; 
that is, there is a
finite subset $F_0 \subset \mathbf{H}_{\Q}$ and there are $t, \omega$ such
that $\sset_{t, \omega} F_0 $ projects onto
$H/\Gamma_H,$ where $\Gamma_H=\mathbf{H}_{\Z}$; equivalently
$H = \sset_{t, \omega} F_0 \Gamma_H.$ The sets $\sset_{t,
  \omega} F_0$ do not represent $\Gamma_H$-cosets uniquely, in fact the
map $\sset_{t, \omega} F_0  \to H/\Gamma_H$ is far from being
injective. Nevertheless the formulas \equ{eq: decompose h} and
\equ{eq: def rho} make it possible to make explicit computations with
the restriction of Haar measure to $\sset_{t, \omega} F_0  $, and in
particular to show that Siegel sets have finite Haar measure.

An important observation is that the set $\bigcup_{a \in
A_t}a \omega a^{-1}$ is bounded, because of the definition of $M$ and
$N$ and because of the compactness of $\omega$. This
means that a Siegel set is contained in a set of the form $\omega' A_t$, where
$\omega'$ is a bounded subset of $H$.

    \subsection{The integrability exponent of an auxiliary function
      on $\LLN_n$}\name{subsec: alpha integrability}   We will
    specialize the discussion in \S \ref{subsec: reduction 
    theory} to the specific choices of $H/\Gamma_H$
    that arise in our application. 
    Let $H$ be as above, let $\sset_{t, \omega}$ be a Siegel set
    and let $F_0 \subset
    \mathbf{H}_{\Q}$ be a finite subset for which $\sset_{t, \omega}
    F_0 \Gamma_H=H$. Given functions $\varphi_1, \varphi_2$ defined on
    $H$, we will write 
   $\varphi_1 \ll \varphi_2$ if there is a constant $c$ 
    such that for all $x
    \in \sset_{t,\omega} F_0$ we have $\varphi_1(x) \leq c
    \varphi_2(x)$. The constant $c$ is called the {\em implicit
      constant}. We will also
    write $\varphi_1 \asymp \varphi_2$ 
    if $\varphi_1 \ll \varphi_2$ and $\varphi_2 \ll \varphi_1$. In
    general these relations on functions depend on the choice of Siegel set (i.e., the
    choice of $t$) and the choice of the finite set $F_0$, but in the
    case we will be interested in, when 
    $\varphi_1, \varphi_2$ are actually lifts of function defined on
    $H/\Gamma_H$, this notion does not depend on choices. 

  We now define an auxiliary function, and
    compute its integrability exponent. 
 Given a nonzero discrete subgroup $\LL' \subset \R^n$
 (not necessarily of rank $n$), we denote by $\covol(\LL')$ the
 volume of a fundamental domain for $\LL'$ in $\spa_{\R}(\LL')$
 (with respect to Lebesgue measure on $\spa_\R(\LL')$, normalized
 using the standard inner product on $\R^n$). For $g \in \SL_n(\R)$ and $\LL =g\Z^n \in
 \LLN_n$, 
 define 
 \eq{eq: def fn alpha}{
\hat \alpha(g) = \alpha(\LL) \df \max\left\{\covol(\LL')^{-1}:  \LL' \subset \LL, \,
  \LL' \neq \{0\}\right\}.
}
Recall that  $\underline{\mathbf{X}}_1 = \underline{H}_1/\underline{\Gamma}_1$
is embedded in $\LLN_n$ as the closed orbit $\underline{\mathbf{X}}_1 =
\underline{H}_1 \Z^n$, and so we can consider the restrictions of $\alpha$ and $\hat
\alpha$ to $\underline{\mathbf{X}}_1$ and to $\underline{H}_1$.

\begin{prop}\name{prop: alpha in Lp}
    In the two cases $\mathbf{G} \cong \SL_k, \ \mathbf{G} \cong \Symp_{2k}$,
    let    $p< r_0 \df  \rank_{\KK} (\mathbf{G})+1$ (see \equ{eq: K
      rank}). Then
    \eq{eq: integrability exponent}{
      \alpha \in L^p\left(\underline\mu \right) \sm
    L^{r_0}\left(\underline{\mu} \right),}
    where $\underline{\mu}$ is the
    $\underline{H}_1$-invariant probability measure on $\underline{\mathbf{X}}_1$.
\end{prop}

\begin{proof}
Let $\lambda_i = \lambda_i(\LL), \ i=1, \ldots, n$ be the successive
minima of a lattice $\LL$, and let $i_0 = i_0(\LL)$ be the index for which
$\lambda_{i_0} (\LL) \leq 1 < \lambda_{i_0+1}(\LL)$. Then it is easy to see
using Minkowski's second theorem (see e.g.\,\cite[\S VIII.2]{Cassels})
that (as functions on $\LLN_n$), 
\eq{eq: as functions}{
  \alpha(\LL) \asymp \left(\lambda_1 \cdots
  \lambda_{i_0(\LL)}(\LL)\right)^{-1}.}
As a consequence, for any $C \subset \SL_n(\R)$ bounded, we have
$$
\forall u \in C, \ \ \ \ \alpha(u\LL) \asymp \alpha(\LL)
$$
(with the implicit constant depending on $C$).

Let $\mathbf{T}$ denote the diagonal subgroup of $\SL_n(\R)$, let
$T  = \mathbf{T}_\R^\circ$ and let $\mathfrak{t}$ be the Lie algebra of
$T$. In what follows we will replace $T$ by its conjugate over
$\SL_n(\Q)$, where the conjugate will be conveniently chosen with
respect to $\underline{H}_1$ and its subgroups. The reader should note
that the statements to follow 
about $T$ are not affected by such conjugations in $\SL_n(\Q)$.

It is easy to check that for the lattice $\Z^n$ and for $a =
\exp(\diag(X_1, \ldots, X_n)) \in T$, we have
$\lambda_i(a\Z^n) = e^{X_{j(i)}}$ where $i \mapsto j(i)$ is a permutation
giving $X_{j(1)} \leq X_{j(2)} \leq \cdots \leq
X_{j(n)}$, and hence
\eq{eq: formula 1}{
  \hat{\alpha}(a)  = \alpha(a\Z^n) \asymp
\exp\left(-\sum_{X_i <0} X_i \right)
.}
Furthermore, for an element $f_0 \in \SL_n(\Q)$ we have that 
$\lambda_i(af_0\Z^n) \asymp e^{X_{j(i)}}$, where implicit constants 
depend on $f_0$, and thus $\hat\alpha (a) \asymp \hat \alpha(af_0)$.

Recall the notation $D = \deg(\KK/\Q)$ from Theorem
\ref{thm: MS classification}. We first prove the proposition
under the assumption $D=1$. That is, we have $\KK = \Q, \
\underline{H}_1 = \SL_{k}(\R)$ and $n=k$ in case $\mathbf{G} \cong \SL_k,$ and
$n=2k, \ \underline{H}_1 = 
\Symp_{2k}(\R)$ in case $\mathbf{G} \cong \Symp_{2k}$. 
Now consider a Siegel set for $H =
\underline H_1$,
  and suppose
  $A_t$ is the corresponding subset of the maximal 
  $\Q$-split torus of $\underline{H}_1$. Since $\mathbf{T}$ is a
  maximal $\Q$-split torus of 
  $\SL_n(\R)$, by \cite[Thm. 15.14]{Borel1}, applying a 
  conjugation in $\SL_n(\Q)$ we can assume that $\mathbf{A} \subset \mathbf{T}$ and the order on the
  roots $\Phi$ is consistent with the standard order on the group of
  characters on $\mathfrak{t}$; that is, $A_t \subset T_{t'}$ for some
  $t'$, as can be observed  by an elementary computation (see
  \cite[Ex. 11.15]{Borel_arithmetiques} for a description of
  $\mathbf{A}$ in the symplectic case). In particular, for $a
  =\exp(\diag(X_j)) \in A_t$ we have 
  $\exp(X_j) \ll \exp(X_{j+1})$ for $j=1, \ldots, n-1$.
Then from \eqref{eq: formula 1}, for $a \in A_t$ and
$f_0\in F_0$, where $F_0$ is a finite subset 
of $\left(\underline{H}_1\right)_\Q$,  we have
\eq{eq: alpha beta}{
\hat\alpha(af_0)  \asymp \max_{1 \leq j \leq n-1} \exp \left(-\beta_j(X) \right),  
}
where
\eq{eq: def betaj}{
  \beta_j(X) \df
\sum_{i=1}^{j}X_i, \ \ \ \  \ X = \diag(X_\ell).
}
Since a Siegel set $\sset_{t, \omega}$ is contained in a set of the form
$\omega' A_t$, where $\omega'$ is a compact subset of $H$, this
implies that
$$
\hat\alpha(kmanf_0) \ll\max_{1 \leq j \leq n-1}\exp(-\beta_j(X)).$$
We will first show the following:
\begin{itemize}
\item[(i)]
For any $j$, and any $X \in \mathfrak{a}$ for which $\exp(X) \in A_t$,
we have $\left(2 \rho - r_0 \beta_j \right) (X)  \ll 1.$
\item[(ii)] The number $r_0$ is the largest number for which the
  conclusion of (i) remains valid. 
\end{itemize}
For $\ell=1, \ldots, n-1$ let $\chi_\ell$ denote the simple roots on $\mathfrak{t}$,
that is, 
\eq{eq: def simple roots}{
\chi_\ell: \mathfrak{t} \to \R,  \ \ \ \chi_\ell \left(\diag\left(X_1,
    \ldots, X_{n} \right) \right)\df X_{\ell+1}-X_\ell. 
}
In order to show (i), since the $\chi_\ell$ are bounded above
on
$A_t$, it suffices to show
that if we write $2\rho = \sum a_\ell \chi_\ell$ and $\beta_j =
\sum b^{(j)}_\ell \chi_\ell$, then $r_0 b^{(j)}_\ell \leq
a_\ell $. In order to show (ii) it suffices to check that there are 
some  $j, \ell $ for which equality holds, i.e., $r_0 b^{(j)}_\ell =
a_\ell $. This can be checked using the tables of  
\cite[pp. 265-270, Plates I \& III]{Bourbaki_lie_groups} (note that
the restrictions of the $\beta_j$ to $A$ are
the fundamental weights
%
in both cases). Namely, for $\mathbf{G} = \SL_k$ we have
\[
  a_\ell = 
  \ell (k-\ell), \ \ \ \ 
r_0 \, b^{(j)}_\ell =  
\left\{ \begin{matrix} \ell(k-j) & \text{ if } \ell < j \\j (k
     -\ell)   & \text{ if } \ell \geq j \end{matrix} \right.
  \]
  and we have the desired inequality, with equality when $\ell =
  j$.
  If $\mathbf{G} = \Symp_{2k}$ we have 
  \[
    a_\ell = \left\{\begin{matrix} \ell (2k - \ell +1) & \text{ if } \ell < k \\\frac{k(k+1)}{2} & \text{
          if } \ell = k \end{matrix}  \right. 
    \ \ \ \
    r_0 \, b^{(j)}_\ell = \left\{ \begin{matrix} \ell (k+1)& \text{ if }
        \ell < j\\  j (k+1)& \text{ if } j
        \leq \ell < k \\ \frac{j(k+1)}{2} & \text{
          if } \ell = k \end{matrix} \right. 
    \]
and again the inequality holds, with equality when $\ell = j = k$.

Now to see that $\alpha \in L^p\left(\underline{\mu}\right)$, since a
Siegel set is contained in $\omega' A_t$ with $\omega'$ bounded, and
by \equ{eq: form of haar measure}, it 
suffices to prove that for $f_0 \in F_0$
we have $\int_{A_t} \hat {\alpha}^p(af_0) \rho_0(a) da <
\infty$. Using the preceding discussion, if we let $\mathfrak{a}_t$
denote the cone in $\mathfrak{a}$ with $A_t = \exp(\mathfrak{a}_t)$
(where $A_t$ is as in \equ{eq: def At}),
and use that $da$ is the pushforward under the exponential map of $dX$, we have
\[
\begin{split}
& \int_{A_t} \hat {\alpha}^p(af_0) \rho_0(a) \, da \ll
\int_{\mathfrak{a}_t} \max_j \exp \left( -p \beta_j(X) \right) \cdot
\exp \left(2\rho(X) \right) \, dX \\  = & \int_{\mathfrak{a}_t} \max_j
\exp \left[\left( \frac{p}{r_0} \left(2\rho - r_0 \beta_j\right)+ \left(1-\frac{p
    }{r_0} \right) 2\rho \right) (X)\right] 
\, dX \\
 \stackrel{(i)}{\ll} & \int_{\mathfrak{a}_t} \exp \left [ 2\rho(X)\right]^{1-\frac{p}{r_0}}
 \, dX <\infty,
  \end{split}
  \]
where the integral is finite as the integrand is the 
exponential of a linear functional
which is strictly decreasing along the
cone $\mathfrak{a}_t$. 
The same computation and (ii) show that we have a
corresponding lower bound $\int_{A_t} \hat {\alpha}^{r_0}(af_0) \rho_0(a)
da \gg \int_{\mathfrak{a}_t} \exp \left(\tau(X) \right) \, dX$, where
$\tau$ is a linear functional which is constant along a face of
$\mathfrak{a}_t$. We have shown \equ{eq: integrability exponent} for
$D=1$.

\medskip


Now suppose $D>1$. Our strategy will be to show that we can repeat the
computations used for the case $D=1$, with the only difference being
that in 
some of the formulas, the characters $\rho$ and $\beta_j$ are
multiplied  by a factor 
of $D$.
Write  $G_1 \df
{}^{\sigma_1}\mathbf{G}_\R$, let $V$ be as in the statement of Theorem
\ref{thm: MS classification}, a $\KK$-subspace of $\R^n$.
Let
\begin{equation}\label{eq: def t}
t \df  \left\{ \begin{matrix} k & \text{ if } \mathbf{G} \cong \SL_k \\
      2k &
    \ \, \text{ if } \mathbf{G} \cong \Symp_{2k}, \end{matrix}\right.
  \end{equation}
so that $\dim V = t$.  Let $\mathbf{A}_1$
denote a maximal $\KK$-split torus in $\mathbf{G}$, and let
$\mathfrak{a}_1$ denote its Lie algebra. Then, with
respect to a
suitable basis of $V_{\KK}$, we can write elements of $\mathfrak{a}_1$
as matrices $\diag(X_1, \ldots, X_t)$, where $\sum X_i =0$ when
$\mathbf{G} \cong \SL_k$ and $X_{i+k} = -X_i$ when $\mathbf{G} \cong
\Symp_{2k}$.

Let $\mathbf{B} \df
\Res_{\KK/\Q}(\mathbf{A}_1)$, 
and let $\mathbf{A}$ denote a maximal $\Q$-split torus in
$\underline{H}_1$. The dimension of $\mathbf{A}_1$ is the number of
independent one-parameter multiplicative $\KK$-subgroups (morphisms
$\KK^\times \to \mathbf{A}_1$), 
and, applying restriction of scalars, each such one-parameter group
gives rise to a one-parameter $\Q$-subgroup 
$\Q^\times \to \mathbf{B}$. This implies that $\mathbf{B}$ contains a
$\Q$-split torus of dimension equal to $\dim \mathbf{A}_1$.  
Since
the $\Q$-rank of $\mathbf{H}$ is the same as the $\KK$-rank of
$\mathbf{G}$, see \cite[6.21 (i)]{BT_reductif}, 
the dimensions of these groups coincide. Since all maximal $\Q$-split
tori in $\mathbf{H}$ are conjugate over $\mathbf{H}_{\Q}$, we can
assume that 
 $ \mathbf{A} \subset \mathbf{B},$
and by conjugating $\SL_n(\R)$ by an element of $\SL_n(\Q)$, we can
also assume that $\mathbf{A} \subset \mathbf{T}$ and the order on the
roots $\Phi$ is consistent with the order on the roots of
$\mathfrak{t}$. We claim that after these conjugations, the elements
of $A = \mathbf{A}^\circ_{\R}$ are of the form
\eq{eq: torus embedding2}{
\diag\left(\underbrace{X_{j(1)}, \ldots, X_{j(1)}}_{D \text{ times}}, \ldots,
  \underbrace{X_{j(t)}, \ldots, X_{j(t)} }_{D \text{ times} }\right), 
}
where $\diag(X_1, \ldots,
X_t)$ ranges over the elements of $\mathfrak{a}_1$ in the
above-chosen basis, and 
$i \mapsto j(i)$ is a permutation guaranteeing $\exp\left(X_{j(1)}
\right) \ll \cdots \ll \exp \left( X_{j(t)} \right)$. 

\ignore{
Let $\mathfrak{a},
\mathfrak{a}_1, \underline{\mathfrak{h}}_1$ denote the Lie algebras of
$A, A_1, \underline{H}_1$, and let $\Delta : \mathfrak{a}_1 \to
\mathfrak{sl}_{n}(\R)$ be defined by
\eq{eq: torus embedding}{
\begin{split}
  &
  \Delta(\diag(X_1 , \ldots, X_t)   ) \df  \diag
\\  &
  \left(\underbrace{X_1,
      \ldots, X_t, \ldots, X_1, \ldots, X_t}_{r \text{ times} }
,
 \underbrace{X_1, X_1,  \ldots, X_t, X_t, \ldots, X_1, X_1, \ldots,
    X_t, X_t}_{s \text{ times} }
\right),
 \end{split}
}
where $t = k $ in case $\mathbf{G} \cong \SL_k$ and $t = 2k$ in case
$\mathbf{G} \cong \Symp_{2k}$. 
We claim that
\begin{equation}\label{eq: claim z}\begin{split}
  &\text{there  is } g_0 \in \SL_n(\R), \text{ a bounded 
  }
  C\subset \SL_n(\R) \text{ and a permutation matrix } \\ & P  
 \text{ such that } 
Pg_0Ag_0^{-1}P^{-1}\subset 
C \cdot \exp(\Delta(\mathfrak{a}_1)).
    \end{split}
  \end{equation}
  }
We first assume the validity of \eqref{eq: torus embedding2}, and conclude the
proof of the case $D>1$.
\ignore{
  Note that the image of \equ{eq: torus embedding} might not satisfy the 
condition of being in a good position; in order to achieve this one
needs to conjugate with a permutation matrix, which transforms 
\equ{eq: torus embedding} into the map
\eq{eq: torus embedding2}{\diag(X_1, \ldots, X_t) \mapsto
\diag\left(\underbrace{X_{j(1)}, \ldots, X_{j(1)}}_{D \text{ times}}, \ldots,
  \underbrace{X_{j(t)}, \ldots, X_{j(t)} }_{D \text{ times} }\right), 
}
with $i \mapsto j(i)$ a permutation guaranteeing $\exp\left(X_{j(1)}
\right) \ll \cdots \ll \exp \left( X_{j(t)} \right)$. So the claim
means that $g_0Ag_0^{-1}$ is within a bounded distance of the image of \eqref{eq:
  torus embedding2}. 
}
We will use \equ{eq: torus embedding2} to compare
characters on $A_1$ with characters on $A$. First, comparing 
 the character
$\rho$ appearing in \equ{eq: def rho} for the two groups
$\underline{H}_1\, , G_1,$ we see that each
real field embedding $\sigma_i, i\leq r$ contributes one 
dimension to the dimension of a root space, and each pair $\sigma_i,
\bar \sigma_i, \, i>r$ of
conjugate non-real embedding contributes two dimensions. Alternatively:
in $G_1$ the root spaces are one dimensional and defined over $\KK$,
since $G_1$ is $\KK$-split. The root spaces in $H_1$ are obtained from
the root spaces 
in $G_1$ by applying the restriction of scalars operation to each one
individually. This implies
that the character $\rho$ for $\underline{H}_1$ is obtained from the
corresponding character for $G_1$ by a multiplication by
$D$. Similarly, it is clear from \equ{eq: torus embedding2} that the
characters $\beta_j$ appearing in \equ{eq: def betaj} 
for $\underline{H}_1$ are obtained from the same characters $\beta_j$
for $G_1$, multiplied by $D$. Thus, the computations guaranteeing
\equ{eq: integrability exponent} for $D=1$, imply the same property
for general $D$.

It remains to prove \eqref{eq: torus embedding2}.
\ignore{
  Let $\KK^\times$ denote the
multiplicative group of the field $\KK$. Since $\mathfrak{a}$ is
$\KK$-split and $\KK^\times \subset \R$ is
dense, it is enough to show \eqref{eq: claim z} for $X_i \in \log
\KK^\times$. 
We first prove this under the
additional assumption that $D=r$, that is, the field $\KK$ is totally
real. In this case the Galois extension $\mathbb{L}$ is contained in
$\R$ and so there is  $g_0
\in \SL_n(\R)$ such that $g_0Ag_0^{-1} \subset T$.

Let $\mathbf{S}$ be a $\QQ$-anisotropic $\Q$-torus in $\mathbf{H}$
which commutes with $\mathbf{A}$ and such that
$\mathbf{S}\mathbf{A}$ is a maximal torus  
of $\mathbf{H}$. The existence of $\mathbf{S}$ follows from \cite[\S
8.15]{Borel1}. Then the orbit 
$S \Z^n \subset \LLN_n$ is bounded, where $S = \mathbf{S}_{\R}$, and
there is a bounded subset $C_1\subset S$ such 
that $S\Z^n = C_1\Z^n$. 

}
Recall that $\mathbf{B} = \Res_{\KK/\Q}(\mathbf{A}_1)$, which we wish
to describe explicitly using the discussion in \S
\ref{subsec: number fields}. For $\vec{y} \in \KK^t$ we define
$ a_1(\vec{y}) \df \diag(y_1, \ldots,
y_t) \in A_1(\KK);
$ that is, these are matrices acting on $V$ which are diagonal with
respect to a $\KK$-basis 
of $V$, and the $y_i$ satisfy $y_1+\cdots + y_t =0$ for $\mathbf{G}
\cong \SL_k$ and 
$y_i = -y_{2k-i+1}$ for $\mathbf{G} \cong \Symp_{2k}$. 
Each $y \in \KK$ has a representative which is a matrix in 
$\Mat_{D \times D}(\Q)$. If we take $y \in \Q$ then the corresponding
representative 
matrix is the scalar matrix $y \cdot \mathrm{Id}_D$. The elements of
$\mathbf{B}$ can be considered as $t \times t$ matrices, whose
entries are elements of $\Mat_{D \times D}$. In particular, for
$\vec{y} \in \Q^{t}$, we get matrices $a_2(\vec{y}) \in \Mat_{n
  \times n}(\Q)$, which are simultaneously diagonalizable,
with each $y_i$ appearing as an eigenvalue $D$ times. That is, up to
permuting the coordinates, the matrices $a_2(\vec{y})$ are as in
\eqref{eq: torus embedding2}, with $X_i \in \Q$. The map
$a_1(\vec{y}) \mapsto a_2(\vec{y})$ is a polynomially defined group
homomorphism. Letting $\mathbf{A}_2$ denote the Zariski closure of
$\{a_2(\vec{y}): \vec{y} \in \Q^t, \, a_1 (\vec{y}) \in A_1\}$, we see that
$\mathbf{A}_2$ is a torus in $\mathbf{B}$ whose group of real points
$(\mathbf{A}_2)_{\R}$ satisfies the description \eqref{eq: torus
  embedding2}, and with $\dim \mathbf{A}_2 = \dim \mathbf{A}_1 = \dim
\mathbf{A}$. Also, 
$\mathbf{A}_2$ is $\Q$-split since the maps 
$a_2(\vec{y}) \mapsto y_i$ are $\Q$-characters. 
  Thus, $\mathbf{A}_2$ is a maximal $\Q$-split
  torus of $\mathbf{H}$, and by the uniqueness of the maximal
  $\Q$-split torus in the torus $\mathbf{B}$ (see
  \cite[Prop. 10.6]{Borel_arithmetiques}), we 
  must have $\mathbf{A} = \mathbf{A}_2$. (See also the related
  discussion in \cite[Example, 
  p. 54]{Rapinchuk_Platonov}, giving an explicit description of
  a maximal $\Q$-anisotropic torus in $\mathbf{B}$ as a product of norm-tori.)
\end{proof}

\subsection{An upper bound for the Siegel transform}
We will now state and prove a result implying Theorem \ref{thm:
  reduction theory}. For a function $F$ on $\R^n$, a measure $\bar
\mu$ on $\ALN_n$, and $\LL \in
\ALN_n$, in analogy with \equ{eq: Siegel Veech transform} we denote
\eq{eq: if define 2}{
\widehat F (\LL)  = \left\{\begin{matrix} \sum_{x \in \LL \sm \{0\}} F(x)
  & \bar \mu \text{ is linear} \\ \sum_{x \in \LL} F(x) & \bar \mu
    \text{ is affine } \end{matrix} \right. 
}
\begin{thm}\name{thm: reduction theory more general}
Let $\bar \mu$ be the $H$-homogeneous measure on $
\ALN_n$ as in Theorem \ref{thm: MS classification}, 
and let $q = q_{\bar \mu}$ be as in \equ{eq: p bounded}. 
Then for any $F \in C_c(\R^n)$ and any $p<q$ we have $\widehat{F} \in
L^p(\bar \mu)$. Moreover, there are $F \in
  C_c(\R^n)$ for which $\widehat F \notin L^{q}(\bar \mu)$.
\end{thm}

We will prove Theorem \ref{thm: reduction theory more general} separately in the
linear and affine cases. In the linear case, we will first
show, using Proposition \ref{prop: restrictions linear}, that the
Siegel-Veech transform \equ{eq: 
  if define 2} can be bounded in terms of a 
Siegel transform of a function on $\LLN_n$. The latter can be bounded
in terms of the function $\alpha$ considered in \S \ref{subsec: alpha
  integrability}.
\begin{proof}[Proof of Theorem \ref{thm: reduction theory more general}, linear
    case]
 Suppose  that $\bar \mu$ satisfies (i) of Proposition \ref{prop:
  restrictions linear}, i.e., $\bar \mu$ is supported on $\overline
\iota (\LLN_n)$. Then we can
assume that the
cut-and-project scheme involves lattices in $\LLN_n$, rather than
grids. Moreover, $H = \iota\circ \pi(H), \, \underline g_1 = g_1, \,
\underline H_1 = \iota \circ \pi (H_1)$, and
the function $\widehat F$ is a Siegel-Veech transform of a Riemann 
integrable function on $\R^n$, for a  homogeneous subspace  of $\LLN_n$. 
It is known that the function $\alpha$ defined in \equ{eq: def fn
  alpha} describes the growth rate of the 
Siegel transforms of functions on $\LLN_n$.
Namely (see \cite[Lemma 3.1]{EMM} or \cite[Lemma 5.1]{KSW}), 
for any Riemann integrable
function $F$ on $\R^n$, for
any $\LL \in \LLN_n$, $\widehat{F}(\LL) \ll 
\alpha(\LL)$. Furthermore, if  $F$ is the indicator of a ball around
the origin then $\widehat{F}(\LL)\gg
\alpha(\LL)$. 
Thus, the conclusion of Theorem 
  \ref{thm: reduction theory more general} in this 
case follows from Proposition \ref{prop: alpha in Lp}.

Now assume that case (ii) of Proposition \ref{prop: restrictions
  linear} holds. 
We cannot use Proposition \ref{prop: alpha
  in Lp} since $\widehat F$ is a function on $\ALN_n$. To remedy this, we
define for each $\LL \in H  \LL_1$ the lattice
$\LL' = \LL'(\underline \pi(\LL))$  appearing in assertion (ii) of Proposition \ref{prop: restrictions
  linear}, and set
$$
\widehat{\underline{F}} (\underline \pi(\LL)) \df \sum_{x \in \LL'(\underline \pi (\LL)) \sm \{0\}} F(x).  
$$
Then the bounds given in Proposition \ref{prop: restrictions
  linear} imply that 
$
\widehat F(\LL) \ll \widehat{ \underline{F}}(\underline \pi (\LL)),
$
with a reverse inequality
$
\underline{\widehat{ F}}(\underline \pi (\LL)) \ll \widehat F(\LL)
$
 for positive $F$. Since $\underline{\widehat{ F}}$ is the Siegel-Veech transform
 of a function on $\R^n$ with respect to a measure on $\LLN_n$, we can
 apply Proposition \ref{prop: alpha in Lp} to conclude the proof in
 this case as well. 
  \end{proof}

 For the affine case, we will need the following additional interpretation of the function
 $\alpha$ defined in \equ{eq: def fn alpha}.

 \begin{prop}\name{prop: alpha again}
  Let $\underline \LL \in \LLN_n$, let $\mathbb{T}^n_{\underline{\LL}} = \mathbb{T}^n =
  \underline{\pi}^{-1}(\underline \LL)
  \cong \R^n/\underline \LL$ be the quotient torus, equipped with its
  invariant measure element  $d\LL$. Then for any ball
  $B \subset \R^n$ and any $p>1$ we have 
   \eq{eq: without affecting}{
\int_{\mathbb{T}^n} |B \cap \LL|^p \, d\LL \asymp \alpha(\underline\LL)^{p-1},
   }
where the implicit constants depend on the dimension $n$, on $p$, and
on the radius of $B$.
\end{prop}

\begin{proof} Let $\lambda_1, \ldots, \lambda_n$ be the Minkowski
  successive minima of $\underline \LL$. Using Korkine-Zolotarev reduction, 
  let $v_1, \ldots, v_n$ be a basis for $\underline \LL$ satisfying
  $\|v_i \|\asymp \lambda_i$ (where implicit constants are allowed to
  depend on the dimension $n$), and let $u_i \df
  \frac{v_i}{\|v_i\|}$. For a vector $\vec{s}$ of positive numbers
  $s_1, \ldots, 
  s_n$ define 
  $$P_{\vec{s}} \df 
  \left\{\sum a_i u_i : |a_i| \leq \frac{s_i}{2} \right\}.
  $$
 Setting $\vec{v}_0 = (\|v_1\|, \ldots, \|v_n\|)$, we have that
 $P_{\vec{v}_0} = \left\{\sum b_i v_i: |b_i|\leq \frac12 \right\}$ is a fundamental
 parallelepiped for $\LL$, and we can 
 identify $\mathbb{T}^n$ with this parallelpiped via the bijection
 $$
P_{\vec{v}_0} \to \mathbb{T}^n, \ \  x \mapsto \LL_x \df \underline \LL +x,
$$
which sends the Lebesgue measure on $P_{\vec{v}_0}$ to the Haar measure
$d \vol$ on $\mathbb{T}^n$. 

Now set
  $$P_r \df P_{\vec{r}} \ \text{ where } \vec{r} = (r, \ldots, r).
  $$
  We can translate $B$ so that it is centered at the origin without
  affecting the integral in \equ{eq: without affecting}, and 
  since there is a lower bound on the angles between the $v_i$, there
  are $r_1 \asymp R \asymp r_2$ such that $P_{r_1} \subset B \subset P_{r_2}
  $. Thus, we can replace $B$ with $P_R$. Furthermore, the lower bound
  on the angles between the $u_i$ implies 
  $$
  d \vol (x)  
  \asymp  dx_1 \cdots dx_n, \ \
\text{ where } x = \sum x_i u_i.
  $$

  Writing each
  vector $y \in \R^n$ in the form $y= \sum_i c_i u_i$, and reducing each
  $c_i$ modulo $\|v_i\| \cdot \Z$, it is easy to verify that for $x
  \in P_{\vec{v}_0}$ we have:
  \begin{itemize}
  \item
    if $ \frac{R}{2} < |x_i| < \frac{\|v_i\|-R}{2} $ for some
    $i$, then $P_R \cap  
    \LL_x = \varnothing$; and
  \item
    if $|x_i | \leq
  \frac{R}{2}$ or $|x_i| \geq \frac{\|v_i\|-R}{2} $ for all $i$, then
  $|P_R \cap \LL_x| \asymp \prod_{\|v_i\|< R}\left(
  \frac{R}{\|v_i\|} \right)
  $.
\end{itemize}
Since
$$\prod_{\|v_i\| <R} 
  \frac{R}{\|v_i\|}  \asymp \prod_{\lambda_i(\underline
    \LL) <1} \frac{1}{\lambda_i(\underline \LL)} \stackrel{\equ{eq: as functions}}{\asymp}
  \alpha(\underline \LL),$$
we obtain
  \[\begin{split}
       & \int_{\mathbb{T}^n} |B \cap \LL|^p \, d\LL \asymp  
  \int_{P_{\vec{v}_0}} |P_R \cap \LL_x|^p \, d\vol(x)  \\ \asymp \ & 
  \alpha(\underline \LL)^p \cdot \vol\left( \left\{x \in P_{\vec{v}_0} :
   |x_i| \leq   \frac{R}{2}\right \} \right) \\ \asymp \  & \alpha (\underline
\LL)^p \cdot \prod_{\|v_i\|\leq R} \|v_i\|\cdot \prod_{\|v_i\|>R} R 
\asymp \
\alpha(\underline \LL)^p \cdot \prod_{\lambda_i (\underline \LL) < 1}
\lambda_i(\underline \LL)
\stackrel{\equ{eq: as functions}}{\asymp} \alpha(\underline \LL)^{p-1} .
\end{split}\]
\end{proof}

\begin{proof}[Proof of Theorem \ref{thm: reduction theory more general}, affine case]
By decomposing $F$ into its positive and negative parts, we see that
it suffices to prove $\widehat F \in L^p(\mu)$
when $F$ is the indicator of a
ball in $\R^n$. By Theorem \ref{thm: MS classification} we have that
in the affine case, the translation group $\R^n$ is contained in $H_1$,
which implies that we can decompose the measure $\bar \mu$ as 
$$\int_{\mathbf{X}_1} \varphi (\LL) \, d\bar \mu (\LL) = \int_{\underline{\mathbf{X}}_1}
\int_{\mathbb{T}^n_{\underline{\LL}}} \varphi(\LL_x) \, d\vol(x) \, d\underline
\mu(\underline \LL) , \ \ \forall \varphi \in L^1(\mathbf{X}_1, \bar \mu).$$
Now the statement follows from Propositions \ref{prop:
alpha in Lp} and
\ref{prop: alpha again}. The case of equality $p = q_\mu$ follows
similarly, taking for $F$ the indicator of a ball in $\R^n$.
\end{proof}

\begin{proof}[Proof of Theorem \ref{thm: reduction theory}]
Let $f \in C_c( \R^d)$ and let $\hat
f$ be as in \equ{eq: Siegel Veech transform}. Let $\mu$ be an
RMS measure on $\Cl(\R^d)$ associated with a cut-and-project scheme
involving grids in $\ALN_n$, a decomposition $\R^n = \Vphys \oplus
\Vint$, and a window $W \subset \Vint$. Let 
$\mathbf{1}_W$ be the indicator function of $W$ and let $\bar
\mu$ be an $H$-homogeneous 
measure, supported  on the orbit $H  \LL_1 \subset \ALN_n$ such
that $\mu = \Psi_* \bar \mu$ (where we have replaced $\mu$ by its
image under a rescaling map to simplify notation). 

Define
\eq{eq: if define 1}{
F: \R^n \to \R, \ \ F(x) = \mathbf{1}_W(\piint(x)) \cdot
f(\piphys(x)),}
and define $\widehat F$ via \equ{eq: if define 2}.
Then it is clear
from the definition of $\Psi$ and \equ{eq: Siegel Veech transform}
that
$
\hat f (\Psi(\LL)) = \widehat F (\LL)
$ provided $\LL$ satisfies {\bf (I)}, and, in the linear case,
provided all
nonzero vectors of $\LL$ project to nonzero vectors in 
$\Vphys$; the last assumption is equivalent to requiring that 
$$
\LL \notin \mathcal{N} \df \{\LL' \in H  \LL_1 : \LL' \cap \Vint \not \subset
\{0\} \}.
$$
The condition that $\LL$ satisfies {\bf (I)} is valid for $\bar
\mu$-a.e.\  $\LL$ by definition of an RMS measure. We claim further that in the linear
case $\bar \mu
(\mathcal{N})=0$.
 Indeed, since $\bar \mu$ is induced by the Haar
measure of $H$, otherwise we would have
some fixed $v \in \LL_1 \sm \{0\}$ such that $H^{\mathcal{N}, v} \df \{h \in
 H:  h  v \in \Vint\}$ has positive Haar measure.
Recall
that for analytic varieties $\mathcal{V}_1, \mathcal{V}_2$, with
$\mathcal{V}_1$ connected, if
$\mathcal{V}_1 \cap \mathcal{V}_2$ has positive measure with respect
to the smooth measure on $\mathcal{V}_1$, then $\mathcal{V}_1 \subset
\mathcal{V}_2$. Since
$H^{\mathcal{N}, v}$ is an analytic subvariety in $H$, if 
it has positive measure with respect to the Haar measure on $H$, it
must coincide with $H$. This contradicts 
Lemma  \ref{lem: dense orbits H'}. This contradiction shows that
$\bar\mu$-almost surely we have $\widehat f \circ \Psi = \widehat
F$. Since $\mu = \Psi_* \bar \mu$, the first assertion that $\hat f
\in L^p(\mu)$ for $p < q_\mu$ now follows from
the first assertion of Theorem \ref{thm: reduction theory more
  general}. 

For the second assertion, let $f $ be a nonnegative
continuous function whose support contains a ball around the
origin. Since we have assumed that $W$ contains a ball around the
origin in $\Vint$, the support of the function $F$ 
also contains a ball around the origin in $\R^n$, so $\hat f$ is
bounded below by the Siegel-Veech 
transform of the indicator of a ball in $\R^n$, and we have that such
functions do not belong to $L^{q_\mu}(\bar \mu)$. 
\end{proof}

  \section{Integral formulas for the Siegel-Veech transform}\name{sec:  Weil}
 In this section we will prove Theorem \ref{thm: weil}. We begin with
 its special case $p=1$, i.e., with a derivation of \equ{eq: Siegel
   summation MS}. This will illustrate the method of Weil \cite{Weil}
 which we will use. Note that \equ{eq: Siegel
   summation MS} was
  first proved by Marklof and Str\"ombergsson in \cite{MS} following
  an argument of Veech 
  \cite{Veech_siegel_measures}. Their argument does not rely on an
  integrability bound such as our Theorem \ref{thm: reduction theory},
  and instead, uses the result of  
  Shah \cite{Shah}, Theorem~\ref{thm: Shah}. 
    
  \subsection{A derivation of a `Siegel summation
    formula'}\name{subsec: derivation}
  Given $f \in C_c(\R^d)$,
  define $F$ via \equ{eq: if define 1}, and define $\widehat{F}(\LL)$ via \equ{eq: if define 2}. 
  We can bound $F$ pointwise from above by a compactly supported
  continuous function on $\R^n$, and hence, 
  by Theorem \ref{thm: reduction theory more general}, $\widehat{F}
  \in L^1(\bar \mu)$. Therefore
$f \mapsto \int_{\mathbf{X}_1} \widehat F \, d \bar \mu$ is a positive 
linear functional on $C_c(\R^d)$. By the Riesz representation theorem,
there is some Radon Borel measure $\nu$ on $\R^d$ such that
$\int_{\mathbf{X}_1} \widehat F \,d \bar \mu = \int_{\R^d} f \, d\nu$.
From the equivariance relation
\equ{eq: equivariance Psi}, $\nu$ is invariant under
$\ASL_d(\R)$ in the affine case and under $\SL_d(\R)$ in the linear
case. Lebesgue measure is the unique (up to scaling) $\ASL_d(\R)$-invariant
Radon Borel measure on $\R^d$, and for $\SL_d(\R)$, the only
additional invariant measure is $\delta_0$, the Dirac mass at the
origin. Thus, there
are constants $c_1, c_2$ such that
\eq{eq: establishing}{
\nu = \left\{ \begin{matrix}  c_1 \Leb & \ \ \ \  \bar \mu \text{ is
      affine} \\ c_1 \Leb+ c_2 \delta_0
  & \ \ \ \ \bar \mu \text{ is linear.} \end{matrix} \right.
}
As we have seen in the proof of Theorem \ref{thm: reduction theory},
we have that $\widehat{F}= \widehat{f} \circ \Psi$ holds
$\bar \mu$-a.e.
Since $\mu = \Psi_* \bar \mu$, this implies that
$$
\int_{\Cl(\R^d)} \widehat{f} \, d\mu=\int_{\mathbf{X}_1} \widehat F \, d\bar \mu =
\int_{\R^d} f \, d\nu .
$$
In combination with \equ{eq: establishing}, this establishes \equ{eq:
  Siegel summation MS} in the affine case, 
and gives
\eq{eq: this gives}{
\int_{\Cl(\R^d)} \hat f \, d\mu = c_1 \int_{\R^d} f \, d\vol + c_2 
f(0), \ \ \forall f \in C_c(\R^d) }
in the linear case. It remains to show that $c_2=0$. 


Let $B_r = B(0,r)$ be the ball in $\R^d$ centered at the origin, 
let $f \in C_c(\R^d)$ satisfy $\mathbf{1}_{B_1} \leq f \leq
\mathbf{1}_{B_2}$, and let $f_r = f\left( \frac{x}{r}\right)$. Thus,	
as $r \to 0$, the functions $f_r$ have smaller and smaller support
around the origin.
%
By \equ{eq: Siegel Veech transform} and discreteness of $\Lambda$ we
have that $\widehat{f}_r 
(\Lambda)\to_{r \to 0} 0$ for any 
$\Lambda$. The functions $f_r$ vanish outside the ball $B_{2r}$, and
for $r \leq 1$, the functions $\widehat{f}_r$ are dominated 
by $\widehat{f}_1$. Therefore
$$
0=\lim_{r \to 0}  \int_{\Cl(\R^d)} \widehat{f}_r \, d\mu
\stackrel{\equ{eq: this gives}}{=} \lim_{r\to 0} \left[
    c_1\int_{\R^d} f_r \, d\vol +c_2 \cdot 1\right]= c_2.
  $$
\qed

  \subsection{A formula following Siegel-Weil-Rogers }\name{subsec:
    Weil}
  In this section we state and prove a generalization of Theorem
  \ref{thm: weil}. Let the notation be as in \ref{thm: MS
    classification}, so that $\bar \mu$ is an $H$-homogeneous measure
  on $\ALN_n$. Let $p \in \N$ and let $\R^{np} = \underbrace{
  \R^n \oplus \cdots \oplus \R^n}_{p \text{ copies}}$. For $f \in
C_c(\R^{np})$ and $\LL \in \ALN_n$, define
\eq{eq: the sum}{
 \fphat
(\LL) \df \left\{ \begin{matrix} \displaystyle{\sum_{v_1, \ldots, v_p \in
      \LL \sm \{0\}} f(v_1,
\ldots, v_p) }& \bar \mu \text{ is linear} \\ \displaystyle{\sum_{v_1, \ldots, v_p \in \LL} f(v_1,
\ldots, v_p)} & \bar \mu \text{ is affine.} \end{matrix} \right. 
}
Let $J \subset \ASL(np, \R)$ be a real algebraic group and let
$\theta$ be a  locally finite Borel measure on $\R^{np}$. We say
that $\theta$ is
{\em $J$-algebraic} if $J$
preserves $\theta$ and has an orbit of full $\theta$-measure
(in this
case $\theta$ can be described in terms of the Haar measure of $J$, see
\cite[statement and proof of Lemma 1.4]{Raghunathan}). 
  \begin{thm}\name{thm: weil more general}
Let $p \in \N$ and assume that $ p
 <q_{\bar \mu}$ where $q_{\bar \mu}$ is as in \equ{eq: p 
    bounded}. Then there is a countable collection $\{\bar{\tau}_{\ee} : \ee
  \in \E\}$ of $H$-algebraic Borel measures on $\R^{np}$
such that 
$\bar \tau \df \sum \bar \tau_\ee$ is locally finite and 
for every $f \in
L^1(\bar \tau )$ we have
\eq{eq: general summation formula}{
\int_{\ALN_n} \fphat \,
d\bar \mu = \int_{\R^{np}} f \, d\bar{\tau} .
}
    \end{thm}
As we will see in the proof, in the affine (resp. linear) case, the
indexing set $\E$ is naturally identified with the set of  
$\Gamma_{H_1}$-orbits in the set of $p$-tuples of (nonzero) vectors in
$\Z^{n}$.

We will need a by-now standard
result of Weil, which is a generalization of the Siegel summation
formula and is proved via an argument similar to the one used in \S
\ref{subsec: derivation}.
Let $G_1 \subset G_2$ be unimodular locally compact groups, let $\Gamma_2 \subset
G_2$ be a lattice in $G_2$ and let $m_{G_2/\Gamma_2}$ denote the
unique $G_2$-invariant Borel probability measure on
$G_2/\Gamma_2$. Since $G_1, G_2$ are unimodular, there is
a unique (up to scaling) locally finite $G_2$-invariant measure on
$G_2/G_1$, which we denote by $m_{G_2/G_1}$ (see
e.g. \cite[Chap. I]{Raghunathan}). Define $\Gamma_1 \df \Gamma_2
\cap G_1$, and for any $\gamma \in 
\Gamma_2$, denote its coset $\gamma \Gamma_1 \in \Gamma_2/\Gamma_1$
by $[\gamma]$. 
With this notation, Weil showed the following:

\begin{prop}[\cite{Weil_brazil}]\name{prop: Weil lemma}
Assume that $\Gamma_1$ is a lattice in $G_1$.
Then we can rescale $m_{G_2/G_1}$ so that the following holds.  For any $F \in
L^1(G_2/G_1, m_{G_2/G_1})$, define 
\eq{eq: formula transform weil}{
  \til F (g\Gamma_2) \df \sum_{[\gamma] \in \Gamma_2/\Gamma_1}
F(g\gamma).}
Then $\til F \in L^1(G_2/\Gamma_2, m_{G_2/\Gamma_2})$ and
$$
\int_{G_2/\Gamma_2} \til F \, dm_{G_2/\Gamma_2} = \int_{G_2/G_1} F
\, dm_{G_2/G_1}.
$$
\end{prop}

\begin{proof}[Proof of Theorem \ref{thm: weil more general}]
Consider the map which sends $f \in C_c(\R^{np})$ to $\int \fphat \,
d\bar \mu$. This is well-defined by
Theorem \ref{thm: reduction theory more general}, and defines a positive linear
functional on $C_c(\R^{np})$. Thus, by the Riesz representation
theorem, there is a locally finite measure $\bar{\tau}$ on $\R^{np}$ such
that
\eq{eq: locally finite measure}{
\forall f \in C_c(\R^{np}), \ \ \int_{\ALN_n} \fphat \,
d\bar \mu = \int_{\R^{np}} f \, d\bar{\tau}.
}
Our goal will be to present $\bar{\tau}$ as a countable linear combination
of $H$-algebraic measures. Note that since $C_c(\R^{np})$ is a dense
linear subspace of $L^1(\bar \tau)$, for any locally finite measure
$\bar \tau$,
it suffices to prove \equ{eq: general
  summation formula}  for functions in 
$C_c(\R^{np})$. 

Let $
H, \, g_1, \, \LL_1 = g_1\Z^n, \,
H_1 = g_1^{-1} H g_1
, \, \Gamma_{H_1} = H_1 \cap
\ASL_n(\Z)$ be as in \S\ref{subsec: homogeneous measures}, so that
$\Gamma_{H_1}$ is a lattice in $H_1$ and 
$\bar \mu$ is an $H$-homogeneous measure supported on $H  \LL_1
\cong H_1/\Gamma_{H_1}$. In the affine (respectively
linear) case, let $\Z^{np}$ denote the 
countable collection of ordered 
$p$-tuples of vectors in $\Z^n$ (respectively, in $\Z^n \sm
\{0\}$). Let $\E$ denote the collection of $\Gamma_{H_1}$-orbits in
$\Z^{np}$. 
For each $\ee \in \E$, define the restriction of the
sum \equ{eq: the sum} to the orbits $H  \LL_1 \subset \ALN_n$
and to the orbit $\ee$, by  
\eq{eq: orbit sum e}{
\fphat_\ee(h \LL_1) \df \sum_{(x_1, \ldots, x_p) \in
 \ee}
f\left(hg_1x_1, \ldots, hg_1x_p \right),
}
so that on $H \LL_1$ we have
\eq{eq: sum all tau}{
\fphat = \sum_{\ee \in \E} \fphat_\ee.
}
If $f$ is a non-negative function then $\fphat_\ee  \leq
\fphat$ everywhere on $H  \LL_1$, and in particular
$\fphat_\ee \in L^1(\bar \mu)$. Thus, the assignment sending $f \in
C_c(\R^{np})$ to 
\eq{eq: corresponding bar tau}{
\int f \, d\bar \tau_\ee \df \int \fphat_\ee \, d\bar \mu
}
 is a 
positive linear functional and hence, via the Riesz representation
theorem, defines 
the locally finite Borel
measure $\bar{\tau}_\ee$ on $\R^{np}$. By \equ{eq: sum all tau},
$\sum_{\ee \in \E} \bar{\tau}_\ee = 
\bar{\tau}.$ It remains to show that each $\bar{\tau}_\ee$ is $
H$-algebraic.

For each $\ee \in \E$, choose a representative $p$-tuple $\vec x _\ee =
(x_1, \ldots, x_p) \in \ee$ and let 
$$G_{1,\ee} \df \{h \in H_1: h x_i = x_i, \,
i=1, \ldots, p\}.$$
We will apply Proposition \ref{prop: Weil
  lemma} with $G_2 = H_1, \, \Gamma_2 = \Gamma_{H_1}, \,  
G_1 = G_{1, \ee}, \, \Gamma_1 = \Gamma_2 \cap G_1$, and with $F
(h_1G_1 ) \df  f(g_1h_1\vec{x}_\ee)$. Comparing \equ{eq: formula
  transform weil} and \equ{eq:
  orbit sum e} we
see that these choices imply that $\til
F(h_1\Gamma_2) = \fphat_\ee(h\LL_1)$, for $h=g_1 h_1 g_1^{-1} \in H$. We will see below that 
$\Gamma_1$  is a lattice in $G_{1}$. Assuming this, 
we apply Proposition \ref{prop: Weil lemma} to obtain
\[\begin{split}
& \int_{\R^{np}} f \, d\bar{\tau}_\ee = \int_{\ALN_n} \fphat_\ee \, d \bar \mu
= \int_{G_2/\Gamma_2} \til F \, dm_{G_2/\Gamma_2} \\ = & \int_{G_2/G_1}
f(g_1h_1 \vec{x}_\ee) \,
dm_{G_2/G_1} (h_1G_1).
\end{split}\]
This shows that $\bar{\tau}_\ee$ is the pushforward of $m_{G_2/G_1}$ under the map
$$
G_2/G_1 \to \R^{np}, \ \ \ h_1G_1 \mapsto g_1h_1 \vec{x}_\ee.
$$
In particular, since $H = g_1 H_1 g_1^{-1}$, $\bar{\tau}_\ee$ is
$H$-algebraic.

It remains to show that $\Gamma_1$ is a lattice in $G_1$. 
To see this, note that $G_2$ is a real algebraic
 group defined over $\Q$, and $G_1$ is 
 the stabilizer in $G_2$ of a finite collection of 
  vectors in $\Z^n$. Thus, $G_1$  is also defined over $\Q$. 
By the theorem of Borel and Harish-Chandra (see
\cite[\S 13]{Borel_arithmetiques}), if
$G_1$ has no nontrivial characters then 
$\Gamma_1 = G_1 \cap \ASL_n(\Z)$ is a
  lattice in $G_1$. Moreover, a real algebraic group generated by unipotents has no
  characters. Thus, to conclude the proof of the claim, it suffices to
  show that $G_1$ is generated by unipotents. We verify this by
  dividing into the various cases arising in Theorem \ref{thm: MS classification}. 

  We first reduce to the case that $G_1$ is a subgroup of $\SL_n(\R)$.
In the linear case we simply identify $G_2$ with its isomorphic image
$\pi(G_2)$, where $\pi: 
\ASL_n(\R) \to \SL_n(\R)$ is
the projection in \equ{eq: def projection}, and thus we can assume
$G_1 \subset \SL_n(\R)$. In the affine 
case, since the property of being generated by unipotents is invariant
under conjugations in $\ASL_n(\R)$, we may conjugate by a
translation to assume that one of the vectors in $\vec{x}_\ee$ is the
zero vector, so that $G_1 \subset \SL_n(\R)$. Thus, in both cases we
may assume that $G_2 = \underline H_1 $ is the group of real points of
$\Res_{\KK/\Q}(\mathbf{G})$, and $G_1$ is the stabilizer in $G_2$ of
the finite collection $x_1, \ldots, x_p$, where these are vectors in
the standard representation on 
$\R^n$. 

Suppose first that $\mathbf{G} = \SL_k$. Then, in the notation of
\equ{eq: shape of restriction}, we have that $G_2 = 
{}^{\sigma_1}\mathbf{G}_\R \times \cdots
\times{}^{\sigma_{r+s}}\mathbf{G}_\R$, 
where for $i=1, \ldots, i$ (respectively, for $i=r+1, \ldots, r+s$) we
have that ${}^{\sigma_i}\mathbf{G}_\R$ is isomorphic to $\SL_k(\R)$
(respectively to $\SL_k(\CC)$ as a real algebraic group). Furthermore,
as in \S \ref{subsec: number fields}, there is a
decomposition 
$$\R^n = V_1 \oplus \cdots \oplus V_{r+s},$$ 
where $V_i \cong \R^{k}$ (resp., $V_i \cong \R^{2k}$) for $i=1,
\ldots, r$ (resp., for $i=r+1, \ldots, r+s$), and such that the action of
$G_2$ on $\R^n$ is the product of the standard action of each
${}^{\sigma_i}\mathbf{G}_{\R}$ on $V_i$. Let $P_i: \R^n \to V_i$ be the
projection with respect to this direct sum decomposition. Then the
stabilizer in $G_2$ of $x_1, \ldots, x_p$ is the direct product of the
stabilizer, in ${}^{\sigma_i}\mathbf{G}_{\R}$, of $P_i(x_1), \ldots,
P_i(x_p)$. So it suffices to show that each of these stabilizers is
generated by unipotents. In other words, we are reduced to the
well-known fact that for $\SL_k(\R)$ acting on $\R^k$ in the standard
action, and for $\SL_k(\CC)$ acting on $\R^{2k} \simeq \CC^k$ in the
standard action, the stabilizer of a finite collection of vectors is
generated by unipotents. 

Now suppose that $\mathbf{G} = \Symp_{2k}$, and let $\mathbb{F} = \R$
or $\mathbb{F}= \CC$. Then by a similar
argument, we are reduced to the statement that for the standard action
of $\Symp_{2k}(\mathbb{F})$ on $\mathbb{F}^{2k}$, the stabilizer of a finite collection
of vectors is generated by unipotents. This can be shown as
follows. Let $\omega$ be the symplectic form preserved by $\Symp_{2k}$,
let $V = \spa(x_1, \ldots, x_p) \subset \mathbb{F}^{2k}$, and let
$$Q \df \{g
\in \Symp_{2k}(\mathbb{F}) : \forall v \in V, \, gv =v \}.$$
We need to show that $Q$ 
is generated by unipotents. We can write $V = V_0 \oplus V_1$, where
$V_0 = \ker \left(\omega|_V \right)$ is Lagrangian, and $V_1$ is
symplectic. Let $2\ell = \dim V_1$, where $\ell \leq k$. Since any
element of $ Q$ fixes $V_1$ pointwise, it leaves 
$V_1^\perp$ invariant, and it also fixes pointwise the subspace $V_0
\subset V_1^\perp$. Thus, $Q$ is isomorphic to 
$$
 \{g \in \Symp(V_1^\perp): \forall v
\in V_0, \, gv=v\} \subset \Symp(V_1^\perp) \cong \Symp_{2m}(\mathbb{F}),  
$$
 where $ m \df k - \ell$.
This means we can reduce the problem to the case in which $V_1 =
\{0\}$, i.e., $\omega(x_i, x_j)=0$ for all $i,j$. We can apply a
symplectic version of the Gram-Schmidt orthogonalization procedure to
assume that $x_1, y_1,
\ldots, x_p, y_p, x_{p+1}, y_{p+1}, \ldots, x_m, y_m$ is a symplectic
basis and $V_0 = \spa (x_1, \ldots, x_p)$. Let
$$V_2 \df \spa (x_{p+1}, y_{p+1}, \ldots, x_m, y_m)   \ \text{ and } \ 
V_3 \df V_0
\oplus V_2. $$
Then $V_2$ is symplectic and the subgroup of $Q$ leaving $V_2$ invariant is
isomorphic to $\Symp_{2m-2p}(\mathbb{F})$, hence generated by
unipotents. 
Also, for $i=1, \ldots, p$, by considering the
identity
$$\omega(gy_i, x_j) = \omega(gy_i, gx_j) = \omega (y_i, x_j) \ \
(j=1, \ldots, p)$$ one sees that any $g \in Q$ must map
the $y_i$ to vectors in  $y_i+ V_3$. This implies that $Q$ is
generated by symplectic matrices leaving $V_2$ invariant, and transvections
mapping $y_i$ to 
elements of $y_i + V_3$. In particular, $Q$ is generated by unipotents.
\end{proof}

\begin{dfn}\name{def: c and p algebraic measure}
Given a real algebraic group $J \subset \ASL_n(\R)$, we will say that
a locally finite measure $\tau$ on $\R^{dp}$ is {\em
  $J$-c\&p-algebraic} if there is a $J$-algebraic measure $\bar \tau$ on  
$\R^{np}$ such that for every $f \in C_c(\R^{dp})$ we have
$$
\int_{\R^{dp}} f \, d\tau = \int_{\R^{np}} F \, d\bar \tau, 
$$
where
$
F: \R^{np} \to \R$ is defined by 
\eq{eq: c and p measure functions}{
  F(x_1,  \ldots, x_p) \df \left\{ \begin{matrix}
      f\left(\piphys(x_1), \ldots, \piphys(x_p) \right) & \forall i,
      \ \piint(x_i) \in W \\
    0 & \text{ otherwise .} \end{matrix}
  \right.
}
We will say $\tau$ is {\em c\&p algebraic} if it is $J$-c\&p
algebraic for some $J$. 
\end{dfn}

It is easy to check that for $p=1$, the measure $\tau$ in
Definition \ref{def: c and p algebraic measure} is the pushforward
under $\piphys$ of the restriction of $\bar \tau$ to
$\pi^{-1}_{\mathrm{int}}(W)$. For general $p$, define projections
$${}^{p}\piphys: \R^{np} \to \R^{dp},\ \ {}^{p}\piphys(x_1,
\ldots, x_p) \df \left(\piphys(x_1), \ldots, \piphys(x_p) 
\right),
$$
and
$$
{}^{p}\piint: \R^{np}\to\R^{mp}, \ \ {}^{p}\piint(x_1, \ldots, x_p)
\df \left(\piint(x_1), \ldots, \piint(x_p) \right).$$ 
Then the measures $\tau, \, \bar\tau$ satisfy
\eq{eq: the measures satisfy}{
\tau =
{}^{p}\piphys
_*\left(
\bar\tau|_{\mathcal{S}}
\right)
, \text{ where
} \mathcal{S} \df
{}^{p}\pi^{-1}_{\mathrm{int}}
\left(
\underbrace{W \times \cdots
  \times W}_{p \text{ copies }}
\right)
.
}
  \begin{proof}[Proof of Theorem \ref{thm: weil}]
By Theorem \ref{thm: MS surjective}, after a rescaling of $\R^d$,
there is a homogeneous measure 
$\bar \mu$ on $\ALN_n$ such that $\mu = \Psi_* \bar \mu$. Suppose $h
\in H$ satisfies that $\piphys|_{h\LL_1}$ is 
injective, and in the linear case, assume also that $h \LL_1 \cap
\Vint \subset \{0\}$. Since $\mu$ is an RMS measure, and in the linear
case, arguing as in the proof of Theorem \ref{thm: reduction theory} using Lemma
\ref{lem: dense orbits H'}, we see that this holds
for a.e.\  $h \in H$. For such 
$h$, letting $\Lambda_h \df \Psi(h\LL_1)$, 
we can rewrite the function $\fphat$
defined in \equ{eq: def fphat} more succinctly in the form
\[
  \fphat(\Lambda_h) =
  \sum_{(x_1, \ldots,
    x_p) \in \LL_1^p}
 F(hx_1,  \ldots, hx_p) ,
\]
where $F$ is as in \equ{eq: c and p measure functions}. Thus, Theorem
\ref{thm: weil} is reduced to Theorem \ref{thm: weil more general}.
\end{proof}

\begin{remark}
  The assignment $\ee \mapsto \bar \tau_\ee$ implicit in
  the proof of Theorem
  \ref{thm: weil} is not injective, nor is it finite-to-one. To see
  this, take $p=1$ and consider the RMS measure corresponding
  to the Haar-Siegel measure on $\LLN_n$. Then $H_1 = \SL_n(\R), \
  \Gamma_{H_1} = \SL_n(\Z)$, and there are countably many
  $\Gamma_{H_1}$-orbits on $\Z^n$, where two integer vectors belong to the
  same orbit if and only if the greatest common divisor of their
  coefficients is the same. On the other hand, as the proof of formula
  \equ{eq: Siegel summation MS} shows, there are two c\&p-algebraic
  measures, namely Lebesgue 
  measure on $\R^d$ and the Dirac measure at 0. The Dirac measure is
  associated with the orbit of $0 \in \Z^n$, and all the other orbits
  of nonzero vectors in $\Z^n$ give rise to multiples of Lebesgue
  measure on $\R^d$.

 Nevertheless, we will 
continue using the symbol $\E$ for both the  
collection of $\Gamma_{H_1}$-orbits in
  $\Z^{np}$, and for the indexing set for the
  countable collection of measure arising in Theorem \ref{thm:
    weil}. This should cause at most mild confusion.
    \end{remark}

%

  \section{The Rogers inequality on moments}\name{sec: Rogers}
In this section we will prove Theorem \ref{thm: Rogers bound RMS}.
We will need more information about the measures $\tau_\ee$ appearing
in Theorem \ref{thm: weil}, in case $p=2$. We begin our discussion with
some properties that are valid for all $p \leq d$. Some of the results
of \S \ref{subsec: normalization} will be given in a greater level of
generality than required for our counting
results. They are likely to be of use in understanding higher
moments for RMS measures.

    \subsection{Normalizing the measures}\name{subsec: normalization}
   For any $k$, denote the normalized Lebesgue measure on $\R^k$ by
   $\vol^{(k)}$.  Some of the c\&p-algebraic measures $\tau$ on $\R^{dp}$ which arise in
   Theorem \ref{thm: weil} are the globally supported Lebesgue measures on
   $\R^{dp}$, i.e., multiples of $\vol^{(dp)}$. Indeed, such a measure
   arises if in Definition
   \ref{def: c and p 
     algebraic measure} we take
   $\bar \tau$ equal to a multiple of Lebesgue measure on
   $\R^{np}$. 
These measures give a main term in the counting problem we will
consider in \S \ref{sec: counting}. We write $\tau_1 \Propto 
\tau_2$ if $\tau_1, \tau_2$ are proportional, recall the measures
$\{\tau_\ee\}$
defined 
in the
proofs of Theorems \ref{thm: weil} and \ref{thm: weil more general}, and set
$$ 
\E^{\mathrm{main}} \df \left \{\ee \in \E: \tau_\ee \Propto
  \vol^{(dp)} \right\}, \ \ \ \ \ 
\tau_{\mathrm{main}} \df \sum_{\ee \in 
  \E^{\mathrm{main}}} \tau_\ee.
$$
We define constants $c_{\mu, p}
$ by the condition
$$
\tau_{\mathrm{main}} = c_{\mu, p}\,
\vol^{(dp)}.$$
   
The next
   result identifies the normalizing constants 
   $c_{\mu, p}
   $.
 Recall from Theorem \ref{thm: MS surjective} that
   an RMS measure $\mu$ is of the form $\mu = \rho_{c*} \bar \mu$ where
   $\bar \mu$ is a homogeneous measure on $\ALN_n$, $c$ is the
   constant of \equ{eq: for scaling}, and $\mu$-a.e.\,$\Lambda$ is of the form $\Lambda
  = \Lambda(\LL, W)$ for a grid $\LL$ with $\covol(\LL) = c^n$.
We denote this almost-sure value of $\covol(\LL)$ by $\covol(\mu)$. 
   Recall also that the function $\Lambda \mapsto D(\Lambda)$
   defined in \equ{eq: density exists} is measurable and invariant,
   and hence is a.e.\,constant, and denote its almost-sure value by $D(\mu)$.

   \begin{prop}\name{prop: normalizing constant}
     For any RMS measure $\mu = \rho_{c*} \Psi_* \bar \mu$ satisfying
     \equ{eq: new condition} (i.e., $\mathbf{G} = \SL_k$ or $\mu$ is
     affine), we have
    \eq{eq: cmu1}{
    c_{\mu,1} = D(\mu) = \frac{\vol^{(m)}(W)}{\covol(\mu)},} 
  and  for $p  \in \N$ 
  satisfying $p<q_\mu$ and
   $  p \leq d
  $ 
  we have 
  \eq{eq: cmuhigher}{
    c_{\mu,p} = c_{ \mu,1}^p.}
\end{prop}
Note that the normalizing constant $c_{\mu,1}$ discussed here is the
same as the constant 
denoted by $c_1$ in \equ{eq: this gives} and by $c$ in \equ{eq: Siegel
  summation MS}. 

With the identification $\R^{\ell p}
\cong M_{\ell, p}(\R)$ in mind, we say that a subspace $V
\subset \R^{\ell p}$ is an {\em annihilator subspace} if it is the common
annihilator of a collection of vectors in $\R^p$; that is, there is a collection
$\mathrm{Ann} \subset \R^p$ such that 
\[\begin{split}
& V  \, =  \mathcal{Z}(\mathrm{Ann}) \\ & \df  \left\{ (v_1, \ldots, v_p) \in \R^{\ell
    p}: \forall i, \, v_i \in \R^\ell \, \& \, 
  \forall (a_1, \ldots, a_p) \in 
\mathrm{Ann}, \ \sum a_i v_i =0  \right \}. 
\end{split}\]
Note that the meaning of $\mathcal{Z}(\mathrm{Ann})$ depends on the
choice of the ambient space $\R^\ell$ containing the vectors $v_i$;
when confusion may arise we will specify 
the ambient space explicitly.

Suppose $\ell \in \N$ and
$(v_1, \ldots, v_p)$ is a
$p$-tuple in $\R^{\ell p}$. In the linear case, let 
$$\mathrm{Ann} (v_1, \ldots, v_p) \df \{(a_1, \ldots, a_p) \in \R^p:
\sum a_i v_i =0\},$$
and in the affine case, let
$$
\mathrm{Ann} (v_1, \ldots, v_p) \df \{(a_1, \ldots, a_{p-1}) \in
\R^{p-1} : \sum a_i (v_{i}-v_p)=0\}. 
$$
Let
$$
L(v_1, \ldots, v_p) \df \mathcal{Z}(\mathrm{Ann}(v_1, \ldots, v_p)),
$$
an annihilator subspace in $\R^{\ell p}$, 
say that $v_1, \ldots, v_p$ are {\em independent} if
$\mathrm{Ann}(v_1, \ldots, v_p) = \{0\}$, and let
$$\rank(v_1, \ldots,
v_p) \df \left\{ \begin{matrix} p - \dim \mathrm{Ann}(v_1, \ldots, 
v_p) & \ \ \ \mu \text{ is linear} \\ p - 1- \dim \mathrm{Ann}(v_1, \ldots, 
v_p)  & \ \ \  \mu \text{ is affine}. \end{matrix} \right. $$
Note that in the linear case, this is the usual relation between the
rank of a matrix and the dimension of its kernel. 
The dimension of $L(v_1, \ldots, v_p)$ is equal to $\ell \, \rank(v_1,
\ldots, v_p)$.

We recall some notation from \S \ref{subsec: number fields} and from
Step 3 of the proof of Lemma \ref{lem:classification}.
Let $\KK$ be a real
number  field of degree $D = r+2s$, with 
$\sigma_1, \ldots, \sigma_r$ being distinct real embeddings, and
$\sigma_{r+1}, \ldots, \sigma_{s}$ denoting representatives of
conjugate pairs of non-real embeddings. Let $\mathbf{G}$
be isomorphic to either $\SL_k(\R)$ or to $\Symp_{2k}(\R)$, and let
$\mathbf{H} = 
\Res_{\KK/\Q}(\mathbf{G}).$  Let $\mathbf{V}$ be a
$\KK$-vector space of dimension $t$, where $t$ is as in \eqref{eq: def
  t}, and denote
$V_j = {}^{\sigma_j}\mathbf{V}_{\R}$, that is, $V_j \cong \R^t$ if $j=1,
\ldots, r$ and $V_j \cong 
\CC^t \cong \R^{2t}$ if $j=r+1, \ldots, s$. These vector spaces are
chosen so that $\mathbf{V}$ is equipped with the standard action of
$\mathbf{G}$, and taking into account the isomorphism 
\eq{eq: the isomorphism}{
  \R^n \cong \left(\Res_{\KK/\Q}(\mathbf{V})\right)_{\R} =V_1 \oplus \cdots \oplus
V_{r+s}.}
Let ${}^{\sigma_j}\pi: \R^n \to
V_j$ be the corresponding projections. In the notation
\equ{eq: shape of restriction}, let $\pi^j : \mathbf{H}_{\R} \to
{}^{\sigma_j}\mathbf{G}_{\R}$, so that the action of $\mathbf{H}_{\R}$
factors through the action of each ${}^{\sigma_j}\mathbf{G}_{\R}$ on
$V_j$.  We can assume without 
loss of generality (see \S \ref{subsec: c and p sets}) that $V_2 
\oplus \cdots \oplus V_{r+s} \subset \Vint$ and
$\piphys = \piphys \circ {}^{\sigma_1}\pi$.

\begin{lem}\name{lem: for small p} Suppose $\mu$ is an RMS measure of
  higher rank, and let $\mathbf{G}$ be the group appearing in Theorem
  \ref{thm: classification}.
  Let $p <
  q_\mu$, let $\vec{x}_\ee = (x_1, \ldots, x_p) \in \ee$, where $\ee
  \in \E$ is
  as defined before \equ{eq: orbit sum e}, and let $v_i
  \df {}^{\sigma_1}\pi(x_i)$, $i=1, \ldots, p$. Assume that
  $$\rank(v_1, \ldots, v_p) \leq 
  \left \{ \begin{matrix} d & \text{ if } \mathbf{G} =
   \SL_k \\ 1 & \ \, \text{ if } \mathbf{G} =
     \Symp_{2k}. \end{matrix} \right.
  $$
  Let $\bar
  \tau \df \bar
  \tau_\ee$ be the algebraic measure on $\R^{np}$ as in \equ{eq:
    corresponding bar tau}  and
  let $\tau$ be a
c\&p-algebraic measure obtained from $\bar \tau $ as in Definition
\ref{def: c and p algebraic measure}. Then $\tau$ is (up to proportionality) the Lebesgue 
measure on some annihilator subspace of $\R^{dp}.$ This subspace is
equal to $\R^{dp}$ 
if and only if $v_1, \ldots, v_p$ are independent.
  \end{lem}

 \begin{proof}
Let
$\bar \tau$ be as in Definition \ref{def: c and p algebraic  
measure}.
   As in the proof of Theorem \ref{thm: weil more general}),  we have
   that $H(x_1, \ldots, x_p)$ is a dense 
subset of full measure in $\supp \, \bar \tau$. 
We will split the proof 
according to the
various cases arising in Theorem \ref{thm: MS classification}.

\medskip

{\em Case 1: $\mu$ is linear, $\mathbf{G} = \SL_k.$}
In this case, our proof will also show that $\supp \, \bar \tau$ is a
sum of annihilator subspaces, one in each $V_j$; in fact, we first
establish this statement. 

The action of
$H$ on $\R^n$ factors into a product of actions of each 
${}^{\sigma_j}\mathbf{G}_{\R}$ on $V_j$. That is, $H$
acts on $v^j_i \df {}^{\sigma_j}\pi(x_i), \, i=1, \ldots, p$ via its mapping to
${}^{\sigma_j}\mathbf{G}_\R$, i.e., via the standard action of
$\SL_k(\R)$ or $\SL_k(\CC)$ on $\R^k$ or $\CC^k$. It follows from
\equ{eq: p 
  bounded} and \equ{eq: K rank} that
$p< q_\mu=k$. Therefore for each $j$, the rank $R_j$ of $\left\{v^j_i :  i=1,
\ldots, p \right\}$ is less than $k$. For the 
standard action, ${}^{\sigma_j}\mathbf{G}_{\R}$ is transitive on
linearly independent 
$R_j$-tuples. From this, by choosing a linearly independent subset $B_j
\subset \{v^j_1, \ldots, v^j_p\}$ of cardinality $R_j$
and expressing 
any $v^j_i \notin B_j$  as a linear
combination of elements of $B_j$,
one sees that if $(u_1, \ldots, u_p),
(w_1, \ldots, w_p )$ 
are two $p$-tuples in $V_j$ 
\eq{eq: split equivalence}{
  \begin{split}
& \text{ there is } h \in {}^{\sigma_j}\mathbf{G}_\R \text{ such that } h(w_1, \ldots, w_p)
= (u_1, \ldots, u_p) \\
\iff & \mathrm{Ann}(w_1, \ldots, w_p) = \mathrm{Ann}(u_1, \ldots, u_p).
\end{split}
}
This implies that ${}^{\sigma_j}\mathbf{G}_{\R}(v^j_1, \ldots, v^j_p)$ is
open and dense in $L(v^j_1, \ldots, v^j_p)$,
and hence $H(x_1, \ldots,
x_p)$ is open and dense in $L_1^{r+s} \df
\bigoplus_{j=1}^{r+s}L(v^j_1, \ldots, v^j_p)$. 
We have shown that
$\supp \, \bar \tau = L_1^{r+s}$ 
and that $\bar \tau$ is a 
multiple of the Lebesgue measure on $L_1^{r+s}$. 

Since $\piphys = \piphys \circ {}^{\sigma_1}\pi$, we have 
$${}^p\piphys\left(L_1^{r+s} \right)
= {}^p\piphys\left(L(v_1, \ldots, v_p) \right).$$
To simplify notation, write  $H^1 \df
{}^{\sigma_1}\mathbf{G}_\R \cong \SL_k(\R)$, and $v_i \df v^1_i \in V_1$. Let
$$
\mathrm{Ann}_1 \df \mathrm{Ann} \left(v_1, \ldots, v_p \right).
$$
We have 
\eq{eq: now show that}{
  {}^p\piphys(L(v_1, \ldots, v_p)) = \mathcal{Z}(\mathrm{Ann}_1), 
}
seen as an annihilator subspace of $\R^{dp}$. 
%
%
Indeed, the inclusion $\subset$ follows from linearity of
$\piphys$. For the opposite inclusion, recall that we have an
inclusion $\Vphys \hookrightarrow
V_1$, and this induces an inclusion
$\iota: \R^{dp} \hookrightarrow \R^{np}.$
We clearly have
$$\iota \left
(\mathcal{Z} \left(\mathrm{Ann}_1\right)\right) \subset L(v_1, \ldots,
v_p),$$
which implies the inclusion $\supset$ in \eqref{eq: now show that}.

Replacing $x_i$ with elements of $x_i + \Vphys$ does not change the
condition $(x_1, \ldots, x_p) \in \mathcal{S}$, where $\mathcal{S}$ is
as in \equ{eq:
  the measures satisfy}. This shows that
$$\supp \, \tau = {}^p\piphys
\left( L_1^{r+s} \right) =  {}^p\piphys \left
  (L(v_1, \ldots, v_p) \right)$$
is an annihilator subspace, and $\tau$ is a
multiple of Lebesgue measure on this subspace. Moreover, the subspace
is proper if and only if $\mathrm{Ann}_1 \neq \{0\}$, or equivalently,
$v_1, \ldots, v_p$ are  dependent.

\medskip

{\em Case 2: $\mu$ is linear, $\mathbf{G} = \Symp_{2k}, \ d=2.$}
The action of $H$ splits as a Cartesian product of actions of the
groups ${}^{\sigma_j}\mathbf{G}_\R$ on the spaces
$V_j$, for $j=1, \ldots, r+s$. As in Case 1, we will pay attention
to the action on the first summand $V_1 $,
where $H$ acts via $H^1 \df {}^{\sigma_1}\mathbf{G}_\R \cong
\Symp_{2k}(\R)$. We denote by $\omega$ the symplectic form on $V_1$
preserved by $H^1$. 
Let $L \df  \overline{H(x_1, \ldots, x_p)}  =
  \supp \, \bar \tau,$ where $\bar \tau$ is the unique (up to scaling)
  $H$-invariant
  measure with support $L$, and let 
  $L^1 \df L \cap V_1= {}^{\sigma_1}\pi(L) = \overline{H^1\left( v_1, \ldots,
      v_p\right)}
  $, where $v_i \df {}^{\sigma_1}\pi(x_i), \, i=1, \ldots, p$. 
%
%

Let $F \cong \SL_2(\R)$ be as in \equ{eq: def G0}. Then $F \subset
H^1$, and hence $\tau$  is $F$-invariant. Write
$$V^1_{\mathrm{int}} \df \Vint \cap V_1 = V^\perp_{\mathrm{phys}},$$
and abusing
notation slightly, let
$\piphys, \piint$ denote the restrictions of these mappings to $V_1$,
so they are the projections associated with the direct sum
decomposition $V_1 = \Vphys \oplus V^1_{\mathrm{int}}.$ Define $R \df
\rank(v_1, \dots, v_p)$, and define $R'$ as the maximal rank of 
$\{\piphys(hv_1), \ldots, \piphys(hv_p)\}$, as $h$ ranges over
elements of $H^1$. Thus, $0 \leq R' \leq R \leq 1$.

If $R'=0$ this means that $\piphys(hv_i)=0$ for all $h \in H$ and all
$i$, and 
then $\tau$ is the Dirac measure at 0, and there is nothing to
prove. Now suppose $R'=R=1$.
Since $R=1,$ there is some $v_i$ such that
$\piphys(v_i) \neq 0$, and there are coefficients $a_j,\, j \neq i$
so that $v_j = a_j v_i$. 
This implies that for all $h$, $\piphys(hv_j) = a_j \piphys(hv_i)$, that is,
$$
\supp \, \tau \subset {}^p\piphys(L) \subset L' \df \{(u_1, \ldots, u_p) \in \R^{2p} : \forall j
\neq i, \, u_j = a_ju_i\}.
$$
Moreover, since $F$ acts transitively on nonzero vectors in $\Vphys$,
and $\tau$ is $F$-invariant, we actually have equality and $\tau$ is a
multiple of Lebesgue 
measure on the annihilator subspace $L'$, and $L'$ is a proper
subspace of $\R^{2p}$, unless $p=1$.
\medskip

{\em Case 3: $\mu$ is affine.}
The affine case can be reduced to the linear case. Note that the
definition of the annihilator $\mathrm{Ann}(v_1, \ldots, v_p)$ in the
affine case is such that it does not change under the diagonal action
of the group of translations, and that the group of translations in
$H$ is the full group $\R^n$, so that $x_1, 
\ldots, x_p$ can be moved so that $x_p=0$.
Moreover, by Proposition \ref{prop: origin in window}, we can assume
that $0 \in W$. We leave the details to
the diligent reader. 
\end{proof}

Let
\eq{eq: def Erest}{
  \E^{\mathrm{rest}} \df \E \sm \E^{\mathrm{main}}, \ \ \
\tau_{\mathrm{rest}} \df \sum_{\ee \in \E^{\mathrm{rest}}} \tau_\ee.}
The preceding discussion gives a description of the measures $\tau_\ee$
with $\ee \in \E^{\mathrm{rest}}.$ 
\begin{cor}\name{cor: for small p}
Under the conditions of Lemma \ref{lem: for small p}, any measure $\tau_\ee, \, \ee
\in \E^{\mathrm{rest}}$, is Lebesgue measure on a proper subspace of
$\R^{dp}$. 
\end{cor}

\begin{proof}[Proof of Proposition \ref{prop: normalizing constant}]

Let $B_r$ denote the Euclidean ball of radius $r$
around the origin in $\R^d$, let  $\mathbf{1}_{B_r}$ be its
indicator function, and let 
  $\widehat{\mathbf{1}_{B_r}}$ be 
  the function obtained from the summation formula \equ{eq: Siegel
    Veech transform}, so that
  $$
D(\Lambda) = \lim_{r \to \infty}
\frac{\widehat{\mathbf{1}_{B_r}}(\Lambda) }{\vol^{(d)}(B_r)}. 
  $$
Applying \equ{eq: Siegel
  summation MS} we get that for any $r>0$, 
\eq{eq: inside the integral}{\int_{\Cl(\R^d)}
 \frac{\widehat{\mathbf{1}_{B_r}}}{\vol^{(d)}(B_r)} \, d\mu = 
 \frac{c_{\mu,1}}{\vol^{(d)}(B_r)} \int_{\R^d} \mathbf{1}_{B_r} \, d\vol = c_{\mu,1}.}
 
 Suppose $\Lambda = \Lambda(\LL, W)$. 
We claim that for $r \geq 1$,
\eq{eq: preceding}{
\widehat{\mathbf{1}_{B_r}} (\Lambda) \ll \vol^{(d)}(B_r)
\alpha(\underline \LL),
}
where $\Lambda = \Lambda(\LL, W)$ and $\underline \LL = \underline
\pi(\LL)$, and where the implicit constant depends on $d, \, n$ and $W$. Indeed, we
can replace $W$ with a larger convex set containing it, so that
$\widehat{\mathbf{1}_{B_r}}(\Lambda)$ is bounded from above by $\# \, (K
\cap \LL)$, where $K \df B_r \times W$. It
is known (see \cite[Chap. 2 \S 9.4]{GL} or \cite[Prop. 2.9]{Widmer})
that for any dimension $n$, 
for any bounded convex set $K'$ and any lattice $\LL' 
\subset \R^n$, if $K' \cap \LL' $ is not contained in a 
proper affine subspace of $\R^n$, then  
$$
\# \, \left( K' \cap \LL' \right)  \leq n! \, \frac{ \vol(K')}{\covol(\LL')} + n.
$$
For any $\LL$ we let $x_0$ be a translation
vector such that $\LL+x_0 = \underline{\LL}$, set
$$
V \df \spa (\underline \LL \cap (K+x_0)), \ \ \ell \df \dim V, \ \
\LL' \df \underline \LL \cap V, \ \ K' \df V \cap (K+x_0),
$$
and apply this estimate in $V \cong \R^\ell$ with $\ell \leq n$. 
For $r \geq 1$ we have 
$\vol^{(\ell)} (K') \ll r^d$ and $\covol(\LL') \gg
\lambda_1 (\underline \LL)\cdots \lambda_\ell (\underline \LL).$ Thus
$$
\# \, ( K \cap \LL) = \# \, (K' \cap \LL' )  \ll \ell ! \frac{r^d }{\lambda_1 (\LL)
\cdots \lambda_\ell(\LL)} + \ell \ll \vol^{(d)}(B_r)
\alpha(\underline \LL),
$$
establishing \equ{eq: preceding} and proving the claim. Therefore, using
Proposition \ref{prop: alpha in Lp} and the
dominated convergence theorem, we are justified in taking a limit $r \to \infty$ inside the
integral \equ{eq: inside the integral}, finding that  
 $c_{\mu,1} = D(\mu)$. Combining this with \equ{eq: density exists}
 gives \equ{eq: cmu1}. See \cite[Proof of Thm.\,1.5]{MS} for a different
 proof of \equ{eq: cmu1}. 

 Now to prove
 \equ{eq: cmuhigher}, let $Q_r$ and $Q_r^p$
 denote the unit cube of sidelength $r$ in $\R^{d}$ and $\R^{dp}$ respectively, let
 $\mathbf{1}_{Q_r}$ and $\mathbf{1}_{Q_r^p}$ be the indicator
 functions, and define
 ${}^p\widehat{\mathbf{1}_{Q^p_r}}$ via \equ{eq: def fphat}. Then
 we have 
$${}^p\widehat{\mathbf{1}_{Q^p_r}}(\Lambda)  = \# \, \bigtimes^p
\left(Q_r\cap  \Lambda\right); $$
that is, the number of $p$-tuples of elements of $\Lambda$ in the $p$-fold
Cartesian product $Q^p_r .$ 
This implies that for $\mu$-a.e.\,$\Lambda$, 
\eq{eq: the first integral}{
  \lim_{r\to\infty} 
\frac{{}^p\widehat{\mathbf{1}_{Q^p_r}}(\Lambda)  }{r^{dp}}
= \left( \lim_{r\to \infty} \frac{\# \,( Q_r  \cap \Lambda)
  }{\vol^{(d)}(Q_r)}\right)^p= D(\Lambda)^p= c_{\mu,1}^p. 
}


By Theorem \ref{thm: weil} 
we have: 
\eq{eq: thus we are justified}{
  \begin{split}
     c_{\mu,p}   
    =& \frac{1}{r^{dp}}
\int_{\R^{dp}} \mathbf{1}_{Q^p_r} \, d\tau_{\mathrm{main}} 
= \frac{1}{r^{dp}}
\left[ 
 \int_{\R^{dp}} \mathbf{1}_{Q^p_r} \, d\tau -
 \int_{\R^{dp}} \mathbf{1}_{Q^p_r} \, d\tau_{\mathrm{rest}}\right] \\
= & \int
 \frac{{}^p\widehat{\mathbf{1}_{Q^p_r}}(\Lambda)  }{r^{dp}} \, d\mu -
 \frac{1}{r^{dp}} 
 \sum_{\ee \in \E^{\mathrm{rest}}} \int_{\R^{dp}} \mathbf{1}_{Q^p_r} \, d\tau_\ee.
\end{split}
}

Repeating the argument establishing \equ{eq: preceding}, we find
$$
{}^p\widehat{\mathbf{1}_{Q_r}}(\Lambda)  \ll \left(\vol^{(d)}(Q_r)
\right)^p \alpha(\underline \LL)^p, 
$$
and thus the integrable function $\alpha^p$ dominates the
integral in the
second line of \equ{eq: thus we are 
  justified}, independently of $r$.  Moreover, since they differ by a
constant,  $\alpha^p$ also dominates the series in the
second line of \equ{eq: thus we are 
  justified}. Using \equ{eq: the first
  integral}, the first integral gives $c_{\mu,1}^p$, and thus it
remains to show that 
\eq{eq: tends to zero}{\lim_{r \to \infty} 
\frac{1}{r^{dp}}
\int_{\R^{dp}}\mathbf{1}_{Q_r^p} \,
d\tau_\ee =0, \ \ \text{ for every } \ee \in 
\E^{\mathrm{rest}}.
}
From \equ{eq: new condition} and Corollary \ref{cor: for small p} we
have that $\tau_\ee$ is 
(up to proportionality) equal to Lebesgue measure on a
subspace $V' \subset \R^{dp}$, and we have $V' \neq \R^{dp}$ since $\ee
\in \E^{\mathrm{rest}}$. This implies \equ{eq: tends to zero}. 
\end{proof}

\begin{remark}
  One can also work in
  $\R^{np}$ rather than $\R^{dp}$, and define
  analogous normalization constants $\bar{c}_{\bar \mu, p}$ by the
  formula $\bar \tau_{\mathrm{main}} = \bar c_{\bar \mu, p}
  \,\vol^{(np)}.$ Then one can show that $\bar{c}_{\bar \mu, p}=1$ for all
  $p<q_\mu$. We will not need the values of these constants
  and leave the proofs to the interested reader. 
  \end{remark}

\subsection{More details for $p=2$}
We will need to describe the measure $\tau_{\mathrm{rest}}$
in the case $p=2$. 
\begin{prop}\name{prop: description taurest}
Let $\mu$ be an RMS measure so that \equ{eq: new condition} holds. Let
$p=2$, and let $\E^{\mathrm{rest}}, 
\, \tau_{\mathrm{rest}}$ be as in \equ{eq: def Erest}. Then there is a
partition $\E^{\mathrm{rest}} = \E_1^{\mathrm{rest}} \sqcup
\E_2^{\mathrm{rest}}$, and constants $\{a_\ee : \ee \in
\E_2^{\mathrm{rest}}\}, \, \{b_\ee: \ee \in \E_1^{\mathrm{rest}}\}, \,
\{c_\ee: \ee \in \E^{\mathrm{rest}}\},$ such that the following hold.
\begin{enumerate}
  \item
  For all $f \in C_c(\R^{2d})$, we have
  \eq{eq: 1}{
\int_{\R^{2d}} f \, d\tau_{\mathrm{rest}} = \sum_{\ee \in
  \E_1^{\mathrm{rest}}} c_\ee \int_{\R^d} f(x, b_\ee x) \,
d\vol^{(d)}(x) + \sum_{\ee \in \E_2^{\mathrm{rest}}} c_\ee \int_{\R^d} f(a_\ee x, x ) \,
d\vol^{(d)}(x).
}
\item $c_\ee >0$ for all $\ee \in \E^{\mathrm{rest}}$ and
$\sum_{\ee \in \E^{\mathrm{rest}}} c_\ee < \infty.$
\item
  $|a_\ee |\leq 1$ for all $\ee \in \E_2^{\mathrm{rest}}$ and $|
  b_\ee | \leq 1
  $ for all $\ee \in \E_1^{\mathrm{rest}}.$ 

  \end{enumerate}
  
\end{prop}

\begin{proof}
  Lemma \ref{lem: for small p} is applicable in view of \equ{eq: new
    condition}; indeed, when $\mathbf{G}=\SL_k$, we have $p = 2 \leq
  d$, and when $\mathbf{G} = \Symp_{2k}$ and $\mu$ is affine, we have
  $\rank (v_1, v_2) \leq 1$. 
Therefore, for each $\ee \in
\E^{\mathrm{rest}}$, there is an annihilator subspace $V_\ee
\varsubsetneq \R^{dp}$ such
that $\tau_\ee$ is proportional to Lebesgue measure on
$V_\ee$. Repeating the argument of \S\ref{subsec: derivation} we can
see that $\tau_\ee$ is not the Dirac mass at the origin. In other
words $V_\ee$ has positive dimension. Since $p=2$, this
means we can find $\alpha, \beta$, not both zero, such that $V_\ee =
\mathcal{Z}(\alpha, \beta)$. We can rescale so that $\max
(|\alpha|, |\beta|)=\max(\alpha, \beta)=1$ and we define
$$
\E_1^{\mathrm{rest}} \df \{\ee \in \E^{\mathrm{rest}}: \beta=1\}, \ \
\E_2^{\mathrm{rest}} \df \E^{\mathrm{rest}}\sm \E_1^{\mathrm{rest}}.
$$
Then if we set $b_\ee = -\alpha$ for $\ee \in \E_1^{\mathrm{rest}}$
and $a_\ee = -\beta $ for $\ee \in \E_2^{\mathrm{rest}}$, then the
bounds in (3) hold and we have 
$$
V_\ee = \left\{ \begin{matrix} \{(x, b_\ee x): x \in \R^d\} &
  \  \text{ for } \ee
    \in \E_1^{\mathrm{rest}} \\ \{( a_\ee x,x ) : x \in \R^d \} &\ \ \text{
      for } \ee
    \in \E_2^{\mathrm{rest}}. \end{matrix} \right.  
$$
We now define $c_\ee$ by the formula
$$
\forall f \in C_c(\R^{2d}), \ \ \int_{\R^{2d}} f \, d\tau_\ee  =
\left\{ \begin{matrix} c_\ee \int_{\R^d} f(a_\ee x, x) \,
    d\vol^{(d)}(x) & \ \text{ for } \ee
    \in \E_1^{\mathrm{rest}} \\ c_\ee \int_{\R^d} f(x, b_\ee x) \,
    d\vol^{(d)}(x) & \ \ \text{ for } \ee
    \in \E_2^{\mathrm{rest}}. \end{matrix} \right.  
$$
Then clearly \equ{eq: 1} holds, and $c_\ee>0$ for all $\ee \in
\E^{\mathrm{rest}}$.

It remains to show $\sum c_\ee <\infty$. 
Let $\mathbf{1}_B$ be the indicator of a ball in $\R^{2d}$ centered at
the origin. Then there is a positive number $\lambda$ which bounds
from below all the numbers 
$$
\left \{ \int_{\R^d} \mathbf{1}_B(ax, x) \, d\vol^{(d)}(x) : |a|\leq 1
\right\} \bigcup \left \{ \int_{\R^d} \mathbf{1}_B(x, bx) \, d\vol^{(d)}(x) : |b|\leq 1
\right\}.
$$
Since 
$\tau_{\mathrm{rest}}$ is a locally finite measure, we have
$\int_{\R^{2d}} \mathbf{1}_B \, d\tau_{\mathrm{rest}}<\infty$. But
\equ{eq: 1} implies that $\lambda \sum_{\ee \in \E^{\mathrm{rest}}}
c_\ee \leq \int_{\R^{2d}} \mathbf{1}_B \, d\tau_{\mathrm{rest}}. $
Therefore $\sum_{\ee \in \E^{\mathrm{rest}}}
c_\ee < \infty$. 
\end{proof}

\begin{proof}[Proof of Theorem \ref{thm: Rogers bound RMS}]
  Given $f: \R^d \to [0,1]$ as in Theorem \ref{thm: Rogers bound RMS}, define
$$
\varphi: \R^{2d}\to [0,1] \ \ \ \text{ by } \varphi(x, y) \df f(x)f(y).
$$
Clearly $\left(\int_{\R^d} f \, d\vol^{(d)} \right)^2 = \int_{\R^{2d}} \varphi\,
d\vol^{(2d)}$, and 
it follows easily from \equ{eq: Siegel Veech transform} and \equ{eq:
  def fphat} that 
\eq{eq: follows easily from}{
{}^2\widehat{\varphi}(\Lambda) = \hat{f}(\Lambda)^2.
}
  Using \equ{eq: follows easily from}, Theorem \ref{thm: 
    weil} with $p=2$, \equ{eq: Siegel
    summation MS}, and \equ{eq: cmuhigher} we have that 
  \[
    \begin{split}
  & \int_{\Cl(\R^d)}  \left| \hat f(\Lambda) - \int_{\Cl(\R^d)}
  \hat f \, d\mu \right|^2  \,
d\mu(\Lambda)  \\= & \int_{\Cl(\R^d)}
\hat f^2 \, d\mu  - \left[\int_{\Cl(\R^d)} \hat f(\Lambda) \, d\mu \right]^2
 =    \int_{\R^{2d}} \varphi\, d\tau-
\left[c_{\mu,1} \int_{\R^d} f\, d\vol^{(d)} 
\right]^2  \\= & c_{\mu, 2} \int_{\R^{2d}} \varphi\, d\vol^{(2d)}
+\int_{\R^{2d}} \varphi\, d\tau_{\mathrm{rest}} - 
c_{\mu,1}^2 \left[\int_{\R^d} f\, d\vol^{(d)} 
\right]^2 =
\int_{\R^{2d}} \varphi \, d\tau_{\mathrm{rest}}.
\end{split}
\]
It remains to show that
\eq{eq: where now}{
  \int_{\R^{2d}} \varphi \, d\tau_{\mathrm{rest}} \ll
\int_{\R^d} f \, d\vol^{(d)},}
where the implicit constant is allowed to depend on $\mu$. And indeed,
by Proposition \ref{prop: description taurest}, we have
\[\begin{split}
    & \int_{\R^{2d}} \varphi \, d\tau_{\mathrm{rest}} \\
    \stackrel{\equ{eq: 1}}{=}
    & \sum_{\ee \in
  \E_1^{\mathrm{rest}}} c_\ee \int_{\R^d} f(a_\ee x) f( x) \, d\vol(x)
+ \sum_{\ee \in
  \E_2^{\mathrm{rest}}} c_\ee \int_{\R^d} f(x) f( b_\ee x) \,
d\vol(x) \\ 
\stackrel{f \leq 1}{\leq}
& \sum_{\ee \in
  \E_1^{\mathrm{rest}}} c_\ee \int_{\R^d} f( x) \, d\vol(x)
+ \sum_{\ee \in
  \E_2^{\mathrm{rest}}} c_\ee \int_{\R^d} f(x) \,
d\vol(x) \\ = & \left( \sum_{\ee \in \E^{\mathrm{rest}}} c_\ee
\right) \, \int_{\R^d} f \, d\vol^{(d)}.
  \end{split}
  \]

\end{proof}

\section{From bounds on correlations to a.e.\  effective
  counting}\name{section: Schmidt counting}
In this section we present two results which we will use for
counting. The first is due to Schmidt \cite{Schmidt_metrical} but we recast it in
a slightly more general form (see also
\cite[Thm. 2.9]{Kleinbock_Skenderi}).
To simplify notation, for
measurable $S \subset \R^n$, we will write $V_S \df \vol^{(n)}(S)$. 

  \begin{thm}\name{thm:schmidt}
Let $n \in \N$ and let $\mu$ be a probability measure on
$\Cl(\R^n)$. Let $\kappa \in [1,2)$, let  $\Phi = \{B_\alpha : 
\alpha \in \R_+\}$ be an unbounded ordered 
family of Borel subsets of $\R^n$, and let $\psi: \R_+ \to
\R_+$. Suppose the following 
hypotheses are satisfied:
\begin{itemize}
\item[(a)]
The measure $\mu$ is supported on discrete sets, and for each $f\in
L^1(\bR^n,\vol)$, a Siegel-Veech transform as in \equ{eq: Siegel Veech transform}
satisfies that 
$\widehat{f}\in 
L^2(\mu).$ Furthermore, there are positive $a,b$ 
such that for any  function $f: \R^n \to [0,1]$, $f \in L^1(\R^n,
\vol)$, we have 
\begin{equation}\label{eq:expec}
\int \widehat{f} \, d\mu =a\int_{\R^n} f \, d\vol 
\end{equation}
and
\begin{equation}\label{eq:varbound}
\Var_{\mu}(\widehat{f})\df \int \left| \widehat{f}-\int \widehat{f} \,
  d\mu \right|^2\, d\mu \leq b \left( \int_{\R^n} f \, d\vol \right)^\kappa.
\end{equation}
\item[(b)]
The function $\psi$ is non-decreasing, and satisfies
$\int_0^\infty\frac{1}{
  \psi(x)} \, dx  < \infty.$
\end{itemize}
Then for $\mu$-a.e.\ $\Lambda$, for every $S\in \Phi$
\eq{eq: for every}{
\# \, (S\cap
\Lambda) =aV_S+O\left(V_S^{\frac{\kappa}{2}}
  \log(V_S)\, \psi(\log{V_S})^{\frac12}  \right) \ \ \text{ as } V_S \to \infty. }

    \end{thm}

Note that we allow defining $\widehat{f}$ as in either one of the
linear or affine cases of \equ{eq: Siegel Veech transform}, as long as
the conditions in (a) are satisfied. For definiteness we will use the
affine case, namely $\widehat{f} = \sum_{v 
    \in \Lambda} f(v)$, so that $\widehat{\mathbf{1}_S}(\Lambda) =|S
  \cap \Lambda|$ for any subset $S \subset \R^n$ with indicator
  function $\mathbf{1}_S$. In the linear case
  we may have  $\widehat{\mathbf{1}_S}(\Lambda) =|S 
  \cap \Lambda|-1$ or $\widehat{\mathbf{1}_S}(\Lambda) =|S 
  \cap \Lambda|$ (depending on whether or not  $S$ contains 0), and
  the reader will have no difficulty adjusting 
  the proof in this case.

\begin{proof}[Proof of Theorem \ref{thm: Schmidt analogue 1} assuming
  Theorem \ref{thm:schmidt}]
Taking $\kappa=1$ and $\psi(t)=t^{1+\vre}$, \equ{eq: for every} becomes
$$\# \, (S \cap \Lambda) =aV_S+O\left(V_S^{\frac{1}{2}} \, (\log V_S)^{\frac{3}{2}+\vre}
  \right) \ \ \text{ as } V_S \to \infty,
$$
which implies \equ{eq: density with rate}. The hypotheses of Theorem
\ref{thm:schmidt} hold in 
the higher rank case by \equ{eq: Siegel summation MS} and Theorem
\ref{thm: Rogers bound RMS}. 
\end{proof}

Before giving the proof of Theorem \ref{thm:schmidt} we will state the
following more general result. 

\begin{thm}\name{thm:schmidt2}
Let $d,m,n \in \N$ with $n= d+m$, let $\mu$ be a probability measure on
$\Cl(\R^n)$, let $\lambda \in [0,1), \, \kappa
\in [1,2),$ let $\psi: \R_+ \to \R_+$, let $\Phi = \{B_\alpha: \alpha
\in \R_+ \}$ be an unbounded ordered family of
Borel subsets of $\R^d$, and let $\{W_\alpha: \alpha \in
\R_+\}$ be a collection of subsets 
of $\R^m$. Suppose that (a) and (b) of Theorem \ref{thm:schmidt} are
satisfied, and in addition: 
\begin{itemize}
\item[(c)]
For any $N\in\N$ there is
$\alpha$ such that $\vol^{(d)}(B_\alpha)=N$.
\item[(d)]
 Each $W_\alpha$
can be partitioned as a disjoint union 
$W_\alpha=\bigsqcup_{\ell = 1}^{L_\alpha} C_\alpha(\ell)$, where
$L_\alpha\asymp \left(\vol^{(d)}\left(B_\alpha\right)\right)^{\lambda}$, and where 
$w_\alpha \df
\vol^{(m)}(C_\alpha(\ell)) $ is
the same for $\ell =1, \ldots, L_\alpha,$ and is of order $\asymp
\left(\vol^{(d)}(B_\alpha)\right)^{-\lambda}$.

\end{itemize}
Denote 
$\bar \Phi \df
\{B_\alpha \times W_\alpha: \alpha \in \R_+\}$ and  for $S \in \bar \Phi$,
denote $V_S \df
\vol^{(n)}(S)$.
Then for $\mu$-a.e.\ $\Lambda$, for every $S\in \bar \Phi$
\eq{eq: for every2}{
\#\, (S\cap
\Lambda)=aV_S+O\left(V_S^{\frac{\kappa(1-\lambda)}{2}+\lambda}
  \, \log(V_S) \, \psi(\log{V_S})^{\frac12}  \right), \ \text{ as } V_S \to \infty. }

\end{thm}

Note that for $\kappa=1$ and
$\psi(t)=t^{1+\vre}$, \equ{eq: for every2} becomes
\eq{eq: 8.2 becomes}{\#\, \left(S \cap \Lambda \right)
=aV_S+O\left(V_S^{\frac{1+\lambda}{2}}\log{V_S}^{\frac{3}{2}+\vre} \right).
}

Theorems \ref{thm:schmidt} and \ref{thm:schmidt2} both follow from
ideas developed by Schmidt in \cite{Schmidt_metrical}. We begin with
Theorem \ref{thm:schmidt}, for which we need the following Lemmas. 

%
By the definition of an unbounded ordered family, we can assume that
for each $V>0$ there is $\Omega \in \Phi$ such that $\vol(\Omega) = V$. 
For each $N\in\N$, let $S_N \in \Phi$ with $\vol(S
_N)=N$ and let $\rho_N \df \mathbf{1}_{S_{N}}$ denote its indicator
function. Given two integers $N_1<N_2$, let  
\[_{N_1}\rho_{N_2} \df \rho_{N_2} -\rho_{N_1}.
\]
Since the $S_N$ are nested, we have $_{N_1}\rho_{N_2} = \mathbf{1}_{S_{N_2}\sm
    S_{N_1}}. $ 

\begin{lem}[cf. \cite{Schmidt_metrical}, Lemma 2] \name{lem: Schmidt2}
     Let $T \in \N$ and let $K_T$ be the set of all pairs of
     integers $N_1$, $N_2$ satisfying $0\leq N_1<N_2\leq 2^T$,
     $N_1=u2^t$, $N_2=(u+1)2^t$, for integers $u$ and $t\geq0$.  Then
     there exists $c>0$ such that
\begin{equation}\label{eq:lem2}
   \sum_{(N_1,N_2)\in K_T} \Var_\mu(\widehat{_{N_1}\rho_{N_2}})\leq c (T+1)2^{\kappa T}. 
\end{equation}
 \end{lem} 
\begin{proof}
Indeed, \eqref{eq:varbound} yields
$\Var_\mu(\widehat{_{N_1}\rho_{N_2}})\leq b(N_2-N_1)^\kappa$. Each
value of $N_2-N_1=2^t$ for $0\leq t\leq T$ occurs $2^{T-t}$ times,
hence 
\[
\sum_{(N_1,N_2)\in K_T} (N_2-N_1)^\kappa=\sum_{0\leq t\leq T}
2^{T+(\kappa-1)t}\leq (T+1) 2^{\kappa T}. 
\]  
\end{proof}
\begin{lem}[cf. \cite{Schmidt_metrical}, Lemma 3] \name{lem: Schmidt3}
      For all $T \in \N$ there exists a subset
      $\on{Bad}_T \subset \supp \, \mu$ of measure 
\eq{eq: bound on Bad}{
\mu(\on{Bad}_T)\leq c\psi(T\log 2-1)^{-1}
}
such that
\begin{equation}\label{eq:lem3}
(\widehat{\rho_N}(\Lambda)-aN)^2\leq T(T+1)2^{\kappa T}\psi(T\log 2-1)
\end{equation}
for every $N\leq 2^T$ and all $\Lambda\in\supp\, \mu \sm \on{Bad}_T$.
 \end{lem}
\begin{proof}
  Let $\on{Bad}_T$ be the set of $\Lambda\in\supp{\mu}$ for
which  it is not true that
\begin{equation}\label{eq:prooflem3}
\sum_{(N_1,N_2)\in K_T}
\left(\widehat{_{N_1}\rho_{N_2}}(\Lambda)-a(N_2-N_1)\right)^2 \leq
(T+1)2^{\kappa T} \psi(T\log 2-1).
\end{equation}
Then the bound \equ{eq: bound on Bad} follows from Lemma \ref{lem:
  Schmidt2} by Markov's
inequality. Assume $N\leq 2^T$ and $\Lambda\in\supp\, \mu \sm 
\on{Bad}_T$. The interval $[0,N)$ can be expressed as a union of 
intervals of the type $[N_1,N_2)$, where $(N_1,N_2)\in \mathcal{I}_N
\subset K_T$ and $|\mathcal{I}_N| \leq T$. Therefore, 
$\widehat{\rho_N(\Lambda)}-aN=\sum\left(\widehat{_{N_1}\rho_{N_2}}(\Lambda)-a(N_2-N_1)\right)$,
where the sum is over $(N_1,N_2)\in \mathcal{I}_N$. Applying the 
Cauchy-Schwarz inequality to the square of this sum together with the bound
from \eqref{eq:prooflem3} we obtain \eqref{eq:lem3}. 
\end{proof}

\begin{proof}[Proof of Theorem \ref{thm:schmidt}]
Let $\on{Bad}_T$ be the sets from Lemma \ref{lem: Schmidt3}. Since
$\psi^{-1}$ is integrable and monotone, we find by
Borel-Cantelli and \equ{eq: bound on Bad} that for $\mu$-a.e.\
$\Lambda$ there is $T_{\Lambda}$ 
such that for any $T\geq T_{\Lambda}$, $\Lambda\not\in \on{Bad}_T$. Assume
now $N\geq N_{\Lambda}=2^{T_{\Lambda}}$ and let $T$ be the unique
integer for which $2^{T-1}\leq N< 2^T$. By Lemma \ref{lem: Schmidt3},
\begin{equation}\label{eq:schmidtforinteger}
(\widehat{\rho_N}(\Lambda)-aN)^2\leq T(T+1)2^{\kappa T}\psi(T\log 2-1)
 = O \left(N^{\kappa}(\log{N})^2\psi(\log{N}) \right).        
\end{equation}
Given arbitrary $S\in\Phi$, let $N$ be such that $N\leq V_S< N+1$, and
let  $S_N, S_{N+1} \in \Phi$ with
$S_{N}\subset S \subset S_{N+1}$ 
and $\vol(S_N)=N, \ \vol(S_{N+1})=N+1. $
Then
\eq{eq: the LHS}{
\# \, (S_N\cap \Lambda )- a(N+1)\leq  \#\, (S\cap \Lambda) - aV_S \leq
\# \, (S_{N+1}\cap \Lambda ) - aN.
}
From \eqref{eq:schmidtforinteger}, the LHS of \equ{eq: the LHS} is
$O\left(N^{\frac{\kappa}{2}}\, \log{N}\, \psi(\log{N})^{\frac{1}{2}} \right)$ and the RHS is
$O\left((N+1)^{\frac{\kappa}{2}}\, \log{(N+1)}\, \psi(\log{N+1})^{\frac{1}{2}}\right)$, and
these quantities are of the same order
$O\left(V_S^{\frac{\kappa}{2}}\, \log(V_S)\, \psi(\log{V_S})^{\frac{1}{2}}\right)$.  A
similar upper bound for 
$aV_S   - \# \, (S\cap \Lambda)$ is proved analogously. 
\end{proof}

We turn to the proof of Theorem \ref{thm:schmidt2}. Note that the
collection $\bar \Phi$ is not ordered; nevertheless one 
can apply similar arguments to each $\ell$
separately, before applying Borel-Cantelli. We turn to the details.

\begin{proof}[Proof of Theorem \ref{thm:schmidt2}]
Given $N$, using  assumption (c), for each $N$ there is  $\alpha = \alpha(N)$ so that
$\vol^{(d)}\left(B_{\alpha} \right)=N$. It follows that
$\vol^{(n)}(B_\alpha \times W_\alpha)=NL_\alpha w_\alpha\asymp N$. We let 
$\rho_N^{\ell}$ be the characteristic function of $B_\alpha \times
C_\alpha(\ell)$, which is of volume $Nw_\alpha\asymp N^{1-\lambda}$.  
We will take
$_{N_1}\rho^\ell_{N_2}$ to be the characteristic
function of $\left(B_{\alpha(N_1)}\setminus B_{\alpha(N_2)} \right) \times
C_{\alpha(N)}(\ell)$. Note that the dependence of the function
$_{N_1}\rho^\ell_{N_2}$ on $N$ is suppressed from the notation.

The argument proving Lemma \ref{lem: Schmidt2} therefore yields 
\eqref{eq:lem2}, with $\kappa$
replaced by $\kappa'\df \kappa(1-\lambda)$, i.e.,  
\begin{equation}\label{eq:lem2new}
   \sum_\ell\sum_{(N_1,N_2)\in K_T}
   \Var_\mu(\widehat{_{N_1}\rho^\ell_{N_2}}) \leq cL_\alpha (T+1)2^{\kappa' T}. 
\end{equation}

For $S = B_{\alpha(N)} \times W_{\alpha(N)}$,
$N \leq 2^T$,  by the definition of
$\widehat{_{N_1}\rho^\ell_{N_2}}(\Lambda)$ and the 
Cauchy-Schwarz inequality, we have  
\[ \begin{split}
 (\# \, (S\cap\Lambda)-aV_S)^2  = & \left( \sum_\ell \# \, \left((B_\alpha \times
   C_\alpha(\ell) )\cap \Lambda \right) - a \vol^{(n)}(B_\alpha \times
   C_\alpha(\ell)) \right)^2
 \\ = &
\left(\sum_{\ell}\left(\widehat{\rho_N^{\ell}}(\Lambda)-aNw_\alpha\right)\right)^2 \\
= &
\left(\sum_{\ell}\sum_{(N_1,N_2) \in \mathcal{I}_N}\left(\widehat{_{N_1}\rho^\ell_{N_2}}
    (\Lambda)-a(N_2-N_1)w_\alpha\right)\right)^2
\\
\leq &   
T \, L_\alpha\, \sum_{\ell}\sum_{(N_1,N_2)\in
  K_T}\left(\widehat{_{N_1}\rho^\ell_{N_2}}(\Lambda)-a(N_2-N_1)w_\alpha\right)^2. 
\end{split}
\]

As in the proof of Lemma \ref{lem: Schmidt3}, we denote by
$\on{Bad}_T$ the  points $\Lambda$ 
not satisfying the bound 
\begin{equation*}\label{eq:lem3proof2}
\sum_{\ell}\sum_{(N_1,N_2)\in
  K_T}\left(\widehat{_{N_1}\rho^\ell_{N_2}}(\Lambda)
  -a(N_2-N_1)w_\alpha \right)^2\leq L_\alpha(T+1)2^{\kappa' T}\psi(T\log 2-1).
\end{equation*}
Then applying \eqref{eq:lem2new} we get $\mu(\on{Bad}_T)\leq
c'\psi(T\log 2-1)^{-1}$, so that by Borel-Cantelli, a.e.\,$\Lambda$
belongs to at most finitely many sets $\on{Bad}_T$. Also for
$\Lambda\not\in \on{Bad}_T$,  we have
\[
\left| \# \, (S\cap\Lambda)-aV_S \right|^2\leq L_\alpha^2 T(T+1)2^{\kappa' T}\psi(T\log 2-1),
\]
which replaces \eqref{eq:lem3}, and we proceed as before.
\end{proof}

\section{Counting patches \`a la Schmidt}\name{sec: counting}
 
In this section we prove Theorem 
\ref{thm: Schmidt in patches}. We recall some notation and terminology
from the introduction and the statement of the theorem. For a
cut-and-project set $\Lambda \subset \R^d$, $x \in \R^d$ and $R>0$, 
$\PP_{\Lambda, R}(x) = B(0, R) \cap (\Lambda -x)$
is called the $R$-patch of $\Lambda$
at $x$, and 
$$
D(\Lambda, \PP_0) = \lim_{T \to \infty} \frac{\# \{x \in \Lambda \cap
  B(0,T) : \PP_{\Lambda, R}(x)  = \PP_0 \}}{\vol(B(0,T))}
$$ is called the {\em frequency} of $\PP_0$. Suppose $\Lambda$ arises
from a cut-and-project construction with associated dimensions $n =
d+m$ and window $W \subset \R^m$, and is chosen according to an RMS
measure $\mu$ of higher rank. The {\em upper box
  dimension} of $W_0 \subset \R^m$ is 
\eq{eq: upper box dim def}{
\dim_B(W_0) \df \limsup_{r \to 0} \frac{\log N(W_0, r)}{-\log r},
}
where $N(W_0, r)$ is the minimal number of balls of radius $r$ needed to
cover $W_0$. Set
\eq{eq: def lambda zero}{
\lambda_0 \df \frac{m}{m+2\delta}
  }
where  $\delta = m - \dim_B(\partial \, W) >0.$
Our goal is to show that for any $\lambda \in( \lambda_0,1)$, any unbounded
ordered family $\{B_\alpha: \alpha 
\in \R\}$, for $\mu$-a.e.\  $\Lambda$, for any patch $\PP_0 =
\PP_{\Lambda, R}(x_0)$, 
\eq{eq: slightly stronger}{
  \begin{split}
    & \#\left\{x\in B_\alpha \cap \Lambda: \PP_{\Lambda, R}(x)  =  \PP_0
 \right\} \\= &
D(\Lambda, \PP_0)\, \vol(B_\alpha) +O \left(
  \vol(B_\alpha)^{\frac{1+\lambda}{2}}\right) \ \
\text{ as } \vol(B_\alpha) \to \infty,
\end{split}
}
 where the implicit constant depends on $\vre, W, \Lambda, 
\PP_0$. 
 Note that \equ{eq: slightly stronger} implies
 \equ{eq: slightly weaker}. 

 The strategy we will use is similar to that of 
 \cite[Proof of Cor. 4.1]{HKW}. 
\begin{proof}[Proof of Theorem \ref{thm: Schmidt in patches}]
 For every $K\in\N$ and $\ell\in\Z^m$ define the box 
 \[
 Q_K(\ell)=\left[\frac{\ell_1}{K},\frac{\ell_1+1}{K}\right)\times\cdots
 \times\left[\frac{\ell_m}{K},\frac{\ell_m+1}{K}\right).  
\]
It is well-known (see e.g.\,\cite[Chap. 5]{Mattila}) that in \equ{eq: upper box dim def}, we 
are free to replace  $N(W, r)$ with the minimal number
of cubes $Q_K(\ell)$ needed to  cover $W$, where $K  = \left \lfloor 
\frac{1}{r} \right \rfloor.$ 
 We consider
 cut-and-project sets of the form $\Lambda=\Lambda(W,\LL)$, with $\LL
 \in \ALN_n$. Here $W \subset \R^m$ is fixed and satisfies $\dim_B(W)
 < m$, and $\Lambda$ is chosen at random, according to a homogeneous measure
 $\bar \mu$ on $\ALN_n$. Let $\Delta$ be an $R$-patch equivalence
 class in 
 $\Lambda$, that is
 $$\Delta = \{x \in \Lambda:
 \PP_{\Lambda, R} (x) = \PP_0\}$$
 for some $R>0$ and some $\PP_0 =
 \PP_{\Lambda, R}(x_0)$. By a well-known observation (see 
\cite[Cor. 7.3]{BaakeGrimm}), $\Delta$ is itself a cut-and-project
set, and in fact arises from the same lattice via a smaller window,
i.e., there is $W_\Delta \subset W$ such
that 
$$\Delta = \Lambda(W_\Delta,\LL).$$
In particular, for irreducible cut-and-project sets (which is a
property satisfied by 
$\bar \mu$-a.e.\  $\LL$), we have 
 \eq{eq: their densities}{
 D(\Lambda, \PP_0)=D(\Delta) =\frac{\vol(W_\Delta)}{\vol(W)} \,
 D(\Lambda). 
}
In addition, it is
shown in \cite[\S2]{Koivusalo_Walton}
that $W_\Delta$ is the intersection
 of finitely many translations of $W$ and its complement. Since  
 \[
 \partial W_\Delta\subset F+\partial W,
\]
for some finite $F \subset \R^m$, 
 we deduce that the upper box dimension of $\partial W_\Delta$ is
 bounded from above by that of $\partial W$.  
 
 Let $\lambda \in (\lambda_0,1)$, and let $\eta >0$ be small enough so that
 \eq{eq: choice of eta}{
 \max \left( \frac{1+  \lambda_0}{2}+\eta,  1- \frac{\lambda_0(\delta
     -\eta)}{m} \right) < \frac{1+\lambda}{2}.}
 Such $\eta$ exists in light of \equ{eq: def
   lambda zero}.
 Given $\alpha$, we let $K_\alpha\in\N$ so that
 $\vol(B_\alpha)^{\lambda_0}\asymp K_\alpha^m$. Define
 \[
 A_{\alpha}^{(1)} \df \bigcup_{Q_{K_\alpha}(\ell)\subset W_\Delta}Q_{K_\alpha}(\ell),
 \quad
 A_{\alpha}^{(2)}\df\bigcup_{Q_{K_\alpha}(\ell)\cap W_\Delta\neq\emptyset}Q_{K_\alpha}(\ell), 
 \]
 and let $\LL \in \supp \, \bar \mu$ satisfy {\bf (D)} and {\bf (I)}.
 Since  $A_{\alpha}^{(1)} \subset W_\Delta \subset A_{\alpha}^{(2)}$, 
 the associated cut-and-project sets 
 \[
 \Lambda^{(i)}_\alpha\df \Lambda \left(A_{\alpha}^{(i)}, \LL
 \right) \ \ (i=1,2)
 \]
 satisfy that for all $\alpha$,  
 \[
 \# \, \left(\Lambda^{(1)}_\alpha \cap B_\alpha \right) \leq \# \, \left( \Delta \cap
 B_\alpha \right) \leq \# \, \left (  \Lambda^{(2)}_\alpha \cap B_\alpha \right)
\]
and
$$
D\left(\Lambda^{(1)}_\alpha \right) \leq D\left(\Delta
\right) \leq D\left(\Lambda^{(2)}_\alpha \right).
$$
 Moreover, by \equ{eq: their densities}, 
 \eq{eq: for cubes}{
   D\left(\Lambda^{(2)}_\alpha \right) - D\left(\Lambda^{(1)}_\alpha\right) =
\frac{D(\Lambda)}{\vol(W)} \, \left(\vol\left(A^{(2)}_\alpha \right) -
\vol\left(A^{(1)}_\alpha\right) \right).}
 Using the triangle inequality we have 
 \eq{eq: bound separately}{
   \begin{split}
 & \big| \# \, (\Delta\cap B_\alpha)-D(\Delta) \, \vol(B_\alpha)\big|
 \leq \max_{i=1,2} \left|\# \, \left(\Lambda_{\alpha}^{(i)}\cap
 B_\alpha \right)-D(\Delta) \, \vol(B_\alpha) \right| 
\\ \le & \max_{i=1,2} 
 \left|\# \, \left(\Lambda^{(i)}_\alpha \cap
 B_\alpha\right )-D \left(\Lambda^{(i)}_\alpha\right) \vol\left(B_\alpha \right
)\right| +
 \left(D\left (\Lambda^{(2)}_\alpha \right)  -
 D\left(\Lambda^{(1)}_\alpha \right) \right) \, \vol(B_\alpha) .
\end{split}
}
We bound separately the two summands on the RHS of \equ{eq: bound
  separately}. For the first summand we use the case \equ{eq: 8.2 
  becomes} of Theorem \ref{thm:schmidt2},
with $W_\alpha=A^{(i)}_\alpha$ and $C_\alpha(\ell) =
Q_{K_\alpha}(\ell) $. Note that assumption (d) is satisfied by our
choice of $K_\alpha$, with implicit constants depending on $\PP_0$. We
obtain, for $\bar \mu$-a.e.\  $\LL$, that 
$\Lambda^{(i)}_\alpha = \Lambda(A^{(i)}_\alpha, \LL)$ satisfies 
 \[
 \left| \# \, \left(\Lambda^{(i)}_\alpha\cap B_\alpha
 \right)-D\left (\Lambda^{(i)}_\alpha \right)\vol \left(B_\alpha \right) \right|\le
c_1\left( \vol(B_\alpha)^{\frac{1+\lambda_0}{2}} \,
  \left(\log(\vol(B_\alpha))\right)^{\frac{3}{2}+\varepsilon}\right), 
 \]
 where $c_1$, as well as the constants appearing in the following
 inequalities, depends only on $\bar \Phi=\left\{B_\alpha\times
 A^{(i)}_\alpha \right \}$ and on $\LL$. 

For the second summand, recall that $\dim_B(\partial W_\Delta) \leq
m-\delta$. This implies that the number of 
 $\ell\in\Z^n$ with
 $Q_{K_\alpha}(\ell)\cap\partial W_\Delta\neq\emptyset$ is $\ll
 K_\alpha^{m-\delta+\eta}$. Therefore
 $$
\vol\left(A^{(2)}_\alpha \sm A^{(1)}_\alpha \right)
=\sum_{Q_{K_\alpha}(\ell) \cap \partial W_\Delta \neq \varnothing}
\vol \left(
Q_{K_\alpha}(\ell) \right)\ll K_\alpha^{m-\delta+\eta} \,
K_\alpha^{-m} = K_\alpha^{-\delta+\eta}. 
 $$
 This implies via \equ{eq: for cubes} that 
 \[\begin{split}
 \left(D\left (\Lambda^{(2)}_\alpha \right)   -
 D\left(\Lambda^{(1)}_\alpha \right) \right) \, \vol(B_\alpha) =  & \frac{
  D(\Lambda) \, \vol(B_\alpha) }{\vol(W)} \left( \vol \left(A^{(2)}_\alpha\right) - \vol
 \left( A^{(1)}_\alpha \right)\right) \\
\ll & \vol(B_\alpha)  \, K_\alpha ^{-\delta + \eta}  \ll
\vol(B_\alpha)^{1-\frac{\lambda_0(\delta - \eta)}{m}}. 
\end{split}
\]
Plugging these two estimates into \equ{eq: bound separately}, and
using \equ{eq: choice of eta} and the fact that
$\left(\log(\vol(B_\alpha))\right)^{\frac{3}{2}+\varepsilon} \leq
\vol(B_\alpha)^\eta $ for large enough $\vol(B_\alpha)$, we have
that for $\bar \mu$ -a.e.\  $\LL$
$$
 \big| \# \, (\Delta\cap B_\alpha)-D(\Delta) \, \vol(B_\alpha)\big| \ll
 \vol(B_\alpha)^{\frac{1+\lambda}{2}},
 $$
 with implicit constants depending on $\eta$, $\LL$ and $\vre$. 
 This shows \equ{eq: slightly stronger} and completes the proof.
\end{proof}

\ignore{

\appendix 
\section{The group in the restriction of scalars - By Dave Morris}
\name{appendix: morris}
\combarak{Rep U Com 38 the referee
  thinks that there are notations here that are not consistent and
  results that should be strengthened.}
}

\bibliographystyle{alpha}
\bibliography{newbib}

\end{document}